\documentclass[times]{oupau}

\usepackage[T1]{fontenc}
\usepackage{tgtermes}
\usepackage{mathrsfs,array,empheq}
\graphicspath{{pic/}{../pic/}}
\usepackage[bookmarksnumbered,colorlinks,plainpages]{hyperref}

\makeatletter

\let\c@proposition\relax

\makeatother
\newtheorem{proposition}{Proposition}

\newif\ifdraft
\drafttrue
\numberwithin{figure}{section}
\numberwithin{table}{section}

\makeatletter
\g@addto@macro\ps@title{%
  \let\@oddfoot\@empty
  \let\@evenfoot\@empty
}
\makeatother

\renewcommand{\titlesize}{\fontsize{18}{22pt}\selectfont}

\begin{document}

\title{The Global Solutions to Nonlinear Thermoelastic Systems under a Null Condition}
\shorttitle{The Global Solutions to Nonlinear Thermoelastic Systems}
\renewcommand{\authsize}{\fontsize{11.5}{13.5pt}\selectfont}
\author{\mdseries Yating Xu and Yi Zhou}
\abbrevauthor{Y. Xu and Y. Zhou}
\headabbrevauthor{Xu, Y., and Zhou, Y.}

\begin{abstract}
We study the Cauchy problem for a class of homogeneous, isotropic nonlinear thermoelastic systems in three space dimensions.
Under a null condition for purely elastic part of our system which is different from those considered by Sideris \cite{Sideris1996,Sideris2000}, we establish the existence and uniqueness of global classical solutions with smallness assumptions on the initial data.
We use Klainerman--Sideris weighted energy estimates and time-weighted energy estimates for the heat equation to establish global existence.
Moreover, we improve upon the global existence result of Ponce and Racke \cite{PonceRacke1990}, who assumed that the purely elastic nonlinear terms vanish.
\end{abstract}

\maketitle
	
	\section{Introduction}
	
	Thermoelastic systems describe the interaction between elastic deformation and temperature variations: deformation affects the evolution of temperature, while temperature variations influence the mechanical response of the material through thermal stresses.
	The study of homogeneous, isotropic thermoelastic materials therefore requires consideration of both deformation and temperature perturbations.
	For the thermodynamic background of thermoelastic constitutive theory, we refer to the work of Coleman and Noll \cite{ColemanNoll1963} on elastic materials with heat conduction and to the thermoelastic system developed by Biot \cite{Biot1956}.
	For studies of nonlinear elastic waves relevant to the mechanical part, we refer to John \cite{John1984} on the formation of singularities and to the global existence results of Sideris \cite{Sideris1996,Sideris2000}.
	
	The existence, asymptotic stability, and decay of solutions have been extensively studied in linear thermoelasticity.
	Leis \cite{Leis1986} used the Helmholtz decomposition to split the displacement into divergence-free and curl-free components. 
	The divergence-free component satisfies an independent wave equation, whereas the curl-free component is coupled to the temperature variable, and the coupled part is totally dissipative.
	For the three-dimensional Cauchy problem, Ponce and Racke \cite{PonceRacke1990} and Racke \cite{Racke1990Cauchy} also analyzed the decay properties of linearized thermoelastic systems.
	For some results on initial-boundary value problems, see \cite{BiesCieslak2023,BiesEtAl2026,LebeauZuazua1999,RackeShibata1991,Slemrod1981}.
	
	For the nonlinear thermoelastic Cauchy problem on the whole real line, Hrusa and Tarabek \cite{HrusaTarabek1989} established the existence of global smooth solutions for small perturbations of constitutive models satisfying nondegeneracy conditions.
	For the three-dimensional Cauchy problem for nonlinear systems,
	Racke \cite{Racke1990Cauchy} studied a thermoelastic system with cubic nonlinear terms and one quadratic term involving temperature.
	Ponce and Racke \cite{PonceRacke1990} considered a thermoelastic system with quadratic nonlinear terms excluding purely elastic quadratic terms, and established the existence and uniqueness of global classical solutions for sufficiently small initial data.
	Further results can be found in \cite{Chirita1987EnergyPartition,Dafermos1968Linear,MunozRivera1993}.

	The study of nonlinear waves with small initial data further highlights the importance of the nonlinear structure alongside the smallness of the initial data.
	For three-dimensional systems of pure elasticity, John also studied the formation of singularities in radially symmetric elastic waves \cite{John1984} and almost global existence for elastic waves arising from small perturbations \cite{John1988AlmostGlobal}.
	On the other hand, Sideris \cite{Sideris1996,Sideris2000} and Agemi \cite{Agemi2000} established global existence under a null condition using weighted energy estimates.
	These purely elastic models contain no temperature variable, and their structural conditions cannot be directly identified with those used here. Nevertheless, they provide background for analyzing the mechanical nonlinearities in thermoelastic systems.
	
	The purpose of this paper is to establish global existence of solutions to nonlinear thermoelastic systems under the different null condition.
	In thermoelasticity, the mechanical response is characterized by a stored energy function, and temperature evolution is described by Fourier's law of heat conduction. 
	One difficulty in the proof is that heat conduction affects the displacement components differently. The Helmholtz decomposition $u=u^{cf}+u^{df}$ splits the displacement into curl-free and divergence-free components: the divergence-free component retains primarily a wave character, whereas the curl-free component is coupled to the temperature variable at the linear level.
	It is therefore necessary to control these components separately and to account for their interaction in the nonlinear estimates. We aim to use solutions of the coupled equations to transfer dissipation to the curl-free component and thereby control resonant interactions involving this component in the nonlinear terms.
	In addition, the two equations have different space-time scalings, which leads to derivative loss and makes it difficult to close the energy estimates. To address this difficulty, we include one additional spatial derivative in the energy.
	For the divergence-free displacement component, we combine vector fields with Klainerman--Sideris type weighted estimates to obtain wave control; for the curl-free component coupled to temperature, estimates must instead be derived from the coupled system itself, since the weighted estimate associated with the longitudinal wave speed $c_1$ cannot be applied directly.
	We construct coupled weighted energies for the curl-free displacement and temperature variables and combine them with inverse-derivative energies and weighted dissipation estimates to obtain the required control.
	
	We allow the higher-order energy to grow slowly in time while seeking uniform control of the lower-order energy.
	The proof first establishes coupled integral inequalities for the higher- and lower-order energies, then uses smallness of the lower-order energy to absorb nonlinear terms and derive a growth bound for the higher-order energy.
	This bound is then substituted into the lower-order estimate, where integrability of the time weight improves the a priori bound for the lower-order energy.
	A bootstrap argument is then used to establish global existence.
	
	The paper is organized as follows.
	The remainder of the introduction presents the equations of motion, basic notation, and the global existence theorem.
	Section \ref{2} collects commutation relations, the Helmholtz decomposition, null form estimates, Sobolev inequalities, integration by parts for vector fields, and Klainerman--Sideris type estimates.
	Section \ref{3} gives weighted energy estimates for the curl-free displacement and temperature components, and Section \ref{4} develops the proof of the main theorem by combining higher- and lower-order energy estimates.
	
	\subsection{The equations of motion}
	Consider a homogeneous thermoelastic material filling the whole space $\mathbb{R}^3$.
	Assume that its density in the undeformed state is unity.
	Let
	\begin{equation}\nonumber
		y(t,x):\mathbb{R}\times\mathbb{R}^3\longrightarrow\mathbb{R}^3;\qquad y(t,x)=x+u(t,x),
	\end{equation}
	be the smooth deformation of the material that evolves with time.
	It is an orientation-preserving diffeomorphism and the deformation gradient is the matrix $F=\nabla y$, whose components are
	$F_{ij}=\frac{\partial y_i}{\partial x_j}$.
	
	Let $W(F,\theta)$ denote the stored energy function of the thermoelastic material.
	And we require that the stored energy function satisfy the objectivity and isotropy conditions, namely,
	\begin{equation}\nonumber
		W(QF,\theta)=W(F,\theta), \qquad W(FQ,\theta)=W(F,\theta),
	\end{equation}
	for any $Q\in SO(3)$.
	
	Nonlinear thermoelastic systems couple the equations of elasticity with a heat equation describing temperature evolution.
	Let $\theta(t,x)$ denote the temperature field, with equilibrium temperature normalized to unity, so that $\theta(t,x)-1$ represents the temperature perturbation. 
	The heat flux is given by Fourier's law $\vec q=-k\nabla\theta$.
	
	In this paper, we consider the thermoelastic systems (see \cite{LiQin2014}),
	\begin{equation}\label{equ_or}
		\left\lbrace
		\begin{aligned}
			&\frac{\partial^2 u^i}{\partial t^2}-\sum_{j=1}^3\frac{\partial}{\partial x_j}\frac{\partial W(\nabla y,\theta)}{\partial F_{ij}}=0\\
			&\theta\frac{\partial \theta}{\partial t}-k\Delta \theta+\theta\sum_{i=1}^3 \frac{\partial^2 u^i}{\partial t\partial x_i}=0.
		\end{aligned}
		\right.
	\end{equation}
	
	To describe the variational structure of nonlinear thermoelasticity, we express $W(\nabla y,\theta)$ explicitly in terms of the displacement gradient $G=\nabla u$ and the temperature $\theta$.
	The displacement gradient matrix $G=\nabla u$ satisfies
	
	\begin{equation}\nonumber
		G=F-I, \qquad C=FF^T-I=G+G^T+GG^T.
	\end{equation}
	Let $\lambda_1,\lambda_2,\lambda_3$ be the eigenvalues of $C$. The principal invariants are
	\begin{align*}
		r_1&=\lambda_1+\lambda_2+\lambda_3=\operatorname{tr}C,\\
		r_2&=\lambda_1\lambda_2+\lambda_2\lambda_3+\lambda_1\lambda_3=\frac{(\operatorname{tr}C)^2-\operatorname{tr}C^2}{2},\\
		r_3&=\lambda_1\lambda_2\lambda_3=\det C =\frac{(\operatorname{tr}C)^3-3(\operatorname{tr}C)(\operatorname{tr}C^2)+2\operatorname{tr}C^3}{6}.
	\end{align*}
	Consequently, $W(\nabla y,\theta)=\sigma(r_1,r_2,r_3,\theta)$.
	
	In three space dimensions, the global existence of small-amplitude solutions to nonlinear hyperbolic systems hinges on the specific form of the quadratic portion of the nonlinearity.
	Such compatibility conditions are often referred to as null conditions; see \cite{Christodoulou1986,John1981,JohnKlainerman1984,Klainerman1986Null,Sideris1983,SiderisTu2001}.
	From the analytical point of view, therefore, it is enough to truncate \eqref{equ_or} to second order in $u$, the higher order corrections having no influence on the existence of small solutions.
	Accordingly, we will truncate $\sigma (r_1,r_2,r_3,\theta)$ to third order in $u$.
	
	We introduce the change of variables $\theta=\sqrt{1+2\Theta}$ and retain the quadratic terms in the expansion of the stored energy function. 
	Neglecting cubic terms in the resulting equations, and following \cite{ZhaPeng2020}, we obtain the following nonlinear thermoelastic system:
	\begin{equation}\label{equ_1}
		\left\lbrace \begin{aligned} &\frac{\partial^2 u^i}{\partial t^2}-c^2_2\Delta u^i -(c_1^2-c_2^2)\sum_{j=1}^3\frac{\partial^2 u^j}{\partial x_i\partial x_j}+ \frac{\partial \Theta}{\partial x_i}=N^i(u,\Theta)\\
			&\frac{\partial \Theta}{\partial t}-k\Delta \Theta+\sum_{i=1}^3 \frac{\partial^2 u^i}{\partial t\partial x_i}=M(u,\Theta),\\
			&t=0:\quad u=u_0,\quad \partial_t u=u_1,\quad \Theta=\Theta_0\triangleq\frac{\theta_0^2-1}{2}.
		\end{aligned} \right.
	\end{equation}
	where $c_1$ and $c_2$ are respectively the propagation speeds of longitudinal waves and transverse waves, satisfying $c_1^2>\frac{4}{3}c_2^2>0$.
	The nonlinear terms are decomposed as
	\begin{equation*}
		\begin{aligned}
			N(u,\Theta)&=N_{al}(u,u)+N_{th}(\Theta,\Theta)+N_{mi}(u,\Theta),\\
			M(u,\Theta)&=M_1(\Theta,u)+M_2(\Theta,\Theta).
		\end{aligned}
	\end{equation*}
	The elastic part is given by $N_{al}(u,v)=N_1(u,v)+N_2(u,v)$, where
	\begin{equation}\label{equ_N_al}
		\begin{aligned}
			N_1(u,v)={}&3d_1\nabla\bigl((\nabla\cdot u)(\nabla\cdot v)\bigr)+d_2\nabla\bigl((\nabla\wedge u)\cdot(\nabla\wedge v)\bigr)\\
			&-d_2\nabla\wedge\bigl((\nabla\cdot u)(\nabla\wedge v)\bigr)-d_2\nabla\wedge\bigl((\nabla\cdot v)(\nabla\wedge u)\bigr),\\
			N_2^i(u,v) ={}&\left(d_3+\frac{d_4}{2}\right)\Bigl[Q_{ij}(\partial_k u^k,v^j)+Q_{ij}(\partial_k v^k,u^j) -Q_{jk}(\partial_i u^k,v^j)-Q_{jk}(\partial_i v^k,u^j)\Bigr]\\
			&+\frac{d_5}{2}\Bigl[Q_{ij}(\partial_j u^k,v^k)+Q_{ij}(\partial_j v^k,u^k)+2Q_{jk}(\partial_j u^i,v^k)+2Q_{jk}(\partial_j v^i,u^k)\Bigr]\\
			&+\frac{d_5}{2}\Bigl[Q_{jk}(\partial_j u^k,v^i)+Q_{jk}(\partial_j v^k,u^i)\Bigr]\\
			\triangleq{}&\partial_l\left(\tilde{g}^{ijk}_{lmn}\partial_m u^j\partial_n v^k\right),
		\end{aligned}
	\end{equation}
	The null forms are defined by
	\begin{equation*}
		Q_{ij}(v,w)=\partial_i v\,\partial_j w-\partial_j v\,\partial_i w.
	\end{equation*}
	The coefficients in \eqref{equ_N_al} are expressed in terms of partial derivatives of the stored energy function $\sigma$ with respect to the invariants $r_1,r_2,r_3$. 
	These derivatives are denoted by subscripts on $\sigma$ and are all evaluated at the reference state $(r_1,r_2,r_3,\theta)=(0,0,0,1)$.
	Specifically,
	\begin{equation*}
		\begin{aligned}
			&d_1=\frac{4}{3}\sigma_{111}+2\sigma_{11},\\
			&d_2=2(\sigma_{11}-\sigma_{12}),\\
			&d_3=-2\sigma_{11}+4\sigma_{12}+\frac{4}{3}\sigma_3,\\
			&d_4=-\frac{2}{3}\sigma_3,\\
			&d_5=2(\sigma_2-\sigma_3).
		\end{aligned}
	\end{equation*}
	
	The thermal and mixed terms are given, respectively, by
	\begin{equation*}
		\begin{aligned}
			N_{th}(\Theta_1,\Theta_2)&=\nabla(\Theta_1\Theta_2),\qquad N_{mi}^i(u,\Theta)&=\partial_l\bigl(g^{ij}_{lm}\partial_m u^j\Theta\bigr).
		\end{aligned}
	\end{equation*}
	The nonlinear terms in the heat equation are
	\begin{equation*}
		\begin{aligned}
			M_1(\Theta,u)=-\Theta\nabla\cdot\partial_t u,\quad M_2(\Theta_1,\Theta_2)=-k\nabla\cdot(\nabla\Theta_1\Theta_2).
		\end{aligned}
	\end{equation*}
	
	The elastic and mixed nonlinearities can also be written in divergence form:
	\begin{equation}\nonumber
		N_{al}(u,v)^i=\partial_l\bigl(g^{ijk}_{lmn}\partial_m u^j\partial_n v^k\bigr), \qquad N_{mi}^i(u,\Theta)=\partial_l(g^{ij}_{lm}\partial_m u^j\Theta),
	\end{equation}
	where the coefficients are constants and satisfy the following symmetries under interchange of index pairs:
	\begin{equation}
		g^{ijk}_{lmn}=g^{jik}_{mln}=g^{kij}_{nlm},\qquad g^{ij}_{lm}=g^{ji}_{ml}.
	\end{equation}
	Moreover, $(g^{ijk}_{lmn})$ is a sixth-order isotropic tensor and $(g^{ij}_{lm})$ is a fourth-order isotropic tensor.
	
	We say that the nonlinear thermoelastic system \eqref{equ_1} satisfies the null condition if
	\begin{equation}\label{equ_d_2}
		d_2=2(\sigma_{11}-\sigma_{12})=0,
	\end{equation}
	the purely elastic nonlinearity is given by
	\begin{equation*}
		N_{al}^i(u,u)=3d_1\partial_i\bigl((\nabla\cdot u)^2\bigr)+\partial_l\left(\tilde{g}^{ijk}_{lmn}\partial_m u^j\partial_n u^k\right).
	\end{equation*}
	which is equivalent to
	\begin{equation}\label{equ_null_condition_2}
		\tilde{g}^{ijk}_{lmn}\omega_l\omega_m\omega_n=0, \qquad \text{for all }\omega\in S^2 \text{ and } 1\le i,j,k\le 3.
	\end{equation}
	
	Our null condition only rules out the cubic term $(\nabla \times u)^2 \nabla \cdot u$ in the stored energy function, 
	whereas the null condition considered by Sideris \cite{Sideris1996,Sideris2000} rules out the cubic terms $(\nabla \cdot u)^3$.
	Thus we can identify whether materials satisfy the different null conditions by checking the cubic term in the stored energy function directly.
	
	\subsection{Notation}
	
	Partial derivatives will be written as
	$$\partial_0=\partial_t=\frac{\partial}{\partial t} \quad\text{and}\quad \partial_i=\frac{\partial}{\partial x^i}, \qquad i=1,2,3.$$
	Denote the spatial gradient and space-time gradient by $\nabla=(\partial_1,\partial_2,\partial_3)$ and $\partial=(\partial_t,\nabla)$, respectively.
	The angular momentum operators, namely the generators of spatial rotations, are the vector fields
	$$\Omega=(\Omega_1,\Omega_2,\Omega_3)=x\wedge\nabla,$$
	$\wedge$ being the usual vector cross product. 
	Crucially, the spatial partial derivatives can be decomposed into radial and angular components
	
	\begin{equation}\nonumber
		\nabla=\frac{x}{r}\partial_r-\frac{x}{r^2}\wedge\Omega, \qquad r=|x|, \qquad \partial_r=\frac{x}{r}\cdot\nabla.
	\end{equation}
	We define the vector-valued modifications of the angular momentum operators by
	\begin{equation}\nonumber
		\widetilde{\Omega}_{\ell}=\Omega_{\ell}I+U_{\ell},
	\end{equation}
	where $U_{\ell}=\varepsilon_{\ell ij}e_i\otimes e_j$, $\{e_i\}_{i=1}^{3}$ is the standard basis of $\mathbb{R}^3$, and $\varepsilon_{\ell ij}$ is the Levi-Civita symbol.
	The scaling operator is
	\begin{equation}\nonumber
		S=t\partial_t+r\partial_r,
	\end{equation}
	and we will use
	\begin{equation}\nonumber
		\widetilde{S}=S-1.
	\end{equation}
	
	The eight vector fields will be written as
	$\tilde\Gamma=(\tilde\Gamma_0,\ldots,\tilde\Gamma_7)=(\partial,\widetilde{\Omega},\widetilde{S})$ and $\Gamma=(\Gamma_0,\ldots,\Gamma_7)=(\partial,\Omega,S)$.
	It is easy to verify that the commutator of $\tilde\Gamma$ and $\Gamma$ with $\nabla$ is in the span of $\nabla$.
	Schematically, we write
	\begin{equation}\nonumber
		[\tilde\Gamma,\nabla]=\nabla\qquad\text{and}\qquad [\Gamma,\nabla]=\nabla.
	\end{equation}
	By $\Gamma^\alpha$, $\alpha=(\alpha_1,\ldots,\alpha_k)$, we denote an ordered product of $k=|\alpha|$ vector fields $\Gamma^{\alpha_1}\cdots\Gamma^{\alpha_k}$.
	
	Denote the Riesz transformation by
	\begin{equation}\nonumber
		R_k=\frac{\partial_k}{\sqrt{-\Delta}},\qquad k=1,2,3.
	\end{equation}
	We also employ the notation $R=(R_1,R_2,R_3)$.
	We will use the Helmholtz decomposition, which projects any vector field onto curl-free and divergence-free components.
	Specifically, for any function $u\in C^{\infty}(\mathbb{R}^3;\mathbb{R}^3)$ with sufficient decay at infinity, we have
	\begin{equation}\nonumber
		u=u^{cf}+u^{df},
	\end{equation}
	with
	\begin{equation}\nonumber
		u^{cf}\triangleq\mathbb P_{cf}u=-R(R\cdot u), \qquad u^{df}\triangleq\mathbb P_{df}u=R\wedge(R\wedge u).
	\end{equation}
	Moreover, it holds that
	\begin{equation}\nonumber
		\nabla\wedge u^{cf}=0, \qquad \nabla\cdot u^{df}=0,
	\end{equation}
	and
	\begin{equation}\nonumber
		\|u\|_{L_x^2}^2=\|u^{cf}\|_{L_x^2}^2+\|u^{df}\|_{L_x^2}^2.
	\end{equation}
	The proof of the above facts can be found in \cite{ChorinMarsden1993}.
	
	Next, we denote the basic energy corresponding to the linear thermoelasticity operator \eqref{equ_1} by
	\begin{equation}\nonumber
		\begin{aligned}
			E_1(u(t),\Theta(t))&=\frac{1}{2}\int_{\mathbb{R}^3} \left(|\partial_tu|^2+c_2^2|\nabla u|^2+(c_1^2-c_2^2)(\nabla\cdot u)^2+|\Theta^2|\right)(t,x)dx+k\int_0^t\int_{\mathbb R^3}|\nabla \Theta |^2dxd\tau.
		\end{aligned}
	\end{equation}
	We define its higher- and lower-order counterparts by
	\begin{equation}\nonumber
		\begin{aligned}
			&\mathcal I_{\kappa}(t)={}\langle t\rangle^{-\delta} \sum_{|\alpha|\leq\kappa-1}\sum_{|a|\le1} E_1\big(\nabla^a\widetilde\Gamma^\alpha u,\nabla^a\Gamma^\alpha\Theta\big)(t),\\
			&\mathcal J_{\kappa-2}(t)={}\sum_{|\alpha|\leq\kappa-3}\sum_{|a|\le1} E_1\big(\nabla^a\widetilde\Gamma^\alpha u,\nabla^a\Gamma^\alpha\Theta\big)(t).
		\end{aligned}
	\end{equation}
	
	To describe the solution space, we also introduce the time-independent analog of $\tilde\Gamma$.
	Set $ \widetilde\Lambda=(\widetilde\Lambda_1,\ldots,\widetilde\Lambda_7)=(\nabla,\widetilde{\Omega},r\partial_r-1)$.
	Then the $\widetilde\Lambda$'s have the same commutation properties as the $\widetilde\Gamma$'s.
	Define the time-independent spaces
	\begin{equation}\nonumber
		H_{\widetilde{\Lambda}}^{\kappa,p}= \left\{f\in L^p(\mathbb{R}^3;\mathbb{R}^3):\widetilde\Lambda^a f\in L^p(\mathbb{R}^3;\mathbb{R}^3),|a|\leq \kappa \right\},
	\end{equation}
	with the norm
	\begin{equation}\nonumber
		\|f\|_{H_{\widetilde\Lambda}^{\kappa,p}} =\left(\sum_{|a|\leq \kappa}\|\widetilde\Lambda^a f\|_{L_x^p}^p\right)^{\frac1p}.
	\end{equation}
	Similarly, we can define the spaces
	\begin{equation}\nonumber
		H_{\Lambda}^{\kappa,p}=\left\{f\in L^p(\mathbb{R}^3;\mathbb{R}):\Lambda^a f\in L^p(\mathbb{R}^3;\mathbb{R}),|a|\leq \kappa\right\},
	\end{equation}
	with the norm
	\begin{equation}\nonumber
		\|f\|_{H_{\Lambda}^{\kappa,p}}=\left(\sum_{|a|\leq \kappa}\|\Lambda^a f\|_{L_x^p}^p\right)^{\frac1p}.
	\end{equation}
	Here, $\Lambda=(\Lambda_1,\ldots,\Lambda_7)=(\nabla,\Omega,r\partial_r)$. The vector fields in $\Lambda$ have the same commutation properties as those in $\Gamma$.
	
	The solution will be constructed in the space
	$(\dot{H}^{\kappa}_{\widetilde\Gamma}(T),\dot{H}^{\kappa}_{\Gamma}(T))$, which is the closure of the set
	$ C^{\infty}\left([0,T];C_c^{\infty}(\mathbb{R}^3;\mathbb{R}^3)\times C_c^{\infty}(\mathbb R^3;\mathbb R)\right)$
	in the norm
	$$\sup_{0\leq t\leq T}\widetilde{E}_{\kappa}(t),$$
	where
	\begin{equation}\nonumber
		\widetilde{E}_{\kappa}(t)\triangleq \sum_{|\alpha|\le \kappa-1}\Big( \sum_{|a|\le 1}E_1(\nabla^a\widetilde{\Gamma}^\alpha u,\nabla^a \Gamma^\alpha \Theta) +E_1(\nabla^{-1}\widetilde{\Gamma}^\alpha u^{cf},\nabla^{-1} \Gamma^\alpha \Theta)\Big),
	\end{equation}
	with $\nabla^{-1}=\nabla (-\Delta)^{-1}$.
	
	Using the Helmholtz decomposition, we define the following high-order weighted energy functionals:
	\begin{equation}\nonumber
		\begin{aligned}
			&X_{df,\kappa}(t)=\langle t\rangle^{-\delta}\sum_{|\alpha|\le\kappa-2} \|\langle c_2t-r\rangle\partial\nabla\widetilde\Gamma^\alpha u^{df}(t)\|_{L_x^2}^2,\\
			&X_{1,\kappa}(t)=\langle t\rangle^{-\delta}\sum_{|\alpha|\le\kappa-2} \Big[S_1(\nabla^{-1}\widetilde\Gamma^\alpha u^{cf}(t),\nabla^{-1}\Gamma^\alpha\Theta(t)) +\int_0^t S_2(\nabla^{-1}\widetilde\Gamma^\alpha u^{cf}(\tau),\nabla^{-1}\Gamma^\alpha\Theta(\tau))\,d\tau\Big],\\
			&X_{2,\kappa}(t)=\langle t\rangle^{-\delta}\sum_{|\alpha|\le\kappa-2} \Big[\langle t\rangle S_1(\widetilde\Gamma^\alpha u^{cf}(t),\Gamma^\alpha\Theta(t)) +\int_0^t\langle \tau\rangle S_2(\widetilde\Gamma^\alpha u^{cf}(\tau),\Gamma^\alpha\Theta(\tau))d\tau\Big],
		\end{aligned}
	\end{equation}
	where $\langle\cdot\rangle=(1+|\cdot|^2)^{\frac12}$ and
	\begin{equation*}
		\begin{aligned}
			&S_1(w(t),\Sigma(t))=\int_{\mathbb R^3}\left( |\partial_t w|^2 +|\nabla w|^2+|\Sigma|^2+|\partial_t \nabla w|^2+ |\Delta w|^2 + |\nabla\Sigma|^2\right) dx,\\
			&S_2(w(t),\Sigma(t))=\int_{\mathbb R^3}\left(|\nabla\Sigma|^2+|\Delta w |^2+|\partial_t \nabla w|^2+ |\Delta\Sigma|^2\right)dx.
		\end{aligned}
	\end{equation*}
	The corresponding lower-order weighted energies are defined by
	\begin{equation}\nonumber
		\begin{aligned}
			&\mathcal J_{df,\kappa-2}(t)=\sum_{|\alpha|\le\kappa-4} \|\langle c_2t-r\rangle\partial\nabla\widetilde\Gamma^\alpha u^{df}(t)\|_{L_x^2}^2,\\
			&\mathcal J_{1,\kappa-2}(t)=\sum_{|\alpha|\le\kappa-4} \Big[S_1(\nabla^{-1}\widetilde\Gamma^\alpha u^{cf}(t),\nabla^{-1}\Gamma^\alpha\Theta(t)) +\int_0^t S_2(\nabla^{-1}\widetilde\Gamma^\alpha u^{cf}(\tau),\nabla^{-1}\Gamma^\alpha\Theta(\tau))\,d\tau\Big],\\
			&\mathcal J_{2,\kappa-2}(t)=\sum_{|\alpha|\le\kappa-4} \Big[\langle t\rangle S_1(\widetilde\Gamma^\alpha u^{cf}(t),\Gamma^\alpha\Theta(t)) +\int_0^t\langle \tau\rangle S_2(\widetilde\Gamma^\alpha u^{cf}(\tau),\Gamma^\alpha\Theta(\tau))d\tau\Big],
		\end{aligned}
	\end{equation}
	We further define
	\begin{equation}\nonumber
		\begin{aligned}
			\mathcal E_{\kappa}(t)=\sup_{0\le\tau\le t}\widetilde{\mathcal I}_{\kappa}(\tau),\qquad
			\mathcal E^{low}_{\kappa-2}(t)=\sup_{0\le\tau\le t}\widetilde{\mathcal J}_{\kappa-2}(\tau),
		\end{aligned}
	\end{equation}
	where
	\begin{equation}\nonumber
		\begin{aligned}
			&\widetilde{\mathcal I}_{\kappa}(\tau)=\mathcal I_{\kappa}(\tau)+X_{df,\kappa}(\tau)+X_{1,\kappa}(\tau)+X_{2,\kappa}(\tau),\\
			&\widetilde{\mathcal J}_{\kappa-2}(\tau)=\mathcal J_{\kappa-2}(\tau)+\mathcal J_{df,\kappa-2}(\tau)+\mathcal J_{1,\kappa-2}(\tau)+\mathcal J_{2,\kappa-2}(\tau).
		\end{aligned}
	\end{equation}
	
	\subsection{Conventions}
	
	\begin{enumerate}
		\item
		We write $A\lesssim B$ if $A\leq CB$ for some absolute constant $C>0$, and $A\sim B$ if both $A\lesssim B$ and $B\lesssim A$.
		
		\item
		We define the mixed Lebesgue norm by
		$$ \|f\|_{L_x^{p,q}} \triangleq \left(\int_0^\infty \left(\int_{\mathbb S^2}|f(r\theta)|^q d\mu_{\mathbb S^2}(\theta)\right)^{p/q}\,dr\right)^{1/p}, $$
		where $\mu_{\mathbb S^2}$ is the standard surface measure on the sphere.
		
		\item
		For nonnegative integers $i$, $j$, and $k$, we use
		$$ \nabla^{\le i}f\triangleq \{\partial_x^a f:|a|\le i\},\qquad  \Gamma^{\le j}f\triangleq\{\Gamma^\alpha f:|\alpha|\le j\},\qquad \widetilde\Gamma^{\le k}f\triangleq\{\widetilde\Gamma^\beta f:|\beta|\le k\}. $$
		These denote all spatial derivatives of order at most $i$, all compositions of at most $j$ vector fields from $\Gamma$, and all compositions of at most $k$ modified vector fields from $\widetilde\Gamma$, respectively.
		Each collection includes $f$ itself at order zero.
		In norms, these expressions denote sums of the corresponding norms:
		$$\|\nabla^{\le i}f\|_X\triangleq\sum_{|a|\le i}\|\partial_x^a f\|_X,\qquad
		\|\Gamma^{\le j}f\|_X\triangleq\sum_{|\alpha|\le j}\|\Gamma^\alpha f\|_X,\qquad
		\|\widetilde\Gamma^{\le k}f\|_X\triangleq\sum_{|\beta|\le k}\|\widetilde\Gamma^\beta f\|_X.$$
		
		\item
		We use the Japanese bracket notation
		$$\langle z\rangle:=(1+|z|^2)^{1/2}.$$
		
	\end{enumerate}
	
	\subsection{The global existence theorem}
	
	\begin{theorem}\label{thm:global-existence}
		
		Let $\kappa\ge 9$ and $0<\delta<\frac14$.
		Assume that the different null condition \eqref{equ_d_2} is satisfied.
		Then there exists $\varepsilon_0>0$ such that the Cauchy problem \eqref{equ_1} with initial data
		\begin{equation}\nonumber
			\nabla u_0^{df},\ u_1^{df} \in H_{\widetilde\Lambda}^{\kappa,2},\qquad
			\nabla u_0^{cf},\ u_1^{cf} \in H_{\widetilde\Lambda}^{\kappa,\frac65} \cap H_{\widetilde\Lambda}^{\kappa,2},\qquad
			\Theta_0 \in H_\Lambda^{\kappa,\frac65} \cap H_\Lambda^{\kappa,2},
		\end{equation}
		admits a unique global classical solution, provided that
		\begin{equation}
			\widetilde E_{\kappa-2}(0)\exp\left(C_0 \widetilde E_{\kappa}^{\frac12}(0)\right) \leq\varepsilon^2, \qquad 0<\varepsilon\leq\varepsilon_0.
		\end{equation}
		Moreover, the solution satisfies
		\begin{equation}
			\begin{aligned}
				\widetilde E_{\kappa-2}(t) \leq C_0^2\varepsilon^2, \qquad \widetilde E_{\kappa}(t) \leq C_0^2 \widetilde E_{\kappa}(0)\langle t\rangle^{C_0\varepsilon+\delta},
			\end{aligned}
		\end{equation}
		for all $t\geq0$, where $C_0>0$ is independent of $t$ and $\varepsilon$.
	\end{theorem}
	
	\section{Preliminaries}\label{2}
	\subsection{Commutation}
	We first record the basic commutation identities:
	\begin{equation}\nonumber
		\begin{aligned}
			[\partial_a,S]=[\partial_a,\widetilde S]=\partial_a,\quad
			[\partial_k,\Omega_{\ell}]=\varepsilon_{\ell kj}\partial_j,\quad
			\widetilde{\Omega}_{\ell}\nabla\Sigma=\nabla\Omega_{\ell}\Sigma,\quad
			\Omega_{\ell}\nabla\cdot u=\nabla \cdot(\widetilde{\Omega}_{\ell}u).
		\end{aligned}
	\end{equation}
	Writing $T=\partial_t-k\Delta$ and $L=(\partial_t^2-c^2_2\Delta)I-(c_1^2-c_2^2)\nabla \otimes \nabla$,
	we have 
	
	\begin{equation*}
		\begin{aligned}
			&L\widetilde S w+\nabla S\Sigma =\widetilde S\bigl(Lw+\nabla\Sigma\bigr) +2\bigl(Lw+\nabla\Sigma\bigr),\\
			&L\widetilde\Omega_\ell w+\nabla\Omega_\ell\Sigma =\widetilde\Omega_\ell\bigl(Lw+\nabla\Sigma\bigr),\\
			&TS\Sigma+\nabla\cdot\partial_t\widetilde S w =S\bigl(T\Sigma+\nabla\cdot\partial_tw\bigr)-\bigl(T\Sigma+\nabla\cdot\partial_tw\bigr)-k\Delta\Sigma,\\
			&T\Omega_\ell\Sigma+\nabla\cdot\partial_t\widetilde\Omega_\ell w =\Omega_\ell\bigl(T\Sigma+\nabla\cdot\partial_tw\bigr).
		\end{aligned}
	\end{equation*}
	
	Next, we record the commutation identities for the nonlinear terms.
	For rotations, the mechanical terms satisfy
	\begin{equation*}
		\begin{aligned}
			&\widetilde\Omega_\ell N_{al}(w_1,w_2)=N_{al}(\widetilde\Omega_\ell w_1,w_2)+N_{al}(w_1,\widetilde\Omega_\ell w_2),\\
			&\widetilde\Omega_\ell N_{th}(\Sigma_1,\Sigma_2)=N_{th}(\Omega_\ell\Sigma_1,\Sigma_2)+N_{th}(\Sigma_1,\Omega_\ell\Sigma_2),\\
			&\widetilde\Omega_\ell N_{mi}(w,\Sigma)=N_{mi}(\widetilde\Omega_\ell w,\Sigma)+N_{mi}(w,\Omega_\ell\Sigma).
		\end{aligned}
	\end{equation*}
	For the thermal terms, we have
	\begin{equation*}
		\begin{aligned}
			&\Omega_\ell M_1(\Sigma,w)=M_1(\Omega_\ell\Sigma,w) +M_1(\Sigma,\widetilde\Omega_\ell w),\\
			&\Omega_\ell M_2(\Sigma_1,\Sigma_2)=M_2(\Omega_\ell\Sigma_1,\Sigma_2)+M_2(\Sigma_1,\Omega_\ell\Sigma_2).
		\end{aligned}
	\end{equation*}
	For scaling, the displacement terms satisfy
	\begin{equation*}
		\begin{aligned}
			&\widetilde S N_{al}(w_1,w_2)=N_{al}(\widetilde S w_1,w_2)+N_{al}(w_1,\widetilde S w_2)-2N_{al}(w_1,w_2),\\
			&\widetilde S N_{th}(\Sigma_1,\Sigma_2)=N_{th}(S\Sigma_1,\Sigma_2)+N_{th}(\Sigma_1,S\Sigma_2)-2N_{th}(\Sigma_1,\Sigma_2),\\
			&\widetilde S N_{mi}(w,\Sigma)=N_{mi}(\widetilde S w,\Sigma) +N_{mi}(w,S\Sigma)-2N_{mi}(w,\Sigma).
		\end{aligned}
	\end{equation*}
	The thermal terms satisfy
	\begin{equation*}
		\begin{aligned}
			&S M_1(\Sigma,w)=M_1(S\Sigma,w)+M_1(\Sigma,\widetilde{S}w)-M_1(\Sigma,w),\\
			&S M_2(\Sigma_1,\Sigma_2)=M_2(S\Sigma_1,\Sigma_2)+M_2(\Sigma_1,S\Sigma_2)-2M_2(\Sigma_1,\Sigma_2).
		\end{aligned}
	\end{equation*}
	
	These identities lead to the following commuted system.
	\begin{lemma}\label{lem_com}
		The equations take the form
		\begin{equation}
			\left\lbrace 
			\begin{aligned} &L\widetilde{\Gamma}^{\alpha}u+\nabla\Gamma^\alpha\Theta=F_{1,\alpha} \\
				&T\Gamma^{\alpha}\Theta+\nabla\cdot\partial_t\widetilde{\Gamma}^{\alpha}u=F_{2,\alpha},
			\end{aligned}
			\right.
		\end{equation}
		where $\tilde\Gamma=(\tilde\Gamma_0,\ldots,\tilde\Gamma_7)=(\partial,\widetilde{\Omega},\widetilde{S})$, $\Gamma=(\Gamma_0,\ldots,\Gamma_7)=(\partial,\Omega,S)$, and
		\begin{equation}\label{equ_F1_F2}
			\left\lbrace 
			\begin{aligned} 
				&F_{1,\alpha}= \sum_{\gamma+\beta=\alpha}\left[N_{al}(\widetilde{\Gamma}^{\beta}u,\widetilde{\Gamma}^{\gamma}u)+N_{th}(\Gamma^{\beta}\Theta,\Gamma^{\gamma}\Theta)+N_{mi}(\widetilde{\Gamma}^{\beta}u,\Gamma^{\gamma}\Theta)\right]\\
				&F_{2,\alpha}=\sum_{|\beta+\gamma|\le |\alpha|}\left[C_{\beta,\gamma}^1M_1(\Gamma^{\beta}\Theta,\widetilde{\Gamma}^{\gamma}u)+C_{\beta,\gamma}^2M_2(\Gamma^{\beta}\Theta,\Gamma^{\gamma}\Theta) \right]+\sum_{\beta< \alpha}\left( C_{\beta}^2\Delta\Gamma^{\beta}\Theta\right).
			\end{aligned}
			\right.
		\end{equation}
	\end{lemma}
	Finally, we record the following commutation relations between the Helmholtz projections and the vector fields; see \cite[Lemma 2.2]{ZhaPeng2020}.
	\begin{lemma}
		We have
		\begin{equation*}
			(\widetilde{\Gamma}^\alpha u)^{cf}= \widetilde{\Gamma}^\alpha u^{cf},\qquad ( \widetilde{\Gamma}^\alpha u)^{df}= \widetilde{\Gamma}^\alpha u^{df}.
		\end{equation*}
	\end{lemma}
	
	\subsection{Null form estimates, Sobolev inequalities and other lemmas}
	We will use the following two estimates for the null forms.
	\begin{lemma}[\cite{Sideris2000}]\label{lem_null}
		We have
		\begin{equation}\label{eq:2.21}
			|N_2(u,v)|\leq \frac{C}{r}\sum_{|a|\leq 1} \left(|\widetilde{\Omega}^{a}u||\nabla^2v|+|\widetilde{\Omega}^{a}v||\nabla^2u| +|\nabla\widetilde{\Omega}^{a}u||\nabla v|+|\nabla\widetilde{\Omega}^{a}v||\nabla u|\right).
		\end{equation}
	\end{lemma}
	\begin{lemma}[ \cite{Sideris1996}]\label{lem_null_2}
		Let 
		\begin{equation}\nonumber
			N(u,v,w) = g^{ijk}_{lmn} \partial_l u_i\,\partial_m v_j\,\partial_n w_k,
		\end{equation}
		where the coefficients $g^{ijk}_{lmn}$ satisfy \eqref{equ_null_condition_2}. Then, for $r=|x|>0$, we have
		\begin{equation}
			|N(u,v,w)|\lesssim \frac{1}{r}\left(|\widetilde{\Omega}^{\leq 1} u||\nabla v||\nabla w|+|\nabla u||\widetilde{\Omega}^{\leq 1}v||\nabla w|+|\nabla u||\nabla v||\widetilde{\Omega}^{\leq 1} w|\right).
		\end{equation}
	\end{lemma}
	
	We will need the following weighted Sobolev-type inequalities.
	\begin{lemma}\label{lem_sob_1}
		Let $u\in\dot H_{\widetilde{\Lambda}}^{\kappa,2}(T)$ with $\mathcal J_{df,\kappa}(t)<+\infty$.
		Then the following estimates hold:
		\begin{equation}\label{equ_sob_1}
			r^{\frac12}\|\Gamma^a u(t,r\omega)\|_{L_\omega^2}\leq C E_\kappa^{\frac12}(u(t)),\qquad |a|+1\leq\kappa,
		\end{equation}
		\begin{equation}\label{equ_sob_2}
			\langle r\rangle^{\frac12}|\Gamma^a u(t,x)|\leq C E_\kappa^{\frac12}(u(t)),\qquad|a|+2\leq\kappa.
		\end{equation}
		\begin{equation}\label{equ_sob_3}
			\langle r\rangle|\partial\Gamma^a u(t,x)|\leq CE_\kappa^{\frac12}(u(t)), \qquad |a|+3\leq\kappa,
		\end{equation}
		\begin{equation}\label{equ_sob_7}
			\langle r\rangle\langle c_2t-r\rangle | \partial\nabla\Gamma^a u^{df}(t,x)| \leq C\mathcal J_{df,\kappa}^{\frac12}(t), \qquad |a|+4\leq\kappa.
		\end{equation}
		\begin{equation}\label{equ_sob_4}
			\langle c_2t-r\rangle |\partial\Gamma^a u^{df}(t,x)| \lesssim E_\kappa^{\frac12}(t) +\mathcal J_{df,\kappa}^{\frac12}(t), \qquad |a|+3\leq\kappa,
		\end{equation}
		\begin{equation}\label{equ_sob_6}
			\langle t\rangle |\partial\widetilde{\Gamma}^\alpha u^{df}(t,x)| \lesssim E_\kappa^{\frac12}(t) +\mathcal J_{df,\kappa}^{\frac12}(t), \qquad |\alpha|+3\leq\kappa,
		\end{equation}
		\begin{equation}\label{equ_sob_5}
			\left\| \langle c_2t-r\rangle \partial\Gamma^a u^{df}(t,\cdot) \right\|_{L_x^6} \lesssim E_\kappa^{\frac12}(t) +\mathcal J_{df,\kappa}^{\frac12}(t), \qquad |a|+2\leq\kappa,
		\end{equation}
		\begin{equation}\label{equ_sob_8}
			\left\| \langle c_2t-r\rangle^\frac12 \partial\Gamma^a u^{df}(t,\cdot) \right\|_{L_x^3}\lesssim E_\kappa^{\frac12}(t) +\mathcal J_{df,\kappa}^{\frac12}(t), \qquad |a|+2\leq\kappa,
		\end{equation}
	\end{lemma}
	\begin{proof}
		Estimates \eqref{equ_sob_3} and \eqref{equ_sob_7} follow from Lemma 2.6 in \cite{ZhaPeng2020}. Estimates \eqref{equ_sob_4}, \eqref{equ_sob_5} and \eqref{equ_sob_8} are immediate consequences of the interpolation inequalities and Sobolev inequalities,
		\begin{equation}\nonumber
			\|u\|_{L_x^\infty}\leq C\sum_{|\beta|\leq 1}\|\nabla\nabla^\beta u\|_{L_x^2}, \qquad \|u\|_{L_x^6}\leq C\|\nabla u\|_{L_x^2},
		\end{equation}
		respectively.
		Finally, combining \eqref{equ_sob_3} and \eqref{equ_sob_4} with the inequality
		$\langle t\rangle\lesssim \langle r\rangle+\langle c_2t-r\rangle $ yields \eqref{equ_sob_6}.
	\end{proof}
	
	We will use the following lemma in the subsequent analysis.
	\begin{lemma}\label{lem_Gamma_1}
		Let $v,w \in H_{\widetilde\Lambda}^{\kappa,2}$ with $\kappa\ge 1$. Then
		\begin{equation}\nonumber
			\begin{aligned}
				\sum_{|\alpha|=1} \int_{\mathbb R^3}\widetilde\Gamma^\alpha v\cdot w dx ={}& \sum_{|\alpha|=1}C_{\alpha,1} \partial_t\int_{\mathbb R^3}v\cdot w dx +\sum_{|\alpha|=1}C_{\alpha,2} \partial_t\left(t\int_{\mathbb R^3}v\cdot w dx\right)
				+ \sum_{|\alpha|\leq1}C_{\alpha,3} \int_{\mathbb R^3}v\cdot\widetilde\Gamma^\alpha w dx,
			\end{aligned}
		\end{equation}
		where $C_{\alpha,1}$, $C_{\alpha,2}$, and $C_{\alpha,3}$ are constants.
	\end{lemma}
	\begin{proof}
		Considering each type of vector field in $\widetilde\Gamma$, we obtain the desired result by integration by parts.
	\end{proof}
	
	\subsection{Klainerman--Sideris type estimates}
	
	We also need the Klainerman--Sideris type estimate for the wave operator, which was first established in \cite{KlainermanSideris1996}.
	
	\begin{lemma}\label{lem_KM_df}
		For $u\in \dot H_{\widetilde{\Lambda}}^{2,2}(T)$, we have
		\begin{equation}\nonumber
			\|\langle c_2 t-r\rangle \partial\nabla u^{df}\|_{L_x^2} \lesssim \sum_{|\alpha|\leq 1}\|\partial\widetilde\Gamma^{\alpha}u^{df}\|_{L_x^2} +t\|\Box_{c_2}u^{df}\|_{L_x^2},
		\end{equation}
		where $\Box_{c_2}=\partial_{t}^{2}-c_2^{2}\Delta$ is the wave operator with wave speed $c_2>0$.
	\end{lemma}
	We next apply the Klainerman--Sideris estimate in Lemma \ref{lem_KM_df} to control $X_{df,\kappa}(t)$ and $\mathcal J_{df,\kappa-2}(t)$.
	\begin{lemma}\label{lem_xd}
		Let $(u,\Theta)\in \dot H_{\widetilde\Gamma}^{\kappa}(T)\times\dot H_{\Gamma}^{\kappa}(T)$ with $\kappa\ge 9$ be a solution of \eqref{equ_1}. Then
		\begin{equation}\nonumber
			\begin{aligned}
				X_{df,\kappa}(t)&\lesssim \mathcal{I}_{\kappa}(t)+\widetilde{\mathcal{I}}_{\kappa}(t)\mathcal E^{low}_{\kappa-2}(t).
			\end{aligned}
		\end{equation}
	\end{lemma}
	\begin{proof}
		By the weighted Klainerman--Sideris estimate (see Lemma \ref{lem_KM_df}), we obtain
		\begin{align*}
			X_{df,\kappa}(t)\lesssim{} &\mathcal{I}_{\kappa}(t)+\sum_{|\beta|+|\gamma|\leq\kappa-2}\langle t\rangle^{2-\delta}\left( \left\|\nabla^{2}\widetilde{\Gamma}^{\beta}u\nabla\widetilde{\Gamma}^{\gamma}u\right\|_{L_x^2}^2+ \left\|\nabla(\nabla\widetilde{\Gamma}^{\beta}u\Gamma^{\gamma}\Theta)\right\|_{L_x^2}^2\right) .
		\end{align*}
		
		Set
		\begin{align*}
			A(u_{1},u_{2}):={}&\sum_{|\beta|+|\gamma|\leq\kappa-2}\langle t\rangle^{1-\frac{\delta}{2}} \left\|\nabla^{2}\widetilde{\Gamma}^{\beta}u_{1}\nabla\widetilde{\Gamma}^{\gamma}u_{2}\right\|_{L_x^2}.
		\end{align*}
		It remains to estimate $ A(u^{df},u)$, $A(u,u^{df})$ and $A(u^{cf},u^{cf})$.
		Using \eqref{equ_sob_3}, \eqref{equ_sob_7}, and the inequality $\langle t\rangle\lesssim\langle r\rangle+\langle c_2t-r\rangle$,
		we obtain
		\begin{align*}
			A(u^{df},u)\lesssim{}&\langle t \rangle ^{-\frac\delta2}\Big(\sum_{\substack{|\beta|+1\leq\kappa-5\\|\gamma|\leq\kappa-2}} \left\|\langle t \rangle\nabla^{2}\widetilde{\Gamma}^{\beta}u^{df}\right\|_{L_x^\infty}\left\|\nabla\widetilde{\Gamma}^{\gamma}u\right\|_{L_x^2}
			+\sum_{\substack{|\gamma|\leq\kappa-5\\|\beta|\leq\kappa-2}}\left\|\langle c_{2}t-r\rangle\nabla^{2}\widetilde{\Gamma}^{\beta}u^{df}\right\|_{L_x^2}\left\|\langle r\rangle\nabla\widetilde{\Gamma}^{\gamma}u\right\|_{L_x^\infty}\Big)\\
			\lesssim{}&\widetilde{\mathcal I}_{\kappa}^{\frac12}(t)\left(\mathcal E^{low}_{\kappa-2}(t)\right)^{\frac12}.
		\end{align*}
		To estimate $A(u,u^{df})$, we use \eqref{equ_sob_6} and \eqref{equ_sob_5} to obtain
		\begin{align*}
			A(u,u^{df})\lesssim{}& \sum_{\substack{|\beta|+1\leq\kappa-5\\|\gamma|\leq\kappa-2}} \langle t\rangle^{-\frac{\delta}{2}}\left( \left\|\langle c_2 t-r\rangle \nabla\widetilde{\Gamma}^{\gamma}u^{df}\right\|_{L_x^6} \left\|\nabla^2\widetilde{\Gamma}^{\beta}u\right\|_{L_x^3} +\left\|\nabla\widetilde{\Gamma}^{\gamma}u^{df}\right\|_{L_x^2}\left\|\langle r\rangle\nabla^2\widetilde{\Gamma}^{\beta}u\right\|_{L_x^\infty}\right) \\
			&+\sum_{\substack{|\gamma|\leq\kappa-5\\|\beta|\leq\kappa-2}}\left\|\langle t\rangle^{-\frac{\delta}{2}}\nabla^2\widetilde{\Gamma}^{\beta}u\right\|_{L_x^2} \left\|\langle t\rangle\nabla\widetilde{\Gamma}^{\gamma}u^{df}\right\|_{L_x^\infty}\\
			\lesssim{}&\widetilde{\mathcal I}_{\kappa}^{\frac12}(t)\left(\mathcal E^{low}_{\kappa-2}(t)\right)^{\frac12}.
		\end{align*}
		Applying the Sobolev inequality to $A(u^{cf},u^{cf})$, we obtain
		\begin{align*}
			A(u^{cf},u^{cf}) \lesssim{}& \langle t\rangle^{1-\frac\delta2}\bigg(\sum_{\substack{|\beta|+1\leq\kappa-5\\|\gamma|\leq\kappa-2}} \left\|  \nabla\widetilde{\Gamma}^{\gamma}u^{cf} \right\|_{L_x^2} \left\| \nabla^{2}\widetilde{\Gamma}^{\beta}u^{cf} \right\|_{L_x^\infty}+\sum_{\substack{|\gamma|\leq\kappa-5\\|\beta|\leq\kappa-2}}\left\|\nabla\widetilde{\Gamma}^{\gamma}u^{cf}\right\|_{L_x^\infty}\left\|\nabla^{2}\widetilde{\Gamma}^{\beta}u^{cf}\right\|_{L_x^2}\bigg)\\
			\lesssim{}&\widetilde{\mathcal I}_{\kappa}^{\frac12}(t)\left(\mathcal E^{low}_{\kappa-2}(t)\right)^{\frac12}.
		\end{align*}
		From the definition of the weighted energy, we obtain
		\begin{equation}\nonumber
			\begin{aligned}
				\sum_{|\beta|+|\gamma|\leq\kappa-2}\langle t\rangle^{1-\frac{\delta}{2}} \left\| \nabla(\nabla\widetilde{\Gamma}^{\beta}u \Gamma^\gamma \Theta )\right\|_{L_x^2}
				\sim \sum_{|\beta|+|\gamma|\leq\kappa-2}\langle t\rangle^{1-\frac{\delta}{2}}\left( \left\| \nabla^2 \widetilde{\Gamma}^{\beta}u^{cf} \nabla\widetilde{\Gamma}^{\gamma}u\right\|_{L_x^2}
				+\left\|\nabla^2 \widetilde{\Gamma}^{\beta}u \nabla\widetilde{\Gamma}^{\gamma}u^{cf}\right\|_{L_x^2}\right).
			\end{aligned}
		\end{equation}
		Then
		\begin{equation}
			\sum_{|\beta|+|\gamma|\leq\kappa-2}\langle t\rangle^{1-\frac{\delta}{2}} \left( \left\| \nabla^{2}\widetilde{\Gamma}^{\beta}u \nabla\widetilde{\Gamma}^{\gamma}u \right\|_{L_x^2}+\left\| \nabla^{2}\widetilde{\Gamma}^{\beta}u \Gamma^\gamma \Theta \right\|_{L_x^2}\right) \lesssim \widetilde{\mathcal{I}}_{\kappa}^{{}\frac12}(t) \left(\mathcal E^{low}_{\kappa-2}(t)\right)^{\frac12}.
		\end{equation}
		This completes the proof.
	\end{proof}
	
	The following lemma gives a bound for $\mathcal J_{df,\kappa-2}(t)$.
	\begin{lemma}
		Let $(u,\Theta)\in \dot H_{\widetilde\Gamma}^{\kappa}(T)\times\dot H_{\Gamma}^{\kappa}(T)$ with $\kappa\ge 9$ be a solution of \eqref{equ_1}. Then
		\begin{equation}\nonumber
			\begin{aligned}
				\mathcal J_{df,\kappa-2}(t) &\lesssim \mathcal{J}_{\kappa-2}(t)+ \left(\mathcal E^{low}_{\kappa-2}(t)\right)^{2}.
			\end{aligned}
		\end{equation}
	\end{lemma}
	\begin{proof}
		By Lemma \ref{lem_KM_df}, it remains to estimate the following nonlinear terms:
		\begin{equation}\nonumber
			\begin{aligned}
				\sum_{|\alpha|\leq\kappa-4} t^2 \left\| \Box_{c_2} \widetilde\Gamma^\alpha u^{df} \right\|_{L_x^2}^2 \lesssim \sum_{|\beta|+|\gamma|\leq\kappa-4} t^2\left( \left\| \nabla^2\widetilde\Gamma^\beta u \nabla\widetilde\Gamma^\gamma u \right\|_{L_x^2}^2+\left\| \nabla(\nabla\widetilde{\Gamma}^{\beta}u \Gamma^\gamma \Theta )\right\|^2_{L_x^2}\right) .
			\end{aligned}
		\end{equation}
		
		We distinguish the interactions according to the decomposition
		$u=u^{cf}+u^{df}$. Since $\left[\frac{\kappa-3}{2}\right]\leq\kappa-6$, we have
		\begin{equation}\nonumber
			\begin{aligned}
				\sum_{|\beta|+|\gamma|\leq\kappa-4}&t \left\| \nabla^2\widetilde\Gamma^\beta u \nabla\widetilde\Gamma^\gamma u \right\|_{L_x^2} \lesssim{}\sum_{\substack{ |\beta|\leq\kappa-4\\ |\gamma|\leq\kappa-6 }} 
				\left\| \nabla^2\widetilde\Gamma^\beta u \right\|_{L_x^2} \left\| \langle t\rangle \nabla\widetilde\Gamma^\gamma u^{df} \right\|_{L_x^\infty} \\
				&+\sum_{\substack{ |\gamma|\leq\kappa-4\\|\beta|+1\leq\kappa-6 }}
				\left( \left\| \nabla^2\widetilde\Gamma^\beta u \right\|_{L_x^3} \left\| \langle c_2t-r\rangle \nabla\widetilde\Gamma^\gamma u^{df} \right\|_{L_x^6}+ \left\|\langle r\rangle \nabla^2\widetilde\Gamma^\beta u \right\|_{L_x^\infty} \left\| \nabla\widetilde\Gamma^\gamma u^{df} \right\|_{L_x^2}\right) \\
				&+ \sum_{\substack{ |\beta|\leq\kappa-4\\|\gamma|\leq\kappa-6 }} 
				\left\| \langle c_2t-r\rangle \nabla^2\widetilde\Gamma^\beta u^{df} \right\|_{L_x^2} \left\| \langle r\rangle \nabla\widetilde\Gamma^\gamma u \right\|_{L_x^\infty} 
				+\sum_{\substack{ |\gamma|\leq\kappa-4\\|\beta|+1\leq\kappa-6 }} 
				\left\| \langle t \rangle \nabla^2\widetilde\Gamma^\beta u^{df} \right\|_{L_x^\infty} \left\| \nabla\widetilde\Gamma^\gamma u \right\|_{L_x^2}\\
				&+\sum_{\substack{ |\beta|\leq\kappa-4\\|\gamma|\leq\kappa-7 }} t \left\| \nabla^2\widetilde\Gamma^\beta u^{cf} \right\|_{L_x^2} \left\| \nabla\widetilde\Gamma^\gamma u^{cf} \right\|_{L_x^\infty} 
				+\sum_{\substack{ |\gamma|\leq\kappa-4\\|\beta|\leq\kappa-7 }} t \left\| \nabla^2\widetilde\Gamma^\beta u^{cf} \right\|_{L_x^\infty} \left\| \nabla\widetilde\Gamma^\gamma u^{cf} \right\|_{L_x^2} \\
				\lesssim{}&\mathcal E^{low}_{\kappa-2}(t).
			\end{aligned}
		\end{equation}
		Following the argument in Lemma \ref{lem_xd}, we complete the proof.
	\end{proof}
	
	\section{Weighted energy estimates for the curl-free and thermal components }\label{3}
	
	In this section, we derive estimates for the energy $X_{2,\kappa}(t)$ associated with $u^{cf}$ and $\Theta$. To this end, we consider the system
	\begin{equation}
		\left\lbrace 
		\begin{aligned}
			&\partial_t^2 w-c_1^2\Delta w+\nabla\Sigma =F_1(w,\Sigma), &&w\in\mathbb R^3, \\
			&\partial_t\Sigma-k\Delta\Sigma +\partial_t\nabla\cdot w =F_2(w,\Sigma), &&\Sigma\in\mathbb R.
		\end{aligned}
		\right.
	\end{equation}
	Here, $(w,\Sigma)\in \dot H_{\widetilde\Gamma}^{\kappa}(T) \times\dot H_{\Gamma}^{\kappa}(T)$ and $\nabla\times w=0$. Then we have
	\begin{equation*}
		\int_{\mathbb R^3}\left( \partial_t w \cdot F_1+\Sigma F_2\right) dx=\partial_t\int_{\mathbb R^3}\left( \frac{1}{2}|\partial_t w|^2+\frac{c_1^2}{2}|\nabla\cdot w|^2+\frac{1}{2}|\Sigma|^2\right) dx+\int_{\mathbb R^3}k|\nabla\Sigma|^2dx
	\end{equation*}
	\begin{equation}\nonumber
		\begin{aligned}
			&\int_{\mathbb R^3}\left( \partial_t \nabla\cdot w \nabla\cdot F_1+\nabla\Sigma\cdot \nabla F_2\right) dx=\partial_t\int_{\mathbb R^3}\left( \frac{1}{2}|\partial_t\nabla\cdot w|^2+\frac{c_1^2}{2}|\Delta w|^2+\frac{1}{2}|\nabla\Sigma|^2\right) dx+\int_{\mathbb R^3}k|\Delta\Sigma|^2dx\\
			&\int_{\mathbb R^3}\left( \Sigma\nabla\cdot F_1+\partial_t\nabla\cdot w F_2\right) dx=\partial_t\int_{\mathbb R^3}\Sigma\partial_t \nabla\cdot wdx-\int_{\mathbb R^3}\left( c_1^2\Sigma\nabla\cdot\Delta w+|\nabla\Sigma|^2+k\Delta \Sigma\partial_t\nabla\cdot w-|\partial_t \nabla w|^2 \right) dx\\
			&\int_{\mathbb R^3}\nabla\cdot w \nabla\cdot F_1dx=\partial_t\int_{\mathbb R^3}\partial_t \nabla\cdot w \nabla\cdot wdx+\int_{\mathbb R^3}\left( c_1^2|\Delta w|^2+\nabla\cdot w\Delta \Sigma-|\partial_t\nabla w|^2\right) dx,
		\end{aligned}
	\end{equation}
	and combining these identities with suitable coefficients yields
	\begin{equation}\nonumber
		\begin{aligned}
			\partial_tK_1(t)+K_2(t)=K_3(t),
		\end{aligned}
	\end{equation}
	where
	\begin{equation}\nonumber
		\begin{aligned}
			K_{1}(t) ={}& \int_{\mathbb R^3} \Bigg[ \frac{1}{2} \bigl(\Sigma+\partial_t\nabla\cdot w\bigr)^2 +\frac{\sigma}{2} \bigl(\partial_t\nabla\cdot w+\nabla\cdot w\bigr)^2 +\frac{\sigma^{-3}}{2}|\partial_t w|^2 +\frac{c_1^2\sigma^{-3}-\sigma}{2} |\nabla\cdot w|^2 \\
			&\qquad +\frac{\sigma^{-3}-1}{2}|\Sigma|^2 +\frac{\sigma^{-3}-1-\sigma}{2} |\partial_t\nabla\cdot w|^2 +\frac{c_1^2\sigma^{-3}}{2}|\Delta w|^2 +\frac{\sigma^{-3}}{2}|\nabla\Sigma|^2 \Bigg]\,dx, \\
			K_2(t) ={}& \int_{\mathbb R^3} \Bigg[ |a_1\nabla\Sigma+a_2\Delta w|^2 +|b_1\Delta\Sigma-b_2\partial_t\nabla\cdot w|^2 +\bigl(k\sigma^{-3}-a_1^2-1\bigr)|\nabla\Sigma|^2 \\
			&\qquad +\bigl(\sigma c_1^2-a_2^2\bigr)|\Delta w|^2 +\bigl(1-\sigma-b_2^2\bigr) |\partial_t\nabla\cdot w|^2 +\bigl(k\sigma^{-3}-b_1^2\bigr)|\Delta\Sigma|^2 \Bigg]\,dx, \\
			K_3(t) ={}& \int_{\mathbb R^3} \Bigg[\sigma^{-3} \Big(\partial_t\nabla\cdot w\,\nabla\cdot F_1 +\nabla\Sigma\cdot\nabla F_2 +\partial_t w\cdot F_1 +\Sigma F_2 \Big)
			+\Sigma\nabla\cdot F_1 +\partial_t\nabla\cdot w\,F_2 +\sigma\nabla\cdot w\,\nabla\cdot F_1\Bigg]\,dx.
		\end{aligned}
	\end{equation}
	
	For sufficiently small $\sigma>0$, we can choose $a_i,b_i>0$, $i=1,2$, satisfying $2a_1a_2=c_1^2-\sigma$ and $2b_1b_2=k$ so that $K_1(t)$ and $K_2(t)$ are nonnegative. 
	Let
	\begin{equation}\nonumber
		\begin{aligned}
			S_1(t)=S_1(w(t),\Sigma(t))&=\int_{\mathbb R^3}\left( |\partial_t w|^2 +|\nabla\cdot w|^2+|\Sigma|^2 +|\partial_t \nabla\cdot w|^2+ |\Delta w|^2 + |\nabla\Sigma|^2\right) dx,\\
			S_2(t)=S_2(w(t),\Sigma(t))&=\int_{\mathbb R^3}\left(|\nabla\Sigma|^2+|\Delta w |^2+|\partial_t \nabla \cdot w|^2+ |\Delta\Sigma|^2\right) dx.
		\end{aligned}
	\end{equation}
	There exists a sufficiently small $\sigma>0$ such that
	\begin{equation}\label{equ_2.20}
		S_1(t)\leq K_1(t)\leq \sigma^{-4}S_1(t), \qquad \sigma^2S_2(t)\leq K_2(t)\leq \sigma^{-4}S_2(t).
	\end{equation}
	Then
	\begin{equation}\label{equ_cf_01}
		\begin{aligned}
			S_1(t) +\int_0^t S_2(\tau)d\tau \lesssim_{\sigma} \int_0^t K_3(\tau)d\tau +\sigma^{-4}S_1(0).
		\end{aligned}
	\end{equation}
	
	Also, we obtain
	\begin{equation}\nonumber
		\begin{aligned}
			\int_0^t (\tau+\sigma^{-10}) \left(\partial_\tau K_1(\tau)+K_2(\tau)\right)d\tau = \int_0^t (\tau+\sigma^{-10})K_3(\tau)d\tau.
		\end{aligned}
	\end{equation}
	By
	$$ S_1(t) \leq S_2\left(w(t),\Sigma(t)\right) +S_2\left(\nabla^{-1}w(t),\nabla^{-1}\Sigma(t)\right),$$
	and \eqref{equ_2.20}, we have
	\begin{equation}\nonumber
		\begin{aligned}
			(t+\sigma^{-10})& S_1\left(t\right) + \int_0^t \left(\sigma^2\tau+\sigma^{-8}-\sigma^{-4}\right) S_2(\tau)d\tau \\
			&\leq \int_0^t (\tau+\sigma^{-10})K_3(\tau)d\tau+ \int_0^t \sigma^{-4} S_2\left( \nabla^{-1}w(\tau), \nabla^{-1}\Sigma(\tau) \right)d\tau+ \sigma^{-14} S_1\left(0\right).
		\end{aligned}
	\end{equation}
	That is,
	\begin{equation}\label{equ_cf_02}
		\begin{aligned}
			\langle t\rangle &S_1\left(t\right)+ \int_0^t \langle \tau \rangle S_2\left(\tau\right)d\tau
			\lesssim \int_0^t \Big[(\tau+\sigma^{-10})K_3(\tau)+ \sigma^{-4} S_2\left( \nabla^{-1}w, \nabla^{-1}\Sigma \right)(\tau)\Big]d\tau+ \sigma^{-14} S_1\left(0\right).
		\end{aligned}
	\end{equation}
	We next use \eqref{equ_cf_01} and \eqref{equ_cf_02} to estimate $X_{1,\kappa}(t)$, $X_{2,\kappa}(t)$, $\mathcal J_{1,\kappa-2}(t)$, and $\mathcal J_{2,\kappa-2}(t)$.
	
	\subsection{Higher-order weighted energy estimates}
	
	\begin{lemma}\label{lem_i_2}
		Let $\kappa\ge 9$, and let $(u,\Theta)\in \dot H_{\widetilde\Gamma}^{\kappa}(T) \times\dot H_{\Gamma}^{\kappa}(T)$ be a solution of \eqref{equ_1}. 
		Then
		\begin{equation}\nonumber
			X_{2,\kappa}(t)\lesssim \mathcal E_\kappa(t) \left(\mathcal E^{low}_{\kappa-2}(t)\right)^{\frac12} + \widetilde{\mathcal I}_{\kappa}^{\frac12}(t) \widetilde{\mathcal I}_{\kappa-1}^{\frac12}(t)+X_{1,\kappa}(t)+X_{2,\kappa}(0).
		\end{equation}
	\end{lemma}
	\begin{proof}
		Substituting
		$(w,\Sigma)=(\widetilde{\Gamma}^\alpha u^{cf},\Gamma^\alpha\Theta)$, $F_1=F_{1,\alpha}$ and $F_2=F_{2,\alpha}$
		(see \eqref{equ_F1_F2}) into \eqref{equ_cf_02}, we obtain
		\begin{equation}\nonumber
			\begin{aligned}
				X_{2,\kappa}(t) \lesssim{}& \sum_{|\alpha|\leq\kappa-2}\langle t\rangle^{-\delta} \int_0^t (\tau+\sigma^{-10}) K_{3,\alpha}(\tau) d\tau+\sigma^{-14}X_{2,\kappa}(0)+\sigma^{-4}X_{1,\kappa}(t).
			\end{aligned}
		\end{equation}
		Here, we define
		\begin{equation}\label{equ_K_3}
			\begin{aligned}
				K_{3,\alpha}(\tau)=&\int_{\mathbb R^3} \Big[A_1 \big(\nabla\cdot \partial_\tau\widetilde\Gamma^\alpha u^{cf} \nabla\cdot F_{1,\alpha} +\nabla\Gamma^\alpha\Theta \cdot \nabla F_{2,\alpha}+ \partial_\tau\widetilde\Gamma^\alpha u^{cf} \cdot F_{1,\alpha} +\Gamma^\alpha\Theta F_{2,\alpha}\big)\\
				&+A_2\big(\Gamma^\alpha\Theta \nabla\cdot F_{1,\alpha} +\nabla\cdot \partial_\tau\widetilde\Gamma^\alpha u^{cf} \cdot F_{2,\alpha}\big) +A_3\nabla\cdot \widetilde\Gamma^\alpha u^{cf} \nabla\cdot F_{1,\alpha} \Big]dx.
			\end{aligned}
		\end{equation}
		We estimate the integrals on the right-hand side separately.
		For the first term involving $F_{1,\alpha}$, the divergence $\nabla\cdot F_{1,\alpha}$ contains one additional derivative compared with $F_{1,\alpha}$ itself. 
		Applying H\"older's inequality and Lemma \ref{lem_sob_1}, we obtain
		
		\begin{equation}\label{equ_F11_1}
			\begin{aligned}
				\sum_{|\alpha|\leq\kappa-2}&\langle t \rangle^{-\delta} \int_0^t\int_{\mathbb R^3} (\tau+\sigma^{-10}) \nabla\cdot \partial_\tau \widetilde\Gamma^\alpha u^{cf} \nabla \cdot F_{1,\alpha} dx d\tau 
				\lesssim{} \langle t\rangle^{-\delta} \sum_{\substack{ |\alpha|\leq\kappa-2\\
						|\beta|+|\gamma|\leq\kappa-2}} \left\| \langle\tau\rangle^{\frac12}\nabla\cdot \partial_\tau\widetilde\Gamma^\alpha u^{cf} \right\|_{L_\tau^2 L_x^2}\\
				&\times\Biggl(\left\|\langle \tau\rangle^{\frac12} \nabla(\nabla^2\widetilde{\Gamma}^\beta u\nabla\widetilde{\Gamma}^\gamma u) \right\|_{L_\tau^2 L_x^2}+\left\|\langle \tau\rangle^{\frac12} \nabla^2(\nabla\widetilde{\Gamma}^\beta u \Gamma^{\gamma}\Theta) \right\|_{L_\tau^2 L_x^2}+\left\|\langle \tau\rangle^{\frac12} \nabla(\nabla\Gamma^{\beta}\Theta \Gamma^{\gamma}\Theta) \right\|_{L_\tau^2 L_x^2} \Biggr).
			\end{aligned}
		\end{equation}
		For the nonlinear elastic term, the Leibniz rule yields
		
		\begin{equation}\label{equ_F11_1.1}
			\begin{aligned}
				&\langle t\rangle^{-\frac{\delta}{2}}\sum_{\substack{|\beta|+|\gamma|\le\kappa-1,|a_1|+|a_2|\le 1\\ |\beta|+|a_1|\le \kappa-1}}\left\|\langle \tau\rangle^{\frac12} \nabla^{2+a_1}\widetilde{\Gamma}^\beta u\nabla^{1+a_2}\widetilde{\Gamma}^\gamma u \right\|_{L_\tau^2 L_x^2}\\
				&\lesssim  \langle t\rangle^{-\frac{\delta}{2}}\Bigg[\sum_{\substack{|\beta|\le \kappa-1\\|\gamma|\le \kappa-5}} 
				\left\| \langle\tau\rangle^{-\frac\delta2} \nabla^{2}\widetilde\Gamma^\beta u \right\|_{L_\tau^\infty L_x^2}\left\| \langle\tau\rangle^{\frac12+\frac{\delta}{2}} \nabla\widetilde\Gamma^\gamma u \right\|_{L_\tau^2 L_x^\infty}
				+\sum_{\substack{|\beta|\le 5 \\|\gamma|\le \kappa-4}} 
				\left\| \langle\tau\rangle^{\frac12+\frac{\delta}{2}} \nabla^{2}\widetilde\Gamma^\beta u \right\|_{L_\tau^2 L_x^\infty}\left\|\nabla\widetilde\Gamma^\gamma u \right\|_{L_\tau^\infty L_x^2}\\
				&\quad+\sum_{\substack{|\beta|\le \kappa-6 \\|\gamma|\le \kappa-1}} 
				\left\| \langle\tau\rangle^{\frac12+\frac{\delta}{2}} \nabla^{2}\widetilde\Gamma^\beta u \right\|_{L_\tau^2 L_x^\infty}\left\| \langle\tau\rangle^{-\frac\delta2} \nabla\nabla^{\le 1}\widetilde\Gamma^\gamma u \right\|_{L_\tau^\infty L_x^2}
				+\sum_{\substack{|\beta|\le \kappa-5 \\|\gamma|\le 5}} 
				\left\| \nabla^{2}\widetilde\Gamma^\beta u \right\|_{L_\tau^\infty L_x^2}\left\| \langle\tau\rangle^{\frac12}  \nabla\widetilde\Gamma^\gamma u \right\|_{L_\tau^2 L_x^\infty}\Bigg]\\
				&\lesssim \mathcal E_\kappa^\frac12(t) \left(\mathcal E^{low}_{\kappa-2}(t)\right)^{\frac12}.
			\end{aligned}
		\end{equation}
		Here, we use the following estimates:
		\begin{equation}\label{equ_1.3}
			\sum_{|\gamma|\le \kappa -5}\langle t\rangle^{-\frac\delta2}\left\|\langle \tau \rangle^{\frac12+\frac{\delta}{2}} \nabla \widetilde{\Gamma}^\gamma u\right\| _{L_\tau^2 L_x^\infty} \lesssim{}\sum_{|\gamma|\le \kappa-5}\left( \left\|\langle \tau \rangle^{\frac12} \nabla \widetilde{\Gamma}^\gamma u^{cf}\right\| _{L_\tau^2 L_x^\infty}+\left\|\langle \tau \rangle \nabla \widetilde{\Gamma}^\gamma u^{df}\right\| _{L_\tau^\infty L_x^\infty}\right) 
			\lesssim{}\left(\mathcal E^{low}_{\kappa-2}(t)\right)^{\frac12},
		\end{equation}
		and
		\begin{equation}\label{equ_1.4}
			\begin{aligned}
				\sum_{|\gamma|\le \kappa-3}&\langle t\rangle^{-\frac\delta2}\left\|\langle \tau \rangle^{\frac12} \nabla \widetilde{\Gamma}^\gamma u\right\| _{L_\tau^2 L_x^\infty} 
				\lesssim{}\sum_{|\gamma|\le \kappa-3}\left(\langle t\rangle ^{-\frac\delta2} \left\|\langle \tau \rangle^{\frac12} \nabla \widetilde{\Gamma}^\gamma u^{cf}\right\| _{L_\tau^2 L_x^\infty}+\left\|\langle \tau \rangle^{1-\frac\delta2} \nabla \widetilde{\Gamma}^\gamma u^{df}\right\| _{L_\tau^\infty L_x^\infty}\right) 
				\lesssim{}\mathcal E_\kappa^{\frac12}(t).
			\end{aligned}
		\end{equation}
		Next, we consider the nonlinear terms describing the interaction between $u$ and $\Theta$. Using \eqref{equ_1.3}, we obtain
		
		\begin{equation}\label{equ_F11_1.2}
			\begin{aligned}
				\langle &t\rangle^{-\frac{\delta}{2}}\sum_{|\beta|+|\gamma|\le\kappa-2}\left\|\langle \tau\rangle^{\frac12} \nabla^2(\nabla\widetilde{\Gamma}^\beta u \Gamma^\gamma \Theta) \right\|_{L_\tau^2 L_x^2}\lesssim\langle t\rangle^{-\frac{\delta}{2}}\sum_{|\beta|+|\gamma|\le\kappa-1}\left\|\langle \tau\rangle^{\frac12} \nabla(\nabla\widetilde{\Gamma}^\beta u \Gamma^\gamma \Theta) \right\|_{L_\tau^2 L_x^2}\\
				&\lesssim \langle t\rangle^{-\frac\delta2}\sum_{|\beta|+|\gamma|\le\kappa-1}\Bigl( \left\| \langle\tau\rangle^{\frac12} \nabla^2\widetilde{\Gamma}^\beta u \Gamma^\gamma \Theta \right\|_{L_\tau^2 L_x^2}+ \left\| \langle\tau\rangle^{\frac12} \nabla\widetilde{\Gamma}^\beta u \nabla\Gamma^\gamma \Theta \right\|_{L_\tau^2 L_x^2}\Bigr) \\
				&\lesssim \sum_{\substack{|\gamma|\le \kappa-1\\
						|\beta|+1\le \kappa-5}}\langle t\rangle^{-\frac\delta2} \left\| \langle\tau\rangle^{\frac12+\frac\delta2} \nabla^2\widetilde\Gamma^\beta u \right\|_{L_\tau^2 L_x^\infty}\left\| \langle\tau\rangle^{-\frac\delta2} \Gamma^\gamma \Theta \right\|_{L_\tau^\infty L_x^2}+\sum_{\substack{|\beta|\le \kappa-1\\
						|\gamma|\le \kappa-5}} \left\| \langle\tau\rangle^{-\frac\delta2} \nabla^3\widetilde\Gamma^\beta u \right\|_{L_\tau^\infty L_x^2}\left\| \langle\tau\rangle^{\frac12} \Gamma^\gamma \Theta \right\|_{L_\tau^2 L_x^\infty} \\
				&\quad +\sum_{\substack{|\gamma|\le \kappa-1\\
						|\beta|\le \kappa-5}}\langle t\rangle^{-\frac\delta2} \left\| \langle\tau\rangle^{\frac12+\frac\delta2} \nabla\widetilde\Gamma^\beta u \right\|_{L_\tau^2 L_x^\infty} \left\| \langle\tau\rangle^{-\frac\delta2} \nabla\Gamma^\gamma \Theta \right\|_{L_\tau^\infty L_x^2}+\sum_{\substack{|\beta|\le \kappa-1\\
						|\gamma|\le \kappa-5}} \left\| \langle\tau\rangle^{-\frac\delta2} \nabla\widetilde\Gamma^\beta u \right\|_{L_\tau^\infty L_x^2} \left\| \langle\tau\rangle^{\frac12} \nabla\Gamma^\gamma \Theta \right\|_{L_\tau^2 L_x^\infty} \\
				&\lesssim \mathcal E_\kappa^\frac12(t) \left(\mathcal E^{low}_{\kappa-2}(t)\right)^{\frac12}.
			\end{aligned}
		\end{equation}
		For the interaction terms involving $\Theta$, we have
		\begin{equation}\label{equ_F11_1.3}
			\begin{aligned}
				\langle t\rangle^{-\frac{\delta}{2}}&\sum_{|\beta|+|\gamma|\le\kappa-1} \left\|\langle \tau\rangle^{\frac12}\nabla^{\le1}(\nabla\Gamma^{\beta}\Theta \Gamma^{\gamma}\Theta )\right\|_{L_\tau^2 L_x^2}
				\lesssim{} \langle t\rangle^{-\frac{\delta}{2}}\sum_{\substack{|\beta|+1\le\kappa-5\\ |\gamma|\le \kappa-1}}\left\|\langle \tau\rangle^{\frac12}\nabla\nabla^{\le1}\Gamma^{\beta}\Theta \right\|_{L_\tau^\infty L_x^3} \left\|\nabla^{\le1}\Gamma^{\gamma}\Theta \right\|_{L_\tau^2 L_x^6}\\
				&+\langle t\rangle^{-\frac{\delta}{2}}\sum_{\substack{|\gamma|\le \kappa-5\\|\beta|\le \kappa-1}}
				\Big(\left\|\langle \tau \rangle^\frac12 \Gamma^{\gamma}\Theta \right\|_{L_\tau^\infty L_x^\infty} \left\| \nabla\nabla^{\le 1}\Gamma^{\beta}\Theta \right\|_{L_\tau^2 L_x^2}+\left\|\langle \tau \rangle^\frac12\nabla^{\le1} \Gamma^{\gamma}\Theta \right\|_{L_\tau^\infty L_x^3} \left\| \nabla\Gamma^{\beta}\Theta \right\|_{L_\tau^2 L_x^6}\Big)\\
				\lesssim {}& \mathcal E_\kappa^\frac12(t) \left(\mathcal E^{low}_{\kappa-2}(t)\right)^{\frac12}.
			\end{aligned}
		\end{equation}
		Substituting \eqref{equ_F11_1.1}, \eqref{equ_F11_1.2}, and
		\eqref{equ_F11_1.3} into \eqref{equ_F11_1}, we obtain
		\begin{equation}\label{equ_F11_1-1}
			\begin{aligned}
				\sum_{|\alpha|\leq\kappa-2}&\langle t \rangle^{-\delta} \int_0^t\int_{\mathbb R^3} (\tau+\sigma^{-10})\nabla\cdot \partial_\tau \widetilde\Gamma^\alpha u^{cf} \nabla \cdot F_{1,\alpha} dx d\tau \lesssim{} \mathcal{E}_{\kappa}(t) \left(\mathcal E^{low}_{\kappa-2}(t)\right)^{\frac12}.
			\end{aligned}
		\end{equation}
		
		Next, we consider the second term involving $F_{1,\alpha}$.
		Since $F_{1,\alpha}$ has a divergence structure, we integrate by parts in space to make use of the $L^2_\tau L^2_x$ estimate for $u^{cf}$.
		Estimating the resulting terms separately, we obtain
		\begin{equation}\label{equ_F11_2}
			\begin{aligned}
				\langle t\rangle^{-\delta}& \sum_{|\alpha|\leq\kappa-2} \int_0^t\int_{\mathbb R^3}(\tau+\sigma^{-10}) \partial_\tau\widetilde\Gamma^\alpha u^{cf} \cdot F_{1,\alpha} dx d\tau 
				\lesssim{} \sum_{\substack{ |\alpha|\leq\kappa-2\\
						|\beta|+|\gamma|\leq\kappa-2}}\langle t\rangle^{-\delta} \left\|\langle\tau\rangle^{\frac12}\partial_\tau \nabla\widetilde\Gamma^\alpha u^{cf}\right\|_{L_\tau^2 L_x^2}\\
				&\times \biggl( \bigl\|\langle\tau\rangle^{\frac12} \nabla\widetilde\Gamma^\beta u\, \nabla\widetilde\Gamma^\gamma u \bigr\|_{L_\tau^2 L_x^2} + \bigl\|\langle\tau\rangle^{\frac12} \nabla\widetilde\Gamma^\beta u\, \Gamma^\gamma\Theta \bigr\|_{L_\tau^2 L_x^2}+ \bigl\|\langle\tau\rangle^{\frac12} \Gamma^\beta\Theta\,\Gamma^\gamma\Theta \bigr\|_{L_\tau^2 L_x^2} \biggr) \\
				\lesssim{}& \mathcal{E}_{\kappa}(t) \left(\mathcal E^{low}_{\kappa-2}(t)\right)^{\frac12}
			\end{aligned}
		\end{equation}
		Here, we have used the following estimates, which follow from \eqref{equ_1.3},
		\begin{equation}\label{equ_hig_F1_1}
			\begin{aligned}
				\sum_{|\beta|+|\gamma|\leq\kappa-2}& \langle t\rangle^{-\frac\delta2} \bigl\|\langle\tau\rangle^{\frac12} \nabla\widetilde{\Gamma}^{\beta}u\, \nabla\widetilde{\Gamma}^{\gamma}u \bigr\|_{L_\tau^2 L_x^2} \lesssim \sum_{\substack{|\beta|\leq\kappa-2\\|\gamma|\leq\kappa-6}}\langle t\rangle^{-\frac\delta2} \bigl\|\langle\tau\rangle^{-\frac\delta2} \nabla\widetilde{\Gamma}^{\beta}u \bigr\|_{L_\tau^\infty L_x^2} \bigl\|\langle\tau\rangle^{\frac\delta2+\frac12} \nabla\widetilde{\Gamma}^{\gamma}u \bigr\|_{L_\tau^2 L_x^\infty} \\
				\lesssim{}& \mathcal E^\frac12_\kappa(t) \left(\mathcal E^{low}_{\kappa-2}(t)\right)^{\frac12},
			\end{aligned}
		\end{equation}
		and
		\begin{equation}
			\begin{aligned}
				&\sum_{|\beta|+|\gamma|\leq\kappa-2} \langle t\rangle^{-\frac\delta2} \bigl\|\langle\tau\rangle^{\frac12} \nabla\widetilde{\Gamma}^{\beta}u\, \Gamma^\gamma\Theta \bigr\|_{L_\tau^2 L_x^2} \\
				&\quad\lesssim \sum_{\substack{|\beta|\leq\kappa-2\\|\gamma|\leq\kappa-6}} \bigl\|\langle\tau\rangle^{-\frac\delta2} \nabla\widetilde{\Gamma}^{\beta}u \bigr\|_{L_\tau^\infty L_x^2} \bigl\|\langle\tau\rangle^{\frac12} \Gamma^\gamma\Theta\bigr\|_{L_\tau^2 L_x^\infty} + \sum_{\substack{|\beta|\leq\kappa-6\\|\gamma|\leq\kappa-2}} \langle t\rangle^{-\frac\delta2} \bigl\|\nabla\widetilde{\Gamma}^{\beta}u\bigr\|_{L_\tau^\infty L_x^\infty} \bigl\|\langle\tau\rangle^{\frac12} \Gamma^\gamma\Theta\bigr\|_{L_\tau^2 L_x^2} \\
				&\quad\lesssim \mathcal E^\frac12_\kappa(t) \left(\mathcal E^{low}_{\kappa-2}(t)\right)^{\frac12},
			\end{aligned}
		\end{equation}
		and
		\begin{equation}\label{equ_hig_F1_3}
			\begin{aligned}
				&\sum_{|\beta|+|\gamma|\leq\kappa-2} \langle t\rangle^{-\frac\delta2} \bigl\|\langle\tau\rangle^{\frac12} \Gamma^\beta\Theta\,\Gamma^\gamma\Theta \bigr\|_{L_\tau^2 L_x^2} \lesssim \sum_{\substack{|\beta|\leq\kappa-2\\|\gamma|\leq\kappa-6}} \bigl\|\langle\tau\rangle^{-\frac\delta2} \Gamma^\beta\Theta\bigr\|_{L_\tau^\infty L_x^2} \bigl\|\langle\tau\rangle^{\frac12} \Gamma^\gamma\Theta\bigr\|_{L_\tau^2 L_x^\infty} 
				\lesssim \mathcal E^\frac12_\kappa(t) \left(\mathcal E^{low}_{\kappa-2}(t)\right)^{\frac12}.
			\end{aligned}
		\end{equation}
		We estimate the remaining two terms involving $F_{1,\alpha}$ together.
		Using \eqref{equ_F11_1.1}, \eqref{equ_F11_1.2}, and \eqref{equ_F11_1.3}, we obtain
		\begin{equation}\label{equ_F11_3}
			\begin{aligned}
				\langle t\rangle^{-\delta}&\sum_{|\alpha|\leq\kappa-2}\int_0^t\int_{\mathbb R^3}(\tau+\sigma^{-10})\left( A_2\Gamma^\alpha \Theta +A_3\nabla \cdot \widetilde{\Gamma}^\alpha u^{cf}\right) \nabla\cdot F_{1,\alpha} 
				\lesssim{} \mathcal E_\kappa(t) \left(\mathcal E^{low}_{\kappa-2}(t)\right)^{\frac12}.
			\end{aligned}
		\end{equation}
		
		We now turn to the terms involving $F_{2,\alpha}$.
		Integrating by parts in space and combining the preceding estimates with \eqref{equ_F11_1.3}, we obtain
		\begin{equation}\nonumber
			\begin{aligned}
				\langle t\rangle^{-\delta}& \sum_{|\alpha|\leq\kappa-2} \int_0^t\int_{\mathbb R^3}(\tau+\sigma^{-10})\left(A_1\nabla \Gamma^\alpha \Theta\cdot\nabla F_{2,\alpha}+A_2 \nabla\cdot \partial_\tau\widetilde\Gamma^\alpha u^{cf} F_{2,\alpha}\right) dx d\tau \\
				\lesssim{}& \langle t\rangle^{-\delta} \sum_{|\alpha|\leq\kappa-2}\left( \left\| \langle\tau\rangle^{\frac12}\nabla \cdot \partial_\tau \widetilde\Gamma^\alpha u^{cf} \right\|_{L_\tau^2 L_x^2}+\left\| \langle \tau \rangle^{\frac12}\Delta \Gamma^\alpha \Theta\right\| _{L_\tau^2 L_x^2} \right) 
				\Biggl[ \sum_{|\beta|\le \kappa-3}\left\| \langle \tau \rangle ^{\frac12}\Delta \Gamma^\beta \Theta\right\| _{L_\tau^2 L_x^2}\\
				&+\sum_{|\beta|+|\gamma|\leq\kappa-2}\Big(\left\| \langle\tau\rangle^{\frac12} \Gamma^\beta\Theta \nabla\cdot\partial_\tau \widetilde\Gamma^\gamma u^{cf} \right\|_{L_\tau^2 L_x^2}+ \left\|\langle \tau\rangle^{\frac12} \nabla(\nabla\Gamma^{\beta}\Theta \Gamma^{\gamma}\Theta) \right\|_{L_\tau^2 L_x^2}\Big)\Biggr]\\
				\lesssim{}&\mathcal E_\kappa(t) \left(\mathcal E^{low}_{\kappa-2}(t)\right)^{\frac12} +\widetilde{\mathcal I}_{\kappa}^{\frac12}(t) \widetilde{\mathcal I}_{\kappa-1}^{\frac12}(t).
			\end{aligned}
		\end{equation}
		Here, the quadratic term is estimated as follows:
		\begin{equation}\nonumber
			\begin{aligned}
				&\sum_{|\beta|+|\gamma|\leq\kappa-2}\langle t\rangle^{-\frac\delta2} \left\| \langle\tau\rangle^{\frac12} \Gamma^\beta\Theta \nabla\cdot\partial_\tau \widetilde\Gamma^\gamma u^{cf} \right\|_{L_\tau^2 L_x^2}\\
				&\quad\lesssim\langle t\rangle^{-\frac\delta2}\Bigg(\sum_{\substack{ |\beta|\leq\kappa-5\\
						|\gamma|\leq\kappa-2}} \left\|  \Gamma^\beta\Theta \right\|_{L_\tau^\infty L_x^\infty} \left\| \langle\tau\rangle^{\frac12} \nabla\partial_\tau \widetilde\Gamma^\gamma u^{cf} \right\|_{L_\tau^2 L_x^2}+ \sum_{\substack{ |\gamma|+2\leq\kappa-5\\|\beta|\leq\kappa-2}} 
				\left\| \nabla\partial_\tau \widetilde\Gamma^\gamma u^{cf} \right\|_{L_\tau^\infty L_x^2} \left\| \langle\tau\rangle^{\frac12} \Gamma^\beta\Theta \right\|_{L_\tau^2 L_x^\infty} \Bigg)\\
				&\quad\lesssim{}\mathcal E_\kappa^{\frac12}(t)\left(\mathcal E^{low}_{\kappa-2}(t)\right)^{\frac12}.
			\end{aligned}
		\end{equation}
		
		For the last term, we use \eqref{equ_F11_1.3} to obtain
		
		\begin{equation}\label{equ_F12_2}
			\begin{aligned}
				\langle t\rangle^{-\delta}& \sum_{|\alpha|\leq\kappa-2} \int_0^t\int_{\mathbb R^3} (\tau+\sigma^{-10}) \Gamma^\alpha\Theta\cdot F_{2,\alpha} \,dx\,d\tau 
				\lesssim{} \langle t\rangle^{-\delta} \sum_{|\alpha|\leq\kappa-2} \left\| \langle\tau\rangle^{\frac12}\nabla\Gamma^\alpha\Theta \right\|_{L_\tau^2 L_x^2}\\
				\times{}&
				\Bigg[ \sum_{|\beta|+|\gamma|\leq\kappa-2}\left(\left\| \langle\tau\rangle^{\frac12} \Gamma^\beta\Theta \nabla\cdot\partial_\tau \widetilde\Gamma^\gamma u^{cf} \right\|_{L_\tau^2 L_x^{6/5}}
				+ \left\| \langle\tau\rangle^{\frac12}\nabla\Gamma^\beta\Theta \Gamma^\gamma\Theta\right\|_{L_\tau^2 L_x^2}\right) 
				+\sum_{|\beta|\leq\kappa-3} \left\| \langle\tau\rangle^{\frac12}\nabla\Gamma^\beta\Theta \right\|_{L_\tau^2 L_x^2}\Bigg]\\
				\lesssim{}&  \mathcal E_\kappa(t) \left(\mathcal E^{low}_{\kappa-2}(t)\right)^{\frac12} + \widetilde{\mathcal I}_{\kappa}^{\frac12}(t) \widetilde{\mathcal I}_{\kappa-1}^{\frac12}(t),
			\end{aligned}
		\end{equation}
		where
		\begin{equation}\label{equ_4.2}
			\begin{aligned}
				&\sum_{|\beta|+|\gamma|\leq\kappa-2} \langle t\rangle^{-\frac\delta2} \left\| \langle\tau\rangle^\frac12 \Gamma^\beta\Theta \nabla\cdot\partial_\tau \widetilde\Gamma^\gamma u^{cf} \right\|_{L_\tau^2 L_x^{6/5}} \lesssim\sum_{\substack{ |\gamma|\leq\kappa-2\\|\beta|\leq\kappa-5}} 
				\langle t\rangle^{-\frac\delta2} \left\| \Gamma^\beta\Theta \right\|_{L_\tau^\infty L_x^3} \left\| \langle\tau\rangle^{\frac12} \nabla\partial \widetilde\Gamma^\gamma u^{cf} \right\|_{L_\tau^2 L_x^2} \\
				&\qquad+ \sum_{\substack{ |\gamma|+2\leq\kappa-5\\|\beta|\leq\kappa-2}} 
				\left\| \langle\tau\rangle^{-\frac\delta2} \Gamma^\beta\Theta \right\|_{L_\tau^\infty L_x^3} \left\| \langle\tau\rangle^{\frac12} \nabla\partial \widetilde\Gamma^\gamma u^{cf} \right\|_{L_\tau^2 L_x^2} 
				\lesssim \mathcal E_\kappa^{\frac12}(t) \left(\mathcal E^{low}_{\kappa-2}(t)\right)^{\frac12}.
			\end{aligned}
		\end{equation}
		
		Combining \eqref{equ_F11_1-1}, \eqref{equ_F11_2}, \eqref{equ_F11_3} - \eqref{equ_F12_2} completes the proof.
		
	\end{proof}
	
	\begin{lemma}\label{lem_i_1}
		Let $\kappa\ge 9$, 
		and let $(u,\Theta)\in \dot H_{\widetilde\Gamma}^{\kappa}(T) \times\dot H_{\Gamma}^{\kappa}(T)$ 
		be a 
		solution of \eqref{equ_1}. 
		Then
		\begin{equation}\nonumber
			X_{1,\kappa}(t)\lesssim \mathcal E_\kappa(t) \left(\mathcal E^{low}_{\kappa-2}(t)\right)^{\frac12} + \widetilde{\mathcal I}_{\kappa}^{\frac12}(t) \widetilde{\mathcal I}_{\kappa-1}^{\frac12}(t)+X_{1,\kappa}(0).
		\end{equation}
	\end{lemma}
	\begin{proof}
		Substituting
		$(w,\Sigma)=(\nabla^{-1}\widetilde{\Gamma}^\alpha u^{cf},\nabla^{-1}\Gamma^\alpha\Theta)$,
		$F_1=\nabla^{-1}F_{1,\alpha}$, and $F_2=\nabla^{-1}F_{2,\alpha}$
		(see \eqref{equ_F1_F2}) into \eqref{equ_cf_01}, we obtain
		\begin{equation}\nonumber
			X_{1,\kappa}(t)\lesssim{} \langle t\rangle^{-\delta} \sum_{|\alpha|\leq\kappa-2} \int_0^t K_{4,\alpha}(\tau)d\tau +\sigma^{-4}X_{1,\kappa}(0),
		\end{equation}
		where
		\begin{equation}\label{equ_K_4}
			\begin{aligned}
				K_{4,\alpha}(\tau)& ={}\int_{\mathbb R^3} \Bigl[ A_1\Big( \partial_\tau\widetilde\Gamma^\alpha u^{cf} \cdot F_{1,\alpha} +\nabla\nabla^{-1}\Gamma^\alpha\Theta :\nabla\nabla^{-1}F_{2,\alpha}+\partial_\tau\nabla^{-1}\widetilde\Gamma^\alpha u^{cf} :\nabla^{-1}F_{1,\alpha}\\
				+&\nabla^{-1}\Gamma^\alpha\Theta \cdot\nabla^{-1}F_{2,\alpha} \Big) +A_2\Big( \nabla^{-1}\Gamma^\alpha\Theta\cdot F_{1,\alpha} +\partial_\tau\widetilde\Gamma^\alpha u^{cf} \cdot\nabla^{-1} F_{2,\alpha} \Big) +A_3\widetilde\Gamma^\alpha u^{cf}\cdot F_{1,\alpha} \Bigr]dx.
			\end{aligned}
		\end{equation}
		
		We estimate the terms on the right-hand side separately.
		For the nonlinear terms involving $F_{2,\alpha}$, combining \eqref{equ_4.2} and \eqref{equ_F11_1.3}, we obtain
		\begin{equation}\nonumber
			\begin{aligned}
				\langle t\rangle^{-\delta}&\sum_{|\alpha|\leq\kappa-2} \int_0^t\!\int_{\mathbb R^3} \Bigl( A_1\big(\nabla\nabla^{-1}\Gamma^\alpha\Theta \cdot\nabla\nabla^{-1}F_{2,\alpha}+\nabla^{-1}\Gamma^\alpha\Theta \cdot\nabla^{-1} F_{2,\alpha}\big) +A_2\partial_\tau\widetilde\Gamma^\alpha u^{cf} \cdot\nabla^{-1}F_{2,\alpha} \Bigr)\,dx\,d\tau \\
				\lesssim{}& \sum_{|\alpha|\leq\kappa-2} \langle t\rangle^{-\delta} \Bigl( \|\langle \tau\rangle^{-\frac\delta2}\nabla\Gamma^\alpha\Theta\|_{L_\tau^\infty L_x^2} +\|\langle \tau\rangle^{-\frac\delta2}\partial_\tau\widetilde\Gamma^\alpha u^{cf}\|_{L_\tau^\infty L_x^2} +\|\langle \tau\rangle^{-\frac\delta2}\nabla^{-1}\Gamma^\alpha\Theta\|_{L_\tau^\infty L_x^2} \Bigr) \\
				&\quad\times\sum_{|\beta|+|\gamma|\leq\kappa-2} \Biggl( \left\|\langle \tau \rangle^{\frac12} \Gamma^\beta\Theta \nabla\cdot\partial_\tau\widetilde{\Gamma}^\gamma u^{cf} \right\|_{L_\tau^2 L_x^{6/5}}
				+\left\| \langle \tau \rangle^{\frac12} \nabla\Gamma^\beta\Theta\,\Gamma^\gamma\Theta \right\|_{L_\tau^2 L_x^2}\Biggr) \left\| \langle \tau \rangle^{-\frac12+\frac\delta2}\right\| _{L_\tau^2}  \\
				&+ \sum_{|\alpha|\leq\kappa-2}\langle t\rangle ^{-\delta} \Bigl( \|\langle \tau \rangle^\frac12 \nabla^2\Gamma^\alpha\Theta\|_{L_\tau^2 L_x^2} +\|\langle \tau \rangle^\frac12 \partial_\tau\nabla\widetilde\Gamma^\alpha u^{cf}\|_{L_\tau^2 L_x^2} +\|\Gamma^\alpha\Theta\|_{L_\tau^2 L_x^2} \Bigr) \sum_{|\beta|\leq\kappa-3}\|\Gamma^\beta\Theta\|_{L_\tau^2 L_x^2} \\
				\lesssim{}&\mathcal E_\kappa(t) \left(\mathcal E^{low}_{\kappa-2}(t)\right)^{\frac12} + \widetilde{\mathcal I}_{\kappa}^{\frac12}(t) \widetilde{\mathcal I}_{\kappa-1}^{\frac12}(t).
			\end{aligned}
		\end{equation}
		
		We next estimate the nonlinear terms involving $F_{1,\alpha}$.Using its divergence structure to integrate by parts in space. 
		Applying H\"older's inequality together with \eqref{equ_hig_F1_1}-\eqref{equ_hig_F1_3}, we obtain
		
		\begin{equation}\nonumber
			\begin{aligned}
				\langle t\rangle^{-\delta}& \sum_{|\alpha|\leq\kappa-2} \int_0^t\!\int_{\mathbb R^3} \biggl[ B_1\left( \partial_\tau\widetilde\Gamma^\alpha u^{cf}\cdot F_{1,\alpha} + \partial_\tau\nabla^{-1}\widetilde\Gamma^\alpha u^{cf} :\nabla^{-1}F_{1,\alpha} \right)+B_2\nabla^{-1}\Gamma^\alpha\Theta\cdot F_{1,\alpha} +B_3\widetilde\Gamma^\alpha u^{cf}\cdot F_{1,\alpha} \biggr] \\
				\lesssim{}& -\langle t\rangle^{-\delta} \sum_{\substack{|\alpha|\leq\kappa-2\\ \beta+\gamma=\alpha}} \int_0^t\!\int_{\mathbb R^3} \biggl( B_1\partial_\tau\partial_l \widetilde\Gamma^\alpha u_i^{cf} +B_1\partial_\tau\partial_l(-\Delta)^{-1} \widetilde\Gamma^\alpha u_i^{cf} +B_2\partial_l\partial_i(-\Delta)^{-1}\Gamma^\alpha\Theta +B_3\partial_l\widetilde\Gamma^\alpha u_i^{cf} \biggr) \\
				&\hspace{5em}\times \biggl( g_{lmn}^{ijk} \partial_m\widetilde\Gamma^\beta u_j\, \partial_n\widetilde\Gamma^\gamma u_k + g_{lm}^{ij} \partial_m\widetilde\Gamma^\beta u_j\, \Gamma^\gamma\Theta \biggr)\,dx\,d\tau \\
				& +\langle t\rangle^{-\delta} \sum_{\substack{|\alpha|\leq\kappa-2\\ \beta+\gamma\le \alpha}} \int_0^t\!\int_{\mathbb R^3}\biggl( \left|\partial^{\le 1}\nabla\widetilde\Gamma^\alpha u^{cf}\right|+\left|\partial_\tau\nabla^{-1} \widetilde\Gamma^\alpha u^{cf}\right|+\left|\Gamma^\alpha\Theta\right| \biggr) \left|\Gamma^\beta\Theta\,\Gamma^\gamma\Theta\right|dxd\tau \\
				\lesssim{}& \langle t\rangle^{-\delta} \sum_{\substack{|\alpha|\leq\kappa-2\\
						|\beta|+|\gamma|\leq\kappa-2}} \biggl( \bigl\|\langle\tau\rangle^{-\frac\delta2} \partial^{\le1}\nabla\widetilde\Gamma^\alpha u^{cf} \bigr\|_{L_\tau^\infty L_x^2} + \bigl\|\langle\tau\rangle^{-\frac\delta2} \partial_\tau \nabla^{-1}\widetilde\Gamma^\alpha u^{cf} \bigr\|_{L_\tau^\infty L_x^2} + \bigl\|\langle\tau\rangle^{-\frac\delta2} \Gamma^\alpha\Theta\bigr\|_{L_\tau^\infty L_x^2}\biggr)\\
				&\hspace{3em}\times \biggl( \bigl\|\langle\tau\rangle^{\frac12} \nabla\widetilde\Gamma^\beta u\, \nabla\widetilde\Gamma^\gamma u \bigr\|_{L_\tau^2 L_x^2} + \bigl\|\langle\tau\rangle^{\frac12} \nabla\widetilde\Gamma^\beta u\, \Gamma^\gamma\Theta \bigr\|_{L_\tau^2 L_x^2}+ \bigl\|\langle\tau\rangle^{\frac12} \Gamma^\beta\Theta\,\Gamma^\gamma\Theta \bigr\|_{L_\tau^2 L_x^2} \biggr)\\
				\lesssim{}&\mathcal E_\kappa(t)\left(\mathcal E^{low}_{\kappa-2}(t)\right)^{\frac12}.
			\end{aligned}
		\end{equation}

		Combining the above estimates completes the proof.
	\end{proof}
	
	\subsection{Lower-order weighted energy estimates}
	
	In this subsection, we estimate the lower-order energies associated with $u^{cf}$ and $\Theta$.
	
	\begin{lemma}\label{lem_J_2}
		Let $\kappa\ge 9$, 
		and let $(u,\Theta)\in \dot H_{\widetilde\Gamma}^{\kappa}(T) \times\dot H_{\Gamma}^{\kappa}(T)$ 
		be a solution of \eqref{equ_1}. 
		Then
		\begin{equation}\nonumber
			\begin{aligned}
				\mathcal J_{2,\kappa-2}(t)\lesssim{}& \int_0^t \langle\tau\rangle^{-\frac32+\frac{\delta}{2}} \widetilde{\mathcal J}_{\kappa-2}(\tau )\widetilde{\mathcal I}_\kappa^{\frac12}(\tau)\,d\tau+ \left(\mathcal E^{low}_{\kappa-2}(t)\right)^{\frac32}
				+\widetilde{\mathcal J}_{\kappa-2}^\frac12(t)\widetilde{\mathcal J}_{\kappa-3}^\frac12(t)+\mathcal J_{1,\kappa-2}(t)+\mathcal J_{2,\kappa-2}(0).
			\end{aligned}
		\end{equation}
	\end{lemma}
	\begin{proof}
		Following the argument in the proof of Lemma \ref{lem_i_2}, we obtain
		
		\begin{equation}\label{equ_j2}
			\begin{aligned}
				\mathcal J_{2,\kappa-2}(t)\lesssim{}& \sum_{|\alpha|\leq\kappa-4}\int_0^t (\tau+\sigma^{-10}) K_{3,\alpha}(\tau) d\tau+\sigma^{-14}\mathcal J_{2,\kappa-2}(0)+\sigma^{-4}\mathcal J_{1,\kappa-2}(t).
			\end{aligned}
		\end{equation}
		Here, $K_{3,\alpha}(\tau)$ is defined in \eqref{equ_K_3}.
		Unlike the estimate for $X_{2,\kappa}(t)$, the present estimate lacks the decay factor $\langle t\rangle^{-\delta}$.
		Consequently, the nonlinear terms cannot all be controlled solely by the lower-order weighted energy. For certain terms, we also need the higher-order weighted energy, together with the time decay provided by the different null structure. 
		We now estimate the terms on the right-hand side of \eqref{equ_j2} separately.
		
		We begin with the term involving $F_{1,\alpha}$. 
		Using the expression for $F_{1,\alpha}$ and applying H\"older's inequality, we obtain
		\begin{equation}\label{eq:F1-alpha-estimate}
			\begin{aligned}
				&\sum_{|\alpha|\leq\kappa-4} \int_0^t\int_{\mathbb R^3}(\tau+\sigma^{-10}) \partial_\tau\nabla\cdot\widetilde\Gamma^\alpha u^{cf}\, \nabla\cdot F_{1,\alpha}\,dx\,d\tau \\
				&\quad\lesssim{}\mathcal A_{1}(u,u)+ \sum_{\substack{|\alpha|\leq\kappa-4\\ \beta+\gamma=\alpha}} \bigl\| \langle\tau\rangle^{\frac12} \partial\nabla \widetilde\Gamma^\alpha u^{cf} \bigr\|_{L_\tau^2 L_x^2}
				\left[ \bigl\| \langle\tau\rangle^{\frac12} \nabla^2(\nabla\widetilde\Gamma^\beta u\Gamma^\gamma\Theta) \bigr\|_{L_\tau^2 L_x^2}+\left\|\langle \tau \rangle^\frac12 \nabla^2( \Gamma^\beta \Theta\Gamma^\gamma\Theta) \right\| _{L_\tau^2 L_x^2}\right].
			\end{aligned}
		\end{equation}
		where $\mathcal A_1(u,u)$ denotes the contribution from the nonlinear elastic interaction and is defined by
		\begin{equation}\nonumber
			\begin{aligned}
				&\mathcal A_{1}(u,u)\triangleq\sum_{\substack{|\alpha|\leq\kappa-4\\ \beta+\gamma=\alpha}}\sum_{a=1}^3\int_0^t\!\int_{\mathbb R^3}(\tau+\sigma^{-10}) \partial_\tau\partial_a\widetilde\Gamma^\alpha u_i^{cf}\, \partial_a\partial_l \bigl( g^{ijk}_{lmn} \partial_m\widetilde\Gamma^\beta u_j \partial_n\widetilde\Gamma^\gamma u_k \bigr)\,dx\,d\tau.
			\end{aligned}
		\end{equation}
		We first apply the Helmholtz decomposition to $\mathcal A_{1}(u,u)$.	For $\mathcal A_{1}(u^{df},u^{df})$, we integrate by parts in space
		and use Lemma \ref{lem_null} to exploit the null structure. We treat the regions $r\ge \frac{\langle c_2\tau \rangle}{2}$ and $r\le \frac{\langle c_2\tau \rangle}{2}$ separately, i.e.
		\begin{equation}\label{equ_df}
			\begin{aligned}
				&\sum_{|\beta|+|\gamma|\le \kappa-4} \left\|\langle \tau \rangle^{\frac32}N_2(\nabla^{\le 1}\widetilde{\Gamma}^\beta u^{df},\widetilde{\Gamma}^\gamma u^{df})\right\|_{L_x^2}\\
				&\quad\lesssim{}\sum_{\substack{|\gamma|\le \kappa-7\\ |\beta|\le \kappa-4}}\Big(\left\|\langle r\rangle^{\frac12}\nabla^{\le 1}\widetilde{\Gamma}^{\beta+1}u^{df}\right\|_{L_x^{\infty,2}}\left\|\nabla^2\nabla^{\le 1}\widetilde{\Gamma}^\gamma u^{df}\right\|_{L_x^{2,\infty}}
				+\left\|\nabla\nabla^{\le 1}\widetilde{\Gamma}^{\beta+1}u^{df}\right\|_{L_x^2}\left\|\langle r\rangle^\frac12\nabla\nabla^{\le 1}\widetilde{\Gamma}^\gamma u^{df}\right\|_{L_x^\infty}\Big)\\
				&\qquad+\sum_{\substack{|\beta|\le \kappa-7\\ |\gamma|\le \kappa-4}}\Big(\left\|\langle r\rangle^{\frac12}\nabla^{\le 1}\widetilde{\Gamma}^{\beta+1}u^{df}\right\|_{L_x^\infty}\left\|\nabla^2\nabla^{\le 1}\widetilde{\Gamma}^\gamma u^{df}\right\|_{L_x^2}
				+\left\|\langle r\rangle^{\frac12}\nabla\nabla^{\le 1}\widetilde{\Gamma}^{\beta+1}u^{df}\right\|_{L_x^\infty}\left\|\nabla\nabla^{\le 1}\widetilde{\Gamma}^\gamma u^{df}\right\|_{L_x^2}\Big)\\
				&\qquad+\Big(\sum_{\substack{|\alpha|,|\beta|\leq\kappa-4 \\ |\gamma|\leq\kappa-7}}+\sum_{\substack{|\alpha|,|\gamma|\leq\kappa-4\\ |\beta|\leq\kappa-7}}\Big)
				\left\|\langle c_2\tau -r\rangle \nabla^2\nabla^{\le 1}\widetilde{\Gamma}^\beta u^{df}\right\|_{L_x^2} \left\|\langle c_2\tau -r\rangle \nabla\nabla^{\le 1}\widetilde{\Gamma}^\gamma u^{df}\right\|_{L_x^\infty}\\
				&\quad\lesssim{}\langle \tau \rangle^{-\frac12+\frac\delta2}\widetilde{\mathcal I}_\kappa^\frac12(\tau)\widetilde{\mathcal J}_{\kappa-2}^\frac12(\tau)+\mathcal E_{\kappa-2}^{low}(\tau).
			\end{aligned}
		\end{equation}
		Hence,
		\begin{equation}\label{eq:A61-df-df}
			\begin{aligned}
				\mathcal A_{1}&(u^{df},u^{df})  \lesssim\int_0^t \langle\tau\rangle^{-\frac32+\frac{\delta}{2}} \widetilde{\mathcal J}_{\kappa-2}(\tau )\widetilde{\mathcal I}_\kappa^{\frac12}(\tau)\,d\tau +\left(\mathcal E_{\kappa-2}^{low}(t)\right)^\frac32.
			\end{aligned}
		\end{equation}
		
		By the symmetry of $g^{ijk}_{lmn}$, the remaining contributions	to $\mathcal A_{1}(u,u)$ reduce to terms of the form $\mathcal A_{1}(u^{cf},u)$. 
		A direct estimate of the term in which all the commutation derivatives fall on $u^{cf}$ would require more derivatives than are controlled by the lower-order energy.
		To avoid this loss of derivatives, we isolate the highest-order term before estimating the remaining terms.
		We write
		\begin{equation}\nonumber
			\begin{aligned}
				\mathcal A_{1}(u^{cf},u) ={}& -\sum_{|\alpha|\leq\kappa-4}\sum_{a=1}^{3} \int_0^t\!\int_{\mathbb R^3} (\tau+\sigma^{-10}) \partial_\tau\partial_l\partial_a \widetilde\Gamma^\alpha u_i^{cf}\, \partial_a\partial_m \widetilde\Gamma^\alpha u_j^{cf}\, \partial_nu_k\,g^{ijk}_{lmn} \,dx\,d\tau \\
				&+\sum_{|\alpha|\leq\kappa-4}\sum_{a=1}^{3} \int_0^t\!\int_{\mathbb R^3}(\tau+\sigma^{-10}) \partial_\tau\partial_a \widetilde\Gamma^\alpha u_i^{cf}\, \partial_l(g^{ijk}_{lmn}\partial_m \widetilde\Gamma^\alpha u_j^{cf}\, \partial_n\partial_a u_k)dx\,d\tau \\
				&+ \sum_{\substack{|\beta|< |\alpha|\leq\kappa-4\\\beta+\gamma=\alpha,}} 
				\sum_{a=1}^{3} \int_0^t\!\int_{\mathbb R^3} (\tau+\sigma^{-10}) \partial_\tau\partial_a \widetilde\Gamma^\alpha u_i^{cf}\, \partial_a\partial_l \bigl( g^{ijk}_{lmn} \partial_m\widetilde\Gamma^\beta u_j^{cf} \partial_n\widetilde\Gamma^\gamma u_k \bigr)\,dx\,d\tau .
			\end{aligned}
		\end{equation}
		For the highest-order term, we use the symmetry of the
		coefficients and integrate by parts in time. Thus, we obtain
		\begin{equation}\label{eq:A61-cf-u-estimate}
			\begin{aligned}
				\mathcal A_{1}(u^{cf},u) \lesssim{}& -\sum_{|\alpha|\leq\kappa-4}\sum_{a=1}^{3} \left. \int_{\mathbb R^3} \frac12 g^{ijk}_{lmn}(\tau+\sigma^{-10}) \partial_a\partial_l \widetilde\Gamma^\alpha u_i^{cf}\, \partial_a\partial_m \widetilde\Gamma^\alpha u_j^{cf}\, \partial_nu_k \,dx \right|_{\tau=0}^{\tau=t} \\
				&+ \sum_{|\alpha|\leq\kappa-4} \int_0^t\!\int_{\mathbb R^3}\frac12 g^{ijk}_{lmn} \partial_a\partial_l \widetilde\Gamma^\alpha u_i^{cf}\, \partial_a\partial_m \widetilde\Gamma^\alpha u_j^{cf} \Biggl( \partial_nu_k +(\tau+\sigma^{-10}) \partial_\tau\partial_nu_k \Biggr)\,dx\,d\tau \\
			\end{aligned}
		\end{equation}
		\begin{equation}\nonumber
			\begin{aligned}
				&+ \sum_{|\alpha|\leq\kappa-4} \bigl\| \langle\tau\rangle^{\frac12} \partial\nabla \widetilde\Gamma^\alpha u^{cf} \bigr\|_{L_\tau^2 L_x^2} \Big(\bigl\| \langle\tau\rangle^{\frac12} \nabla\widetilde\Gamma^\alpha u^{cf} \bigr\|_{L_\tau^2 L_x^6} \bigl\|\nabla^3u\bigr\|_{L_\tau^\infty L_x^3} +
				\bigl\| \langle\tau\rangle^{\frac12} \nabla^2\widetilde\Gamma^\alpha u^{cf} \bigr\|_{L_\tau^2 L_x^2} \bigl\|\nabla^2u\bigr\|_{L_\tau^\infty L_x^\infty}\Big) \\
				&+ \sum_{\substack{\beta+\gamma=\alpha\\|\beta|<|\alpha|\leq\kappa-4}} \bigl\| \langle\tau\rangle^\frac12 \partial_\tau\nabla \widetilde\Gamma^\alpha u^{cf} \bigr\|_{L_\tau^2 L_x^2} \Bigl( \bigl\| \langle\tau\rangle^{\frac12} \nabla^3\widetilde\Gamma^\beta u^{cf} \bigr\|_{L_\tau^2 L_x^2} \bigl\|\nabla\widetilde\Gamma^\gamma u\bigr\|_{L_\tau^\infty L_x^\infty} \\
				&\qquad +\bigl\| \langle\tau\rangle^{\frac12} \nabla^2\widetilde\Gamma^\beta u^{cf} \bigr\|_{L_\tau^2 L_x^6} \bigl\|\nabla^2\widetilde\Gamma^\gamma u\bigr\|_{L_\tau^\infty L_x^3} +\bigl\| \langle\tau\rangle^{\frac12} \nabla\widetilde\Gamma^\beta u^{cf} \bigr\|_{L_\tau^2 L_x^\infty} \bigl\|\nabla^3\widetilde\Gamma^\gamma u\bigr\|_{L_\tau^\infty L_x^2} \Bigr) \\
				\lesssim{}& \left(\mathcal E^{low}_{\kappa-2}(t)\right)^{\frac32}.
			\end{aligned}
		\end{equation}
		The mixed terms can be bounded as follows
		\begin{equation}\label{eq:A63-u-Theta-estimate}
			\begin{aligned}
				\sum_{\substack{|\beta|+|\gamma|\le \kappa-4\\ |a|\le 2}}& \bigl\| \langle\tau\rangle^{\frac12} \nabla^a(\nabla\widetilde\Gamma^\beta u\Gamma^\gamma\Theta) \bigr\|_{L_\tau^2 L_x^2}\lesssim{} \sum_{|\gamma|,|\beta|\le \kappa-4}\Bigg(\left\| \langle \tau \rangle ^{\frac12}\Gamma^\gamma \Theta\right\| _{L_\tau^2 L_x^\infty}\left\|\nabla \nabla ^{\le 2}\widetilde{\Gamma}^\beta u\right\| _{L_\tau^\infty L_x^2}\\
				&+\left\| \langle \tau \rangle^\frac12 \nabla \Gamma^\gamma \Theta\right\| _{L_\tau^2 L_x^6}\left\| \nabla\nabla^{\le 1}\widetilde{\Gamma}^\beta u\right\| _{L_\tau^\infty L_x^3} +\left\| \langle \tau \rangle^\frac12 \nabla^2 \Gamma^\gamma \Theta\right\| _{L_\tau^2 L_x^2}\left\| \nabla\widetilde{\Gamma}^\beta u\right\| _{L_\tau^\infty L_x^\infty} \Bigg) \\
				\lesssim{}&\mathcal E^{low}_{\kappa-2}(t).
			\end{aligned}
		\end{equation}
		Finally, we estimate the term involving only $\Theta$.
		\begin{equation}\label{eq:F1-alpha-theta-estimate}
			\begin{aligned}
				\sum_{\substack{|\beta|+|\gamma|\le \kappa-4\\ |a|\le 2}}& \left\|\langle \tau \rangle^\frac12 \nabla^a( \Gamma^\beta \Theta\Gamma^\gamma\Theta) \right\| _{L_\tau^2 L_x^2} 
				\lesssim{}\sum_{\substack{|\beta|+|\gamma|\le \kappa-3}} \left\|\langle \tau \rangle^\frac12 \nabla^{\le 2}( \Gamma^\beta \Theta\Gamma^\gamma\Theta) \right\| _{L_\tau^2 L_x^2} \\
				\lesssim{}& \sum_{\substack{|\gamma|\le \kappa-3\\ |\beta|\le \kappa-6}}\Bigg(\left\| \langle \tau \rangle ^{\frac12}\nabla^{\le 2}\Gamma^\beta \Theta\right\| _{L_\tau^\infty L_x^3}\left\|\nabla ^{\le 1}\Gamma^\gamma \Theta\right\| _{L_\tau^2 L_x^6}+ \left\| \langle \tau \rangle ^{\frac12}\Gamma^\beta \Theta\right\| _{L_\tau^\infty L_x^\infty}\left\|\nabla ^{2}\Gamma^\gamma \Theta\right\| _{L_\tau^2 L_x^2}\Bigg) \\
				\lesssim{}&\mathcal E^{low}_{\kappa-2}(t).
			\end{aligned}
		\end{equation}
		Combining \eqref{eq:A61-df-df},
		\eqref{eq:A61-cf-u-estimate}-
		\eqref{eq:F1-alpha-theta-estimate}, and substituting
		the resulting bounds into \eqref{eq:F1-alpha-estimate},
		we obtain the following estimate
		\begin{equation}\label{equ_l_F12_1}
			\begin{aligned}
				\sum_{|\alpha|\leq\kappa-4}& \int_0^t\int_{\mathbb R^3}(\tau+\sigma^{-10})\partial_\tau\nabla\widetilde\Gamma^\alpha u^{cf} \cdot\nabla F_{1,\alpha} dx d\tau \lesssim{} \left(\mathcal E^{low}_{\kappa-2}(t)\right)^{\frac32} +\int_0^t \langle\tau\rangle^{-\frac32+\frac{\delta}{2}} \widetilde{\mathcal J}_{\kappa-2}(\tau )\widetilde{\mathcal I}_\kappa^{\frac12}(\tau)\,d\tau .
			\end{aligned}
		\end{equation}
		
		We next estimate the $F_{1,\alpha}$ term without a spatial divergence.
		Expanding $F_{1,\alpha}$, we obtain
		\begin{equation}
			\begin{aligned}
				\sum_{|\alpha|\leq\kappa-4} &\int_0^t\int_{\mathbb R^3} (\tau+\sigma^{-10}) \partial_\tau\widetilde\Gamma^\alpha u^{cf} \cdot F_{1,\alpha} dx d\tau \\
				\lesssim{}& \sum_{\substack{|\alpha|\leq\kappa-4\\\beta+\gamma=\alpha}} \int_0^t\int_{\mathbb R^3} \langle\tau\rangle \left| g^{ijk}_{lmn}\partial_\tau \partial_l\widetilde\Gamma^\alpha u^{cf}_i \partial_m\widetilde\Gamma^\beta u_j^{df} \partial_n\widetilde\Gamma^\gamma u_k^{df}\right|dx d\tau \\
				&+ \sum_{\substack{|\alpha|\leq\kappa-4\\\beta+\gamma=\alpha}} \int_0^t\int_{\mathbb R^3} \langle\tau\rangle | \partial_\tau\partial_l\widetilde\Gamma^\alpha u_i^{cf}|\left( | g^{ijk}_{lmn} \partial_m\widetilde\Gamma^\beta u_j \partial_n\widetilde\Gamma^\gamma u_k^{cf}|+| g^{ij}_{lm} \partial_m\widetilde\Gamma^\beta u_j \Gamma^\gamma\Theta|\right) dx d\tau \\
				&+ \sum_{|\alpha|,|\beta|,|\gamma|\leq\kappa-4} \left\| \langle\tau\rangle^{\frac12} \partial\widetilde\Gamma^\alpha u^{cf} \right\|_{L_\tau^2 L_x^6} \left\| \langle\tau\rangle^{\frac12} \nabla\Gamma^\beta\Theta \right\|_{L_\tau^2 L_x^2} \left\| \Gamma^\gamma\Theta \right\|_{L_\tau^\infty L_x^3}\\
				\lesssim{}& \left(\mathcal E^{low}_{\kappa-2}(t)\right)^{\frac32} +\int_0^t \langle\tau\rangle^{-\frac32+\frac{\delta}{2}} \widetilde{\mathcal J}_{\kappa-2}(\tau )\widetilde{\mathcal I}_\kappa^{\frac12}(\tau)\,d\tau .
			\end{aligned}
		\end{equation}
		By Lemma \ref{lem_null_2}, we have
		\begin{equation}\label{equ_nu_df_df}
			\begin{aligned}
				&\sum_{|\beta|+|\gamma|\le \kappa-4}\int_{\mathbb R^3}\langle \tau\rangle^{\frac32}\left| g_{lmn}^{ijk} \partial_l w_i \partial_m\widetilde\Gamma^\beta u_j^{df} \partial_n\widetilde\Gamma^\gamma u_k^{df}\right|dx
				\triangleq \sum_{|\beta|+|\gamma|\le \kappa-4}\int_{\mathbb R^3}\langle \tau\rangle^{\frac32}N(w,\widetilde{\Gamma}^\beta u^{df},\widetilde{\Gamma}^\gamma u^{df})dx\\
				&\quad\lesssim{}\|\nabla w\|_{L_x^2}\Big(\sum_{|\beta|,|\gamma|\le\kappa-4}
				\left\|\langle c_2\tau-r\rangle\nabla\widetilde\Gamma^\beta u^{df}\right\|_{L_x^6}
				\left\|\langle c_2\tau-r\rangle^\frac12\nabla\widetilde\Gamma^\gamma u^{df}\right\|_{L_x^3}\\
				&\qquad+\sum_{\substack{|\gamma|\le\kappa-7\\|\beta|\le\kappa-4}}
				\left\|\langle r\rangle^{\frac12}\widetilde\Gamma^{\beta+1} u^{df}\right\|_{L_x^{\infty,2}}
				\left\|\nabla\widetilde\Gamma^\gamma u^{df}\right\|_{L_x^{2,\infty}}
				+\sum_{\substack{|\gamma|\le\kappa-4\\|\beta|\le\kappa-7}}
				\left\|\langle r\rangle^{\frac12}\widetilde\Gamma^{\beta+1} u^{df}\right\|_{L_x^\infty}
				\left\|\nabla\widetilde\Gamma^\gamma u^{df}\right\|_{L_x^2}\Big)\\
				&\qquad+\langle \tau \rangle^{-\frac12}\sum_{\substack{|\beta|\le \kappa-7\\ |\gamma|\le \kappa-4}}\left\| \widetilde{\Omega}^{\le 1}w\right\|_{L_x^2}\left\|\langle \tau \rangle\nabla\widetilde{\Gamma}^\beta u^{df}\right\|_{L_x^\infty}\left\|\nabla\widetilde{\Gamma}^\gamma u^{df}\right\|_{L_x^2}\\
				&\quad\lesssim{}\left(\|\nabla w(\tau)\|_{L_x^2}+\langle \tau\rangle^{-\frac12}\left\| \widetilde{\Omega}^{\le 1}w\right\|_{L_x^2}\right)\widetilde{\mathcal J}_{\kappa-2}(\tau).
			\end{aligned}
		\end{equation}
		Setting $w=\partial_\tau\widetilde{\Gamma}^\alpha u^{cf}$ in the preceding estimate yields
		\begin{equation}
			\sum_{\substack{|\alpha|\leq\kappa-4\\\beta+\gamma=\alpha}} \int_0^t\int_{\mathbb R^3} \langle\tau\rangle \left| g^{ijk}_{lmn}\partial_\tau \partial_l\widetilde\Gamma^\alpha u^{cf}_i \partial_m\widetilde\Gamma^\beta u_j^{df} \partial_n\widetilde\Gamma^\gamma u_k^{df}\right|dx d\tau 
			\lesssim{}\int_0^t \langle\tau\rangle^{-\frac32+\frac{\delta}{2}} \widetilde{\mathcal J}_{\kappa-2}(\tau )\widetilde{\mathcal I}_\kappa^{\frac12}(\tau)\,d\tau+\left(\mathcal E_{\kappa-2}^{low}(t)\right)^\frac32.
		\end{equation}
		For the remaining products, H\"older's inequality gives
		\begin{equation}\nonumber
			\begin{aligned}
				\sum_{|\beta|,|\gamma|\leq\kappa-4}&\left( \left\| \langle \tau \rangle^\frac12 g^{ijk}_{lmn} \partial_m \widetilde{\Gamma}^\beta u_j \partial_n\widetilde\Gamma^\gamma u_k^{cf}\right\| _{L_\tau^2 L_x^2}+\left\| \langle \tau \rangle^\frac12 g^{ij}_{lm} \partial_m \widetilde{\Gamma}^\beta u_j \Gamma^\gamma\Theta\right\| _{L_\tau^2 L_x^2}\right) \\
				\lesssim {}&\sum_{|\beta|,|\gamma|\le \kappa-4}\left\| \nabla \widetilde{\Gamma}^\beta u\right\| _{L_\tau^\infty L_x^3}\left(\left\| \langle \tau \rangle^\frac12 \nabla \widetilde{\Gamma}^\gamma u^{cf}\right\| _{L_\tau^2 L_x^6} +\left\| \langle \tau \rangle^\frac12 \Gamma^\gamma\Theta\right\| _{L_\tau^2 L_x^6}\right)\\
				\lesssim{}& \mathcal E^{low}_{\kappa-2}(t).
			\end{aligned}
		\end{equation}
		Combining the above estimates, we obtain
		\begin{equation}\nonumber
			\begin{aligned}
				\sum_{|\alpha|\leq\kappa-4} &\int_0^t\int_{\mathbb R^3} (\tau+\sigma^{-10}) \partial_\tau\widetilde\Gamma^\alpha u^{cf} \cdot F_{1,\alpha} dx d\tau 
				\lesssim{} \left(\mathcal E^{low}_{\kappa-2}(t)\right)^{\frac32} +\int_0^t \langle\tau\rangle^{-\frac32+\frac{\delta}{2}} \widetilde{\mathcal J}_{\kappa-2}(\tau )\widetilde{\mathcal I}_\kappa^{\frac12}(\tau)\,d\tau .
			\end{aligned}
		\end{equation}
		Finally, we estimate the two remaining integral terms involving $F_{1,\alpha}$ together. 
		Integrating by parts in space and taking absolute values, we combine \eqref{eq:A63-u-Theta-estimate} and \eqref{eq:F1-alpha-theta-estimate} with the interior and exterior estimates used in \eqref{equ_df} to obtain
		\begin{equation}
			\begin{aligned}
				\sum_{|\alpha|\leq\kappa-4}&\int_0^t\int_{\mathbb R^3}(\tau+\sigma^{-10}) \left( A_2\Gamma^\alpha \Theta +A_3\nabla \cdot \widetilde{\Gamma}^\alpha u^{cf}\right) \nabla\cdot F_{1,\alpha}\,dx\,d\tau \\
				\lesssim{}& \sum_{\substack{|\alpha|\le \kappa-4\\ |\beta|+|\gamma|\le \kappa-4}}\int_0^t\int_{\mathbb R^3} \langle \tau \rangle \left|\left(A_2 \partial_i \Gamma^\alpha \Theta+ A_3\Delta\widetilde{\Gamma}^\alpha u^{cf}_i \right)\partial_l\left(g^{ijk}_{lmn}\partial_m \widetilde{\Gamma}^\beta u^{df}_j\partial_n \widetilde{\Gamma}^\gamma u^{df}_k \right)\right|\,dx\,d\tau \\
				&+ \sum_{\substack{ |\alpha|\leq\kappa-4\\
						|\beta|+|\gamma|\leq\kappa-4}} \Biggl( \left\| \langle\tau\rangle^{\frac12} \nabla\Gamma^\alpha\Theta \right\|_{L_\tau^2 L_x^2} + \left\| \langle\tau\rangle^{\frac12} \nabla^2\widetilde\Gamma^\alpha u^{cf} \right\|_{L_\tau^2 L_x^2} \Biggr) \\
				&\quad\times \Biggl(\left\| \langle\tau\rangle^{\frac12}\nabla( \nabla\widetilde\Gamma^\beta u^{cf} \nabla\widetilde{\Gamma}^\gamma u) \right\|_{L_\tau^2 L_x^2}+\left\| \langle\tau\rangle^{\frac12}\nabla( \nabla\widetilde\Gamma^\beta u \Gamma^\gamma \Theta) \right\|_{L_\tau^2 L_x^2} + \left\| \langle\tau\rangle^{\frac12} \nabla(\Gamma^\beta\Theta \Gamma^\gamma\Theta) \right\|_{L_\tau^2 L_x^2} \Biggr)\\
				&\lesssim \left(\mathcal E^{low}_{\kappa-2}(t)\right)^{\frac32}.
			\end{aligned}
		\end{equation}
		where
		\begin{equation}\nonumber
			\begin{aligned}
				\sum_{|\beta|+|\gamma|\le \kappa-4}& \left\| \langle\tau\rangle^{\frac12}\nabla( \nabla\widetilde\Gamma^\beta u^{cf} \nabla\widetilde{\Gamma}^\gamma u) \right\|_{L_\tau^2 L_x^2}\\
				\lesssim {}&\sum_{|\beta|,|\gamma|\le \kappa-4}\left( \left\| \langle \tau \rangle^\frac12 \nabla\widetilde{\Gamma}^\beta u^{cf}\right\| _{L_\tau^2 L_x^6}\left\| \nabla^2\widetilde{\Gamma}^\gamma u\right\| _{L_\tau^\infty L_x^3}+\left\| \langle \tau \rangle ^\frac12\nabla^2\widetilde{\Gamma}^\beta u^{cf}\right\| _{L_\tau^2 L_x^2}\left\| \nabla \widetilde{\Gamma}^\gamma u\right\| _{L_\tau^\infty L_x^\infty}\right) \\
				\lesssim{}&\mathcal E^{low}_{\kappa-2}(t).
			\end{aligned}
		\end{equation}
		
		We next estimate the remaining nonlinear terms. Expanding $F_{2,\alpha}$ and using \eqref{eq:F1-alpha-theta-estimate}, we obtain
		\begin{equation}\nonumber
			\begin{aligned}
				\sum_{|\alpha|\leq\kappa-4} &\int_0^t\int_{\mathbb R^3} (\tau+\sigma^{-10}) \Big( -A_1\Delta\Gamma^\alpha\Theta +A_2\partial_\tau\nabla\cdot \widetilde\Gamma^\alpha u^{cf} \Big) F_{2,\alpha} dx d\tau \\
				\lesssim{}& \sum_{|\alpha|\leq\kappa-4} \Bigg( \left\| \langle\tau\rangle^{\frac12} \Delta\Gamma^\alpha\Theta \right\|_{L_\tau^2 L_x^2} + \left\| \langle\tau\rangle^{\frac12} \partial_\tau\nabla \widetilde\Gamma^\alpha u^{cf} \right\|_{L_\tau^2 L_x^2} \Bigg)
				\Bigg[ \sum_{|\beta|\leq\kappa-5} \left\| \langle\tau\rangle^{\frac12} \Delta\Gamma^\beta\Theta \right\|_{L_\tau^2 L_x^2} \\
				&\quad+ \sum_{|\beta|+|\gamma|\leq\kappa-4} \left\| \langle\tau\rangle^{\frac12} \Gamma^\beta\Theta \nabla\cdot\partial_\tau \widetilde\Gamma^\gamma u^{cf} \right\|_{L_\tau^2 L_x^2}+ \sum_{|\beta|+|\gamma|\leq\kappa-4} \left\| \langle\tau\rangle^{\frac12} \nabla\bigl( \nabla\Gamma^\beta\Theta \Gamma^\gamma\Theta \bigr) \right\|_{L_\tau^2 L_x^2} \Bigg]\\
				\lesssim{}& \left(\mathcal E^{low}_{\kappa-2}(t)\right)^{\frac32} + \widetilde{\mathcal J}_{\kappa-2}^{\frac12}(t) \widetilde{\mathcal J}_{\kappa-3}^{\frac12}(t).
			\end{aligned}
		\end{equation}
		where
		\begin{equation}\nonumber
			\begin{aligned}
				\sum_{|\beta|+|\gamma|\leq\kappa-4}& \left\| \langle\tau\rangle^{\frac12} \Gamma^\beta\Theta \nabla\cdot\partial_\tau \widetilde\Gamma^\gamma u^{cf} \right\|_{L_\tau^2 L_x^2} \lesssim{} \sum_{|\beta|+|\gamma|\leq\kappa-4} \left\| \Gamma^\beta\Theta \right\|_{L_\tau^\infty L_x^\infty} \left\| \langle\tau\rangle^{\frac12} \nabla\cdot\partial_\tau \widetilde\Gamma^\gamma u^{cf} \right\|_{L_\tau^2 L_x^2}
				\lesssim{}\mathcal E^{low}_{\kappa-2}(t) .
			\end{aligned}
		\end{equation}
		
		Finally, combining the above estimates with \eqref{eq:F1-alpha-theta-estimate}, we obtain
		\begin{equation}\label{equ_l_F22_2}
			\begin{aligned}
				&\sum_{|\alpha|\leq\kappa-4} \int_0^t\int_{\mathbb R^3} (\tau+\sigma^{-10}) \Gamma^\alpha\Theta F_{2,\alpha} dx d\tau \\
				&\quad\lesssim{}\sum_{|\alpha|\leq\kappa-4} \left\| \langle\tau\rangle^{\frac12} \Gamma^\alpha\Theta \right\|_{L_\tau^\infty L_x^2}\sum_{|\beta|+|\mu|\leq\kappa-4} \left\| \langle\tau\rangle^{\frac12} \Gamma^\beta\Theta \nabla\partial_\tau \widetilde\Gamma^\mu u^{cf} \right\|_{L_\tau^1 L_x^2} \\
				&\qquad+ \sum_{|\alpha|\leq\kappa-4} \left\| \langle\tau\rangle^{\frac12} \nabla\Gamma^\alpha\Theta \right\|_{L_\tau^2 L_x^2} 
				\Bigg( \sum_{|\beta|+|\gamma|\leq\kappa-4} \left\| \langle\tau\rangle^{\frac12} \nabla\Gamma^\beta\Theta \Gamma^\gamma\Theta \right\|_{L_\tau^2 L_x^2} + \sum_{|\beta|\leq\kappa-5} \left\| \langle\tau\rangle^{\frac12} \nabla\Gamma^\beta\Theta \right\|_{L_\tau^2 L_x^2} \Bigg)\\
				&\quad\lesssim{} \left(\mathcal E^{low}_{\kappa-2}(t)\right)^{\frac32} + \widetilde{\mathcal J}_{\kappa-2}^{\frac12}(t) \widetilde{\mathcal J}_{\kappa-3}^{\frac12}(t).
			\end{aligned}
		\end{equation}
		Here, we have used the divergence structure of $M_1$.
		We also use the following product estimate: 
		\begin{equation}\nonumber
			\begin{aligned}
				\sum_{|\beta|+|\gamma|\leq\kappa-4}& \left\| \langle\tau\rangle^{\frac12} \Gamma^\beta\Theta \nabla\cdot\partial_\tau \widetilde\Gamma^\gamma u^{cf} \right\|_{L_\tau^1 L_x^2} \lesssim{} \sum_{\substack{ |\gamma|, |\beta|\leq\kappa-4 }} \left\| \Gamma^\beta\Theta \right\|_{L_\tau^2 L_x^\infty} \left\| \langle\tau\rangle^{\frac12}\partial\nabla \widetilde\Gamma^\gamma u^{cf} \right\|_{L_\tau^2 L_x^2}
				\lesssim{}\mathcal E^{low}_{\kappa-2}(t).
			\end{aligned}
		\end{equation}
		Finally, substituting \eqref{equ_l_F12_1}-\eqref{equ_l_F22_2} into \eqref{equ_j2} yields the desired estimate and completes the proof of the lemma.
	\end{proof}

	\begin{lemma}\label{lem_J_1}
		Let $\kappa\ge 9$,
		and let $(u,\Theta)\in \dot H_{\widetilde\Gamma}^{\kappa}(T) \times\dot H_{\Gamma}^{\kappa}(T)$ 
		be a solution of \eqref{equ_1}. 
		Then
		\begin{equation}\label{equ_J1}
			\mathcal J_{1,\kappa-2}(t)\lesssim \int_0^t \langle\tau\rangle^{-\frac32+\frac{\delta}{2}} \widetilde{\mathcal J}_{\kappa-2}(\tau )\widetilde{\mathcal I}_\kappa^{\frac12}(\tau)\,d\tau+ \left(\mathcal E^{low}_{\kappa-2}(t)\right)^{\frac32}+\widetilde{\mathcal J}_{\kappa-2}^\frac12(t)\widetilde{\mathcal J}_{\kappa-3}^\frac12(t)+\mathcal J_{1,\kappa-2}(0).
		\end{equation}
	\end{lemma}
	\begin{proof}
		By the definitions and the preceding estimates, we obtain
		\begin{equation}\nonumber
			\begin{aligned}
				\mathcal J_{1,\kappa-2}(t)  \lesssim{}& \sum_{|\alpha|\leq\kappa-4} \int_0^t K_{4,\alpha}(\tau) d\tau+\mathcal J_{1,\kappa-2}(0).
			\end{aligned}
		\end{equation}
		Here, $K_{4,\alpha}(\tau)$ is defined in \eqref{equ_K_4}.
		We estimate the integral terms on the right-hand side separately.
		For the first term involving $F_{1,\alpha}$,
		we use its divergence structure, H\"older's inequality, Lemma \ref{lem_sob_1}, \eqref{eq:A63-u-Theta-estimate}, and \eqref{eq:F1-alpha-theta-estimate} to obtain
		\begin{equation}\label{equ_low_F1_1}
			\begin{aligned}
				&\sum_{|\alpha|\leq\kappa-4} \int_0^t\int_{\mathbb R^3}\left( A_1 \partial_\tau\widetilde\Gamma^\alpha u^{cf} \cdot F_{1,\alpha}+ A_2\nabla^{-1}\Gamma^\alpha\Theta\cdot F_{1,\alpha} + A_3\widetilde\Gamma^\alpha u^{cf}\cdot F_{1,\alpha} \right) dx d\tau \\
				&\quad\lesssim{}\sum_{\substack{|\alpha|\leq\kappa-4\\|\beta|+|\gamma|\leq|\alpha|}} 
				\int_0^t \langle\tau\rangle^{-\frac32}\left( \left\| \langle\tau\rangle^{\frac12}\partial^{\le1}\nabla\widetilde\Gamma^\alpha u^{cf} \right\|_{L_x^2}+\left\| \langle \tau \rangle^\frac12\Gamma^\alpha \Theta\right\| _{L_x^2}\right)
				\left\| \langle \tau \rangle \nabla\widetilde{\Gamma}^\beta u\nabla\widetilde{\Gamma}^\gamma u\right\| _{L_x^2} \\
				&\qquad+\sum_{\substack{|\alpha|\leq\kappa-4\\|\beta|+|\gamma|\leq|\alpha|}}
				\left( \left\| \langle\tau\rangle^{\frac12}\nabla \partial\widetilde\Gamma^\alpha u^{cf} \right\|_{L_\tau^2 L_x^2}+ \left\| \langle \tau \rangle^\frac12\Gamma^\alpha \Theta\right\| _{L_\tau^\infty L_x^2}+\left\|\langle \tau \rangle^\frac12 \nabla \widetilde{\Gamma}^\alpha u^{cf}\right\|_{L_\tau^\infty L_x^2}  \right)\\
				&\qquad\times\left( \left\|\langle \tau\rangle^\frac12 \nabla\widetilde{\Gamma}^\beta u\Gamma^\gamma \Theta\right\| _{L_\tau^2 L_x^2}+\left\| \langle \tau \rangle^\frac12 \Gamma^\beta \Theta \Gamma^\gamma \Theta \right\| _{L_\tau^2 L_x^2}\right) \\
				&\quad\lesssim{}\left(\mathcal E^{low}_{\kappa-2}(t)\right)^{\frac32}.
			\end{aligned}
		\end{equation}
		Here, we use the following estimate:
		\begin{equation}
			\begin{aligned}
				\sum_{|\beta|+|\gamma|\le \kappa-4}&\langle t\rangle \left\| \nabla\widetilde{\Gamma}^\beta u\nabla\widetilde{\Gamma}^\gamma u\right\| _{L_x^2}\\
				\lesssim{}& \sum_{\substack{|\gamma|\le \kappa-7\\|\beta|\le \kappa-4}}\Big( \left\| \langle t\rangle \nabla\widetilde\Gamma^\gamma u^{df} \right\|_{L_x^\infty} \left\|\nabla\widetilde\Gamma^\beta u\right\|_{L_x^2} + \left\| \langle t\rangle^{\frac12} \nabla\widetilde\Gamma^\beta u^{cf} \right\|_{L_x^2} \left\| \langle t\rangle^{\frac12} \nabla\widetilde\Gamma^\gamma u^{cf} \right\|_{L_x^\infty}\Big) \\
				&+\sum_{\substack{|\beta|\le \kappa-7 \\|\gamma|\le \kappa-4}}\left(\left\| \langle c_2\tau -r\rangle\nabla \widetilde{\Gamma}^\gamma u^{df}\right\|_{L_x^6}\left\|\nabla \widetilde{\Gamma}^\beta u\right\|_{L_x^3}+\left\|\nabla \widetilde{\Gamma}^\gamma u^{df}\right\|_{L_x^2}\left\|\langle r\rangle\nabla \widetilde{\Gamma}^\beta u\right\|_{L_x^\infty}\right)\\
				\lesssim {}&\widetilde{\mathcal J}_{\kappa-2}(t).
			\end{aligned}
		\end{equation}
		
		Next, we consider the integral term involving $\nabla^{-1}F_{1,\alpha}$.
		This term requires a more delicate analysis.
		In particular, we separate and estimate the contribution of $N_{al}$, 
		\begin{equation}\nonumber
			\begin{aligned}
				&\sum_{|\alpha|\leq\kappa-4} \int_0^t\!\int_{\mathbb R^3} \partial_\tau\nabla^{-1}\widetilde\Gamma^\alpha u^{cf} \cdot\nabla^{-1}F_{1,\alpha}\,dx\,d\tau \\
				&\quad\lesssim \sum_{\substack{ |\alpha|\leq\kappa-4\\\beta+\gamma=\alpha }} 
				\int_0^t\!\int_{\mathbb R^3} \left| g_{lmn}^{ijk} \partial_\tau\partial_l(-\Delta)^{-1}\widetilde\Gamma^\alpha u^{cf}_i \partial_m\widetilde\Gamma^\beta u^{df}_j \partial_n\widetilde\Gamma^\gamma u^{df}_k\right|\,dx\,d\tau \\
				&\qquad+ \sum_{\substack{ |\alpha|\leq\kappa-4\\|\beta|+|\gamma|=|\alpha| }} \left\| \langle\tau\rangle^{\frac12} \partial_\tau\widetilde\Gamma^\alpha u^{cf} \right\|_{L_\tau^\infty L_x^2}\\
				&\qquad\quad \times\Bigl( \left\| \langle\tau\rangle^{-\frac12} \Gamma^\beta\Theta\,\Gamma^\gamma\Theta \right\|_{L_\tau^1 L_x^{6/5}}
				+\left\| \langle\tau\rangle^{-\frac12} \nabla\widetilde\Gamma^\beta u\nabla\widetilde{\Gamma}^\gamma u^{cf} \right\|_{L_\tau^1 L_x^{6/5}}+ \left\| \langle\tau\rangle^{-\frac12} \nabla\widetilde\Gamma^\beta u\,\Gamma^\gamma\Theta \right\|_{L_\tau^1 L_x^{6/5}} \Bigr)\\
				&\quad\lesssim \sum_{\substack{ |\alpha|\leq\kappa-4\\
						\beta+\gamma=\alpha }} \int_0^t\!\int_{\mathbb R^3}| \partial_\tau\partial_l(-\Delta)^{-1}\widetilde\Gamma^\alpha u^{cf}_i \bigl( g_{lmn}^{ijk} \partial_m\widetilde\Gamma^\beta u^{df}_j \partial_n\widetilde\Gamma^\gamma u^{df}_k \bigr)|\,dx\,d\tau+\left(\mathcal E^{low}_{\kappa-2}(t)\right)^{\frac32},
			\end{aligned}
		\end{equation}
		where
		\begin{equation}\nonumber
			\begin{aligned}
				\sum_{|\beta|+|\gamma|\leq\kappa-4} \left\| \langle\tau\rangle^{-\frac12} \Gamma^\beta\Theta\,\Gamma^\gamma\Theta \right\|_{L_\tau^1 L_x^{6/5}} &\lesssim \sum_{\substack{ |\beta|\leq\kappa-4\\
						|\gamma|\leq\kappa-7 }} \left\| \langle\tau\rangle^{\frac12}\Gamma^\beta\Theta \right\|_{L_\tau^\infty L_x^2} \left\| \langle\tau\rangle^{\frac12}\Gamma^\gamma\Theta \right\|_{L_\tau^\infty L_x^3}\lesssim \mathcal E^{low}_{\kappa-2}(t).
			\end{aligned}
		\end{equation}
		and
		\begin{equation}\nonumber
			\begin{aligned}
				&\sum_{|\beta|+|\gamma|\leq\kappa-4}\left( \left\| \langle\tau\rangle^{-\frac12} \nabla\widetilde\Gamma^\beta u\, \Gamma^\gamma\Theta \right\|_{L_\tau^1 L_x^{6/5}}+\left\| \langle\tau\rangle^{-\frac12} \nabla\widetilde\Gamma^\beta u\nabla\widetilde{\Gamma}^\gamma u^{cf} \right\|_{L_\tau^1 L_x^{6/5}}\right) \\
				&\quad\lesssim \sum_{|\beta|+|\gamma|\leq\kappa-4} \left\|\nabla\widetilde\Gamma^\beta u\right\|_{L_\tau^\infty L_x^2} \left\|\langle\tau\rangle^{-1}\right\|_{L_\tau^{4/3}}\\
				&\qquad\times\left( \left\| \langle\tau\rangle^{\frac12}\Gamma^\gamma\Theta \right\|_{L_\tau^2 L_x^6}^{\frac12} \left\| \langle\tau\rangle^{\frac12}\Gamma^\gamma\Theta \right\|_{L_\tau^\infty L_x^2}^{\frac12}+ \left\| \langle\tau\rangle^{\frac12}\nabla\widetilde{\Gamma}^\gamma u^{cf} \right\|_{L_\tau^2 L_x^6}^{\frac12} \left\| \langle\tau\rangle^{\frac12}\nabla\widetilde{\Gamma}^\gamma u^{cf} \right\|_{L_\tau^\infty L_x^2}^{\frac12}\right) \\
				&\quad\lesssim \mathcal E^{low}_{\kappa-2}(t).
			\end{aligned}
		\end{equation}
		For the elastic nonlinear terms involving two divergence-free factors,
		substituting $w=\partial_\tau(-\Delta)^{-1}\widetilde\Gamma^\alpha u^{cf}$ into \eqref{equ_nu_df_df}, we obtain
		\begin{equation}\label{equ_F1_2}
			\begin{aligned}
				\sum_{|\alpha|\leq\kappa-4}& \int_0^t\int_{\mathbb R^3} \partial_\tau\nabla^{-1}\widetilde\Gamma^\alpha u^{cf} \cdot\nabla^{-1}F_{1,\alpha} dx d\tau \lesssim{} \int_0^t \langle\tau\rangle^{-\frac32+\frac{\delta}{2}} \widetilde{\mathcal J}_{\kappa-2}(\tau )\widetilde{\mathcal I}_\kappa^{\frac12}(\tau)\,d\tau+\left(\mathcal E^{low}_{\kappa-2}(t)\right)^{\frac32}.
			\end{aligned}
		\end{equation}
		
		We next estimate the nonlinear terms involving $F_{2,\alpha}$.
		Using the expansion of $F_{2,\alpha}$ and \eqref{eq:F1-alpha-theta-estimate}, we obtain
		\begin{equation}\label{eq:J-fourth-F2}
			\begin{aligned}
				&\sum_{|\alpha|\leq\kappa-4} \int_0^t\!\int_{\mathbb R^3} \Bigl( A_1\nabla\nabla^{-1}\Gamma^\alpha\Theta \cdot\nabla\nabla^{-1}F_{2,\alpha} +A_2\partial_\tau\widetilde\Gamma^\alpha u^{cf} \cdot\nabla^{-1}F_{2,\alpha} +A_1\nabla^{-1}\Gamma^\alpha\Theta \cdot\nabla^{-1} F_{2,\alpha} \Bigr)\,dx\,d\tau \\
				&\quad\lesssim \sum_{\substack{|\alpha|\leq\kappa-4\\ |\beta|+|\gamma|\leq\kappa-4}} \Bigl( \|\nabla\Gamma^\alpha\Theta\|_{L_\tau^\infty L_x^2} +\|\partial_\tau\widetilde\Gamma^\alpha u^{cf}\|_{L_\tau^\infty L_x^2} +\|\nabla^{-1}\Gamma^\alpha\Theta\|_{L_\tau^\infty L_x^2} \Bigr) \left\| \Gamma^\beta\Theta\, \nabla\partial_\tau\widetilde\Gamma^\gamma u^{cf} \right\|_{L_\tau^1 L_x^{6/5}}  \\
				&\qquad+ \sum_{|\alpha|\leq\kappa-4} \Bigl( \|\langle \tau \rangle^\frac12 \nabla^2\Gamma^\alpha\Theta\|_{L_\tau^2 L_x^2} +\|\langle \tau \rangle^\frac12 \partial_\tau\nabla\widetilde\Gamma^\alpha u^{cf}\|_{L_\tau^2 L_x^2} +\|\Gamma^\alpha\Theta\|_{L_\tau^2 L_x^2} \Bigr)\\
				&\qquad\times\Bigg(\sum_{|\beta|+|\gamma|\leq\kappa-4} \left\|\Gamma^\beta\Theta\,\Gamma^\gamma\Theta \right\|_{L_\tau^2 L_x^2}+ \sum_{|\beta|\leq\kappa-5}\|\Gamma^\beta\Theta\|_{L_\tau^2 L_x^2}\Bigg) \\
				&\quad\lesssim \left(\mathcal E^{low}_{\kappa-2}(t)\right)^{\frac32} + \widetilde{\mathcal J}_{\kappa-2}^{\frac12}(t) \widetilde{\mathcal J}_{\kappa-3}^{\frac12}(t).
			\end{aligned}
		\end{equation}
		where
		\begin{equation}\nonumber
			\begin{aligned}
				\sum_{|\beta|+|\gamma|\leq\kappa-4}& \left\| \Gamma^\beta\Theta\, \nabla\partial_\tau\widetilde\Gamma^\gamma u^{cf} \right\|_{L_\tau^1 L_x^{6/5}} 
				\lesssim{} \sum_{|\beta|,|\gamma|\leq\kappa-4} \left\| \langle\tau\rangle^{\frac12} \nabla\partial_\tau\widetilde\Gamma^\gamma u^{cf} \right\|_{L_\tau^2 L_x^2} \left\| \langle\tau\rangle^{\frac12}\Gamma^\beta\Theta \right\|_{L_\tau^\infty L_x^3}\lesssim \mathcal E^{low}_{\kappa-2}(t).
			\end{aligned}
		\end{equation}
		
		Finally, combining \eqref{equ_low_F1_1}, \eqref{equ_F1_2}, and \eqref{eq:J-fourth-F2} with the initial energy inequality yields the desired estimate and completes the proof.
	\end{proof}
	
	\section{The proof of Theorem \ref{thm:global-existence}}\label{4}
	
	In this section, we complete the proof of Theorem
	\ref{thm:global-existence} via a bootstrap argument.
	Let $\kappa\geq 9$, and suppose that
	$(u,\Theta)\in \dot H_{\widetilde\Gamma}^{\kappa}(T)\times\dot H_{\Gamma}^{\kappa}(T)$
	is a solution of \eqref{equ_1}.
	The proof consists of two key steps.
	We first derive a pair of integral inequalities for the lower-order energy
	$\mathcal E^{low}_{\kappa-2}(t)$ and the higher-order energy
	$\mathcal E_{\kappa}(t)$.
	We then use these inequalities to close the bootstrap argument.
	
	\subsection{Higher order energy estimate}
	
	We next estimate the higher-order energy $\mathcal I_\kappa(t)$.
	The energy method yields
	\begin{equation}\label{eq:I-high-energy-identity}
		\begin{aligned}
			\mathcal I_\kappa(t) ={}&\langle t\rangle ^{-\delta}\mathcal I_\kappa(0) +\langle t\rangle^{-\delta} \sum_{\substack{|\alpha|\leq\kappa-1\\|a|\leq1}} \int_0^t\int_{\mathbb R^3} \Bigl\langle \partial_\tau\nabla^a\widetilde\Gamma^\alpha u,\, L\nabla^a\widetilde\Gamma^\alpha u +\nabla\nabla^a\Gamma^\alpha\Theta \Bigr\rangle\,dx\,d\tau \\
			&\quad +\langle t\rangle^{-\delta} \sum_{\substack{|\alpha|\leq\kappa-1\\|a|\leq1}} \int_0^t\int_{\mathbb R^3} \Bigl\langle \nabla^a\Gamma^\alpha\Theta,\, T\nabla^a\Gamma^\alpha\Theta +\partial_\tau\nabla\cdot \nabla^a\widetilde\Gamma^\alpha u \Bigr\rangle\,dx\,d\tau .
		\end{aligned}
	\end{equation}
	Substituting the commuted equations (see Lemma \ref{lem_com}) into
	\eqref{eq:I-high-energy-identity}, isolating the terms involving
	the highest-order derivatives, and integrating by parts using
	the symmetry conditions, we obtain
	
	\begin{equation}\label{eq:I-high-expanded}
		\begin{aligned}
			\mathcal I_\kappa&(t)-\mathcal I_\kappa(0)\\
			\lesssim{}& -\langle t\rangle^{-\delta} \sum_{\substack{|\alpha|\leq\kappa-1\\|a|\leq1}} \int_{\mathbb R^3} \partial_l\nabla^a\widetilde\Gamma^\alpha u_i \partial_m\nabla^a\widetilde\Gamma^\alpha u_j \biggl[ g^{ijk}_{lmn} \partial_nu_k +\frac12g^{ij}_{lm} \Theta \biggr]\,dx \\
			&+\langle t\rangle^{-\delta} \sum_{\substack{|\alpha|\leq\kappa-1\\|a|\leq1}} \int_0^t\int_{\mathbb R^3} \partial_l\nabla^a\widetilde\Gamma^\alpha u_i\partial_m\nabla^a\widetilde\Gamma^\alpha u_j \biggl[ g^{ijk}_{lmn} \partial_\tau\partial_nu_k +\frac12g^{ij}_{lm} \partial_\tau\Theta \biggr]\,dx\,d\tau \\
			&+\langle t\rangle^{-\delta} \sum_{\substack{|\alpha|\leq\kappa-1\\|a|=1}} \int_0^t\int_{\mathbb R^3} \partial_\tau\nabla^a\widetilde\Gamma^\alpha u_i\biggl[ 2g^{ijk}_{lmn} \partial_l\bigl( \partial_m\widetilde\Gamma^\alpha u_j \partial_n\nabla^au_k \bigr) +g^{ij}_{lm} \partial_l\bigl( \partial_m\widetilde\Gamma^\alpha u_j \nabla^a\Theta \bigr) \biggr]\,dx\,d\tau \\
			&+\langle t\rangle^{-\delta} \sum_{\substack{|\alpha|\leq\kappa-1\\|a|\leq1}} \int_0^t\int_{\mathbb R^3} g^{ij}_{lm} \partial_\tau\nabla^a\widetilde\Gamma^\alpha u_i\, \nabla^a\partial_l \bigl(\partial_m u_j\,\Gamma^\alpha\Theta\bigr) \,dx\,d\tau \\
			&+\langle t\rangle^{-\delta} \sum_{\substack{|\alpha|\leq\kappa-1,\ |a|\leq1\\ |\beta|+|\gamma|\leq|\alpha|\\ \beta\neq\alpha,\ \gamma\neq\alpha}} 
			\int_0^t\int_{\mathbb R^3} \partial_\tau\nabla^a\widetilde\Gamma^\alpha u_i\biggl[ g^{ijk}_{lmn}\nabla^a\partial_l \bigl( \partial_m\widetilde\Gamma^\beta u_j\, \partial_n\widetilde\Gamma^\gamma u_k \bigr) +g^{ij}_{lm}\nabla^a\partial_l \bigl( \partial_m\widetilde\Gamma^\beta u_j\, \Gamma^\gamma\Theta \bigr) \biggr]\,dx\,d\tau \\
			&+\langle t\rangle^{-\delta} \sum_{\substack{|\alpha|\leq\kappa-1,\ |a|\leq1\\|\beta|+|\gamma|\leq|\alpha|}}
			\int_0^t\int_{\mathbb R^3}\left| \partial_\tau\nabla^a\widetilde\Gamma^\alpha u\cdot\nabla^{a} \bigl(\nabla\Gamma^\beta\Theta\,\Gamma^\gamma\Theta\bigr)\right| \,dx\,d\tau \\
			&+\langle t\rangle^{-\delta} \sum_{\substack{|\alpha|\leq\kappa-1, |a|\leq1\\|\beta|+|\gamma|\leq|\alpha|}} 
			C^1_{\beta,\gamma} \int_0^t\int_{\mathbb R^3} \nabla^a\Gamma^\alpha\Theta\, \nabla^a\Bigl( \Gamma^\beta\Theta\, \nabla\cdot\partial_\tau\widetilde\Gamma^\gamma u^{cf} \Bigr)\,dx\,d\tau \\
			&-\langle t\rangle^{-\delta} \sum_{\substack{|\alpha|\leq\kappa-1,\ |a|\leq1\\|\beta|+|\gamma|\leq|\alpha|}} 
			C^2_{\beta,\gamma} \int_0^t\int_{\mathbb R^3} \nabla^{1+a}\Gamma^\alpha\Theta\cdot \nabla^a\bigl(\nabla\Gamma^\beta\Theta\, \Gamma^\gamma\Theta\bigr) \,dx\,d\tau \\
			&+\langle t\rangle ^{-\delta} \sum_{\substack{|\alpha|\leq\kappa-1,\ |a|\leq1\\|\beta|\leq\kappa-2}}
			\left\|\nabla\nabla^a\Gamma^\alpha\Theta\right\| _{L_\tau^2 L_x^2}\left\| \nabla\nabla^a\Gamma^\beta\Theta\right\| _{L_\tau^2 L_x^2} .
		\end{aligned}
	\end{equation}
	
	We next estimate the terms in \eqref{eq:I-high-expanded} in order,
	\begin{equation}\label{eq:I-high-boundary}
		\begin{aligned}
			-\langle t\rangle^{-\delta}& \sum_{\substack{|\alpha|\leq\kappa-1\\|a|\leq1}} \int_{\mathbb R^3}  \partial_l\nabla^a\widetilde\Gamma^\alpha u_i \partial_m\nabla^a\widetilde\Gamma^\alpha u_j \biggl[ g^{ijk}_{lmn} \partial_nu_k +\frac12g^{ij}_{lm} \Theta \biggr]\,dx 
			\lesssim{}\mathcal E_\kappa(t)\left(\mathcal E^{low}_{\kappa-2}(t)\right)^{\frac12}.
		\end{aligned}
	\end{equation}
	Integrating by parts and applying Lemma \ref{lem_sob_1}, we obtain
	\begin{equation}\nonumber
		\begin{aligned}
			\langle t\rangle^{-\delta}&\sum_{\substack{|\alpha|\leq\kappa-1\\|a|\leq1}} \int_0^t\!\int_{\mathbb R^3} \partial_l\nabla^a\widetilde\Gamma^\alpha u_i g^{ijk}_{lmn} \partial_m\nabla^a\widetilde\Gamma^\alpha u_j \partial_\tau\partial_nu_kdxd\tau\\
			\lesssim{}&\sum_{\substack{|\alpha|\le \kappa-1\\|a|\le 1}}
			\left\| \langle \tau \rangle^{-\frac\delta2}\partial\nabla^a\widetilde{\Gamma}^\alpha u \right\| _{L_\tau^\infty L_x^2}^2\left\| \langle \tau \rangle ^{\frac12}\partial \nabla u^{cf}\right\| _{L_\tau^2 L_x^\infty} \\
			&+\sum_{\substack{|\alpha|\le \kappa-1\\|a|\le 1}}\langle t\rangle ^{-\delta}\int_{0}^t\langle \tau\rangle ^{-1+\delta}\left\| \langle \tau \rangle^{-\frac\delta2}\partial\nabla^a\widetilde{\Gamma}^\alpha u \right\|^2_{L_x^2}\left\| \langle \tau\rangle \partial \nabla u^{df}\right\| _{L_x^\infty}\\
			\lesssim{}& \langle t\rangle ^{-\delta}\int_{0}^t\langle \tau\rangle ^{-1+\delta}\widetilde{\mathcal I}_{\kappa}(\tau)\widetilde{\mathcal J}^{\frac12}_{\kappa-2}(\tau)d\tau + \mathcal E_\kappa(t)\left(\mathcal E^{low}_{\kappa-2}(t)\right)^{\frac12}.
		\end{aligned}
	\end{equation}
	The same argument also gives
	\begin{equation}
		\begin{aligned}
			\langle t\rangle^{-\delta}& \sum_{\substack{|\alpha|\leq\kappa-1\\|a|=1}} \int_0^t\int_{\mathbb R^3} \partial_\tau\nabla^a\widetilde\Gamma^\alpha u_i  \partial_l\bigl( g^{ijk}_{lmn} \partial_m\widetilde\Gamma^\alpha u_j \partial_n\nabla^au_k \bigr)\,dx\,d\tau\\
			&\lesssim \langle t\rangle ^{-\delta}\int_{0}^t\langle \tau\rangle ^{-1+\delta}\widetilde{\mathcal I}_{\kappa}(\tau)\widetilde{\mathcal J}^{\frac12}_{\kappa-2}(\tau)d\tau + \mathcal E_\kappa(t)\left(\mathcal E^{low}_{\kappa-2}(t)\right)^{\frac12}.
		\end{aligned}
	\end{equation}
	
	We next consider the terms whose nonlinear parts contain no highest-order factors, namely,
	\begin{equation}
		\begin{aligned}
			\langle t\rangle^{-\delta}&\sum_{\substack{|\beta|,|\gamma|<|\alpha|\leq\kappa-1,\\
					|\beta|+|\gamma|\leq|\alpha|, |a|\leq1}} \int_0^t\!\int_{\mathbb R^3} \partial_\tau\nabla^a\widetilde\Gamma^\alpha u_i g^{ijk}_{lmn} \nabla^a\partial_l\!\left( \partial_m\widetilde\Gamma^\beta u_j\, \partial_n\widetilde\Gamma^\gamma u_k \right)dxd\tau\\
			\lesssim{}&\langle t\rangle^{-\delta}\sum_{\substack{|\beta|,|\gamma|<|\alpha|\leq\kappa-1,\\|\beta|+|\gamma|\leq|\alpha|,|a|=1}} 
			\int_0^t\!\int_{\mathbb R^3}2 g^{ijk}_{lmn} \partial_\tau\nabla^a\widetilde\Gamma^\alpha u_i\Big(\partial_l\partial_m\widetilde\Gamma^\beta u_j \partial_n\nabla^a\widetilde\Gamma^\gamma u_k
			+\partial_l\partial_m\nabla^{a}\widetilde\Gamma^\beta u_j \partial_n\widetilde\Gamma^\gamma u_k\big)dxd\tau\\
			&+\sum_{\substack{|\beta|,|\gamma|<|\alpha|\leq\kappa-1,\\ |\beta|+|\gamma|\leq|\alpha|}}
			\left\| \langle \tau\rangle ^{-\frac\delta2}\partial \widetilde{\Gamma}^\alpha u\right\| _{L_\tau^\infty L_x^2}\left\|\langle \tau \rangle^\frac12 \nabla^2\widetilde{\Gamma}^\beta u\nabla \widetilde{\Gamma}^\gamma u\right\| _{L_\tau^2 L_x^2} .
		\end{aligned}
	\end{equation}
	Since these terms have different derivative structures, we estimate them separately. 
	For the first term, we treat the case $|\gamma|=\kappa-5$ separately. By \eqref{equ_1.3},
	\begin{equation}\nonumber
		\begin{aligned}
			\langle t\rangle^{-\delta}&\sum_{\substack{|\beta|,|\gamma|<|\alpha|\leq\kappa-1,\\
					|\beta|+|\gamma|\leq|\alpha|,|a|=1}} \int_0^t\!\int_{\mathbb R^3}2 g^{ijk}_{lmn} \partial_\tau\nabla^a\widetilde\Gamma^\alpha u_i\partial_l\partial_m\widetilde\Gamma^\beta u_j \partial_n\nabla^a\widetilde\Gamma^\gamma u_k\,dx\,d\tau\\
			\lesssim{}& \sum_{\substack{|\alpha|\leq\kappa-1,\ |\beta|\leq\kappa-2\\
					|\gamma|\leq\kappa-6}}\langle t\rangle^{-\frac\delta2} \bigl\|\langle\tau\rangle^{-\frac{\delta}{2}} \partial\nabla\widetilde\Gamma^\alpha u\bigr\|_{L_\tau^\infty L_x^2} \bigl\|\langle\tau\rangle^{-\frac{\delta}{2}}\nabla^2\widetilde\Gamma^\beta u\bigr\|_{L_\tau^\infty L_x^2} \bigl\|\langle\tau\rangle^{\frac12+\frac\delta2} \nabla^2\widetilde\Gamma^\gamma u \bigr\|_{L_\tau^2 L_x^\infty}\\
			&+ \sum_{\substack{|\alpha|\leq\kappa-1,\ |\beta|=4\\
					|\gamma|=\kappa-5}} \bigl\|\langle\tau\rangle^{-\frac{\delta}{2}} \partial\nabla\widetilde\Gamma^\alpha u\bigr\|_{L_\tau^\infty L_x^2} \bigl\|\langle\tau\rangle^{-\frac{\delta}{2}}\langle r\rangle \nabla^2\widetilde\Gamma^\beta u\bigr\|_{L_\tau^\infty L_x^\infty}\\
			&\quad\times\biggl( \bigl\|\langle\tau\rangle^{\frac12} \nabla^2\widetilde\Gamma^\gamma u^{cf} \bigr\|_{L_\tau^2 L_x^2} + \bigl\|\langle c_2\tau-r\rangle \nabla^2\widetilde\Gamma^\gamma u^{df} \bigr\|_{L_\tau^\infty L_x^2} \biggr) \\
			\lesssim{}& \mathcal E_\kappa(t) \left(\mathcal E^{low}_{\kappa-2}(t)\right)^{\frac12}.
		\end{aligned}
	\end{equation}
	Similarly, we estimate the case $|\beta|+2=\kappa-4$ separately. By \eqref{equ_1.3}, we obtain
	\begin{equation}\nonumber
		\begin{aligned}
			\langle t\rangle ^{-\delta}&\sum_{\substack{|\beta|,|\gamma|<|\alpha|\leq\kappa-1,\\
					|\beta|+|\gamma|\leq|\alpha|,|a|=1}} \int_0^t\!\int_{\mathbb R^3}2 g^{ijk}_{lmn} \partial_\tau\nabla^a\widetilde\Gamma^\alpha u_i\partial_l\partial_m\nabla^{a}\widetilde\Gamma^\beta u_j \partial_n\widetilde\Gamma^\gamma u_k\,dx\,d\tau\\
			\lesssim{}& \sum_{\substack{|\alpha|\leq\kappa-1,\ |\gamma|\leq\kappa-2\\
					|\beta|+2\leq\kappa-5}}\langle t\rangle^{-\frac\delta2} \bigl\|\langle\tau\rangle^{-\frac{\delta}{2}} \partial\nabla\widetilde\Gamma^\alpha u\bigr\|_{L_\tau^\infty L_x^2} \bigl\|\langle\tau\rangle^{\frac12+\frac{\delta}{2}} \nabla^3\widetilde\Gamma^\beta u \bigr\|_{L_\tau^2 L_x^\infty} \bigl\|\langle \tau \rangle^{-\frac\delta2}\nabla\widetilde\Gamma^\gamma u\bigr\|_{L_\tau^\infty L_x^2} \\
			&+\sum_{\substack{|\alpha|\leq\kappa-1,\ |\beta|\leq\kappa-2\\
					|\gamma|\leq\kappa-5}}\langle t\rangle^{-\frac\delta2} \bigl\|\langle\tau\rangle^{-\frac{\delta}{2}} \partial\nabla\widetilde\Gamma^\alpha u\bigr\|_{L_\tau^\infty L_x^2}\bigl\|\langle \tau \rangle^{-\frac\delta2}\nabla^3\widetilde\Gamma^\beta u\bigr\|_{L_\tau^\infty L_x^2}\bigl\|\langle\tau\rangle^{\frac12+\frac{\delta}{2}}\nabla\widetilde\Gamma^\gamma u\bigr\|_{L_\tau^2 L_x^\infty} \\
			&+\sum_{\substack{|\alpha|\leq\kappa-1,\ |\gamma|\leq 5\\
					|\beta|+2=\kappa-4}} \bigl\|\langle\tau\rangle^{-\frac{\delta}{2}} \partial\nabla\widetilde\Gamma^\alpha u\bigr\|_{L_\tau^\infty L_x^2} \bigl\| \nabla^3\widetilde\Gamma^\beta u\bigr\|_{L_\tau^\infty L_x^2} \\
			&\quad\times \biggl(\langle t\rangle^{-\frac\delta2} \bigl\|\langle\tau\rangle^{\frac12} \nabla\widetilde\Gamma^\gamma u^{cf} \bigr\|_{L_\tau^2 L_x^\infty} + \bigl\|\langle \tau \rangle^{-\frac\delta2}\langle \tau\rangle \nabla\widetilde\Gamma^\gamma u^{df} \bigr\|_{L_\tau^\infty L_x^\infty} \biggr) \\
			\lesssim{}& \mathcal E_\kappa(t) \left(\mathcal E^{low}_{\kappa-2}(t)\right)^{\frac12}.
		\end{aligned}
	\end{equation}
	Combining these estimates with \eqref{equ_F11_1.1}, we obtain
	\begin{equation}\label{eq:i-fourth-term-final}
		\begin{aligned}
			& \sum_{\substack{|\beta|,|\gamma|<|\alpha|\leq\kappa-1,\\
					|\beta|+|\gamma|\leq|\alpha|, |a|\leq1}} \langle t\rangle^{-\delta}\int_0^t\!\int_{\mathbb R^3} \partial_\tau\nabla^a\widetilde\Gamma^\alpha u_i g^{ijk}_{lmn} \nabla^a\partial_l\!\left( \partial_m\widetilde\Gamma^\beta u_j\, \partial_n\widetilde\Gamma^\gamma u_k \right)dxd\tau\lesssim \mathcal E_\kappa(t)\left(\mathcal E^{low}_{\kappa-2}(t)\right)^{\frac12}.
		\end{aligned}
	\end{equation}
	When the highest-order derivatives fall on $\Theta$, \eqref{equ_1.3} gives
	\begin{equation}\label{equ_high-order-Theta-remainder}
		\begin{aligned}
			\langle t\rangle^{-\delta}& \sum_{\substack{ |\alpha|\leq\kappa-1\\|a|\leq1 }} 
			\int_0^t\int_{\mathbb R^3} g^{ij}_{lm}\, \partial_\tau\nabla^a \widetilde\Gamma^\alpha u_i\, \nabla^a\partial_l \bigl(\partial_m u_j\,\Gamma^\alpha\Theta\bigr) \,dx\,d\tau \\
			\lesssim{}&\langle t\rangle^{-\frac\delta2} \sum_{|\alpha|\leq\kappa-1} \bigl\| \langle \tau \rangle^{-\frac\delta2}\partial\nabla^{\le1}\widetilde\Gamma^\alpha u \bigr\|_{L_\tau^\infty L_x^2} \Biggl[ \langle t\rangle^{-\frac\delta2}\bigl\|\nabla\nabla^{\le1}\Gamma^\alpha\Theta\bigr\|_{L_\tau^2 L_x^2} \bigl\| \langle\tau\rangle^{\frac12+\frac\delta2} \nabla\nabla^{\le2} u \bigr\|_{L_\tau^2 L_x^\infty}\\
			& + \bigl\|\Gamma^\alpha\Theta\bigr\|_{L_\tau^2 L_x^6} \bigl\| \langle\tau\rangle^{\frac12} \nabla^2\nabla^{\le 1}u^{cf} \bigr\|_{L_\tau^2 L_x^3} 
			+\bigl\| \langle\tau\rangle^{-\frac\delta2} \Gamma^\alpha\Theta \bigr\|_{L_\tau^\infty L_x^2}\bigl\| \langle \tau\rangle \nabla^2 \nabla^{\le1}u^{df} \bigr\|_{L_\tau^\infty L_x^\infty}\left\| \langle \tau \rangle^{-1+\frac\delta2}\right\| _{L_\tau^1} \Biggr] \\
			\lesssim{}& \mathcal E_\kappa(t)\left(\mathcal E^{low}_{\kappa-2}(t)\right)^{\frac12}.
		\end{aligned}
	\end{equation}
	Since $\Theta$ and $u^{cf}$ have similar structures, the same argument, together with \eqref{equ_high-order-Theta-remainder}, yields
	\begin{equation}\label{equ_A11}
		\begin{aligned}
			\langle t\rangle^{-\delta}\sum_{\substack{|\alpha|\le \kappa-1,|a|\le1\\ \beta+\gamma=\alpha}}&\int_0^t\int_{\mathbb R^3} \partial_\tau\nabla^a\widetilde\Gamma^\alpha u_i\biggl[ g^{ijk}_{lmn}\nabla^a\partial_l \bigl( \partial_m\widetilde\Gamma^\beta u_j\, \partial_n\widetilde\Gamma^\gamma u_k \bigr) +g^{ij}_{lm}\nabla^a\partial_l \bigl( \partial_m\widetilde\Gamma^\beta u_j\, \Gamma^\gamma\Theta \bigr) \biggr]\,dx\,d\tau \\    
			\lesssim{}& \langle t\rangle ^{-\delta}\int_{0}^t\langle \tau\rangle ^{-1+\delta}\widetilde{\mathcal I}_{\kappa}(\tau)\widetilde{\mathcal J}^{\frac12}_{\kappa-2}(\tau)d\tau + \mathcal E_\kappa(t)\left(\mathcal E^{low}_{\kappa-2}(t)\right)^{\frac12}.
		\end{aligned}
	\end{equation}
	For the quadratic terms involving only the temperature in the elastic equation, \eqref{equ_F11_1.3} yields
	\begin{equation}
		\begin{aligned}
			&\langle t\rangle^{-\delta} \sum_{\substack{|\alpha|\leq\kappa-1\\|\beta|+|\gamma|\leq|\alpha|}} 
			\int_0^t\int_{\mathbb R^3}\left| \partial_\tau\nabla^{\le 1}\widetilde\Gamma^\alpha u \nabla^{\le1} \bigl(\nabla\Gamma^\beta\Theta\,\Gamma^\gamma\Theta\bigr)\right|dxd\tau
			\lesssim \mathcal E_\kappa(t)\left(\mathcal E^{low}_{\kappa-2}(t)\right)^{\frac12}.
		\end{aligned}
	\end{equation}
	
	We next estimate the $M_1$ term in the heat equation,
	\begin{equation}
		\begin{aligned}
			&\langle t\rangle^{-\delta} \sum_{\substack{|\alpha|\leq\kappa-1,\ |a|\leq1\\
					|\beta|+|\mu|\leq|\alpha|}} \int_0^t\int_{\mathbb R^3} \nabla^a\Gamma^\alpha\Theta\, \nabla^a\Bigl( \Gamma^\beta\Theta\, \partial_\tau\nabla\cdot\widetilde\Gamma^\mu u^{cf} \Bigr)\,dx\,d\tau \\
			&\quad\lesssim \sum_{\substack{|\alpha|,|\beta|\leq\kappa-1\\
					|\mu|+1\leq\kappa-4}} \langle t\rangle^{-\delta} \bigl\|\Delta\Gamma^\alpha\Theta\bigr\|_{L_\tau^2 L_x^2} \bigl\| \Gamma^\beta\Theta\bigr\|_{L_\tau^2 L_x^\infty} \bigl\| \partial_\tau\nabla \widetilde\Gamma^\mu u^{cf}\bigr\|_{L_\tau^\infty L_x^2} \\
			&\qquad+ \sum_{\substack{|\alpha|,|\mu|\leq\kappa-1\\
					|\beta|\leq\kappa-4}} \langle t\rangle^{-\frac{\delta}{2}} \bigl\|\Delta\Gamma^\alpha\Theta\bigr\|_{L_\tau^2 L_x^2} \bigl\|\langle\tau\rangle^{-\frac\delta2} \partial_\tau\nabla \widetilde\Gamma^\mu u^{cf}\bigr\|_{L_\tau^\infty L_x^2} \bigl\| \Gamma^\beta\Theta\bigr\|_{L_\tau^2 L_x^\infty} \\
			&\qquad+ \sum_{\substack{|\alpha|,|\mu|\leq\kappa-1\\
					|\beta|\leq\kappa-5}} \langle t\rangle^{-\frac{\delta}{2}}\bigl\|\langle\tau\rangle^{-\frac\delta2} \partial_\tau\nabla\widetilde\Gamma^\mu u^{cf} \bigr\|_{L_\tau^\infty L_x^2} \bigl\|\Gamma^\alpha\Theta\bigr\|_{L_\tau^2 L_x^6} \bigl\| \Gamma^\beta\Theta\bigr\|_{L_\tau^2 L_x^3} \\
			&\qquad+ \sum_{\substack{|\beta|,|\alpha|\leq\kappa-1\\
					|\mu|\leq\kappa-5}}\langle t\rangle^{-\frac\delta2}\bigl\|\langle\tau\rangle^{\frac12} \partial_\tau\nabla\widetilde\Gamma^\mu u^{cf} \bigr\|_{L_\tau^2 L_x^2} \bigl\|\Gamma^\alpha\Theta\bigr\|_{L_\tau^2 L_x^6} \bigl\|\langle \tau \rangle^{-\frac\delta2} \Gamma^\beta\Theta\bigr\|_{L_\tau^\infty L_x^3} \\
			&\quad\lesssim \mathcal E_\kappa(t)\left(\mathcal E^{low}_{\kappa-2}(t)\right)^{\frac12}.
		\end{aligned}
	\end{equation}
	
	For the last nonlinear term $M_2$ in the heat equation, by \eqref{equ_F11_1.3}, we have
	\begin{equation}\label{eq:i-low-term-10}
		\begin{aligned}
			\langle t\rangle^{-\delta}&\sum_{\substack{ |\alpha|\leq\kappa-1,|a|\le 1\\|\beta|+|\gamma|\leq\kappa-1 }} 
			\int_0^t\int_{\mathbb R^3} \nabla^{1+a}\Gamma^\alpha\Theta\cdot \nabla^a\bigl( \nabla\Gamma^\beta\Theta\, \Gamma^\gamma\Theta \bigr)\,dx\,d\tau \\
			\lesssim{}&\sum_{\substack{ |\alpha|\leq\kappa-1,|a|\le 1\\|\beta|+|\gamma|\leq\kappa-1 }}
			\langle t\rangle^{-\delta} \left\|\nabla^{1+a}\Gamma^\alpha\Theta\right\|_{L_\tau^2 L_x^2}\left\| \nabla^a\bigl( \nabla\Gamma^\beta\Theta\,\Gamma^\gamma\Theta \bigr)\right\|_{L_\tau^2 L_x^2} \\
			\lesssim{}& \mathcal E_\kappa(t)\left(\mathcal E^{low}_{\kappa-2}(t)\right)^{\frac12}.
		\end{aligned}
	\end{equation}
	
	Combining \eqref{eq:I-high-boundary},
	\eqref{eq:i-fourth-term-final} - \eqref{eq:i-low-term-10} with \eqref{eq:I-high-expanded}
	and applying the quadratic inequality, we obtain
	\begin{equation}
		\begin{aligned}
			\mathcal I_{\kappa}(t) \lesssim{}& \langle t\rangle ^{-\delta}
			\int_{0}^t\langle \tau\rangle ^{-1+\delta}\widetilde{\mathcal I}_{\kappa}(\tau)\widetilde{\mathcal J}^{\frac12}_{\kappa-2}(\tau)d\tau 
			+\mathcal E_\kappa(t)\left(\mathcal E^{low}_{\kappa-2}(t)\right)^{\frac12}+ \mathcal I_{\kappa}(0) +\mathcal I^{\frac12}_{\kappa}(t)\mathcal I^{\frac12}_{\kappa-1}(t).
		\end{aligned}
	\end{equation}
	
	\subsection{Lower order energy estimate}
	
	The energy method yields
	\begin{equation}
		\begin{aligned}
			J_{\kappa-2}(t)-J_{\kappa-2}(0) ={}&\sum_{\substack{|\alpha|\leq\kappa-3\\|a|\leq1}} \int_0^t\!\int_{\mathbb R^3} \Bigl\langle \partial_\tau\nabla^a\widetilde\Gamma^\alpha u,\, L\nabla^a\widetilde\Gamma^\alpha u +\nabla\nabla^a\Gamma^\alpha\Theta \Bigr\rangle\,dx\,d\tau\\
			&+\sum_{\substack{|\alpha|\leq\kappa-3\\|a|\leq1}} \int_0^t\!\int_{\mathbb R^3} \nabla^a\Gamma^\alpha\Theta \Bigl( T\nabla^a\Gamma^\alpha\Theta +\partial_\tau\nabla\cdot \nabla^a\widetilde\Gamma^\alpha u \Bigr)\,dx\,d\tau.
		\end{aligned}
	\end{equation}
	Using \eqref{equ_null_condition_2}, we obtain
	\begin{equation}\label{eq:J-low-expanded}
		\begin{aligned}
			J_{\kappa-2}&(t)-J_{\kappa-2}(0)\\
			\lesssim{}& -\sum_{\substack{|\alpha|\leq\kappa-3\\|a|\leq1}} \int_{\mathbb R^3} \partial_l\nabla^a\widetilde\Gamma^\alpha u_i \partial_m\nabla^a\widetilde\Gamma^\alpha u_j \left.\left[ g^{ijk}_{lmn} \partial_nu_k+\frac12 g^{ij}_{lm} \Theta \right]dx\right| _{\tau=0}^{\tau=t}\\
			&+\sum_{\substack{|\alpha|\leq\kappa-3\\|a|\leq1}} \int_0^t\!\int_{\mathbb R^3} \partial_l\nabla^a\widetilde\Gamma^\alpha u_i \partial_m\nabla^a\widetilde\Gamma^\alpha u_j\left( g^{ijk}_{lmn} \partial_\tau\partial_nu_k+\frac12 g^{ij}_{lm} \partial_\tau\Theta\right) dxd\tau\\
			&+\sum_{\substack{|\alpha|\leq\kappa-3\\|a|=1}} \int_0^t\!\int_{\mathbb R^3} \partial_\tau\nabla^a\widetilde\Gamma^\alpha u_i\left[ 2 g^{ijk}_{lmn} \partial_l\!\left( \partial_m\widetilde\Gamma^\alpha u_j\, \partial_n\nabla^a u_k \right)+g^{ij}_{lm} \partial_l\!\left( \partial_m\widetilde\Gamma^\alpha u_j\, \nabla^a\Theta \right)\right] dxd\tau\\
			&+\sum_{\substack{|\alpha|\leq\kappa-3\\|a|\leq1}} \int_0^t\!\int_{\mathbb R^3}\left| g^{ij}_{lm} \partial_\tau\partial_l\widetilde\Gamma^\alpha u_i\, \nabla^{2a}\left( \partial_m u_j\,\Gamma^\alpha\Theta \right)\right| dx\,d\tau\\
			&+\sum_{\substack{|\alpha|\leq\kappa-3,\ |a|\leq1\\
					|\beta|+|\gamma|\leq|\alpha|\\
					\beta\ne\alpha,\ \gamma\ne\alpha}} \int_0^t\!\int_{\mathbb R^3} \partial_\tau\nabla^a\widetilde\Gamma^\alpha u_i\left[ g^{ijk}_{lmn} \nabla^a\partial_l\!\left( \partial_m\widetilde\Gamma^\beta u_j\, \partial_n\widetilde\Gamma^\gamma u_k \right)+ g^{ij}_{lm}\, \nabla^a\partial_l\!\left( \partial_m\widetilde\Gamma^\beta u_j\, \Gamma^\gamma\Theta \right)\right] \,dx\,d\tau\\
			&+\sum_{\substack{|\alpha|\leq\kappa-3,\ |a|\leq1\\
					|\beta|+|\gamma|\leq|\alpha|}} \int_0^t\!\int_{\mathbb R^3}\left| \partial_\tau\nabla^a\widetilde\Gamma^\alpha u \cdot\nabla^{a}\left(\nabla\Gamma^\beta\Theta\, \Gamma^\gamma\Theta \right)\right| dx\,d\tau\\
			&+\sum_{\substack{|\alpha|\leq\kappa-3,\ |a|\leq1\\
					|\beta|+|\gamma|\leq|\alpha|}} \int_0^t\!\int_{\mathbb R^3}\nabla^a\Gamma^\alpha\Theta\left[ C^1_{\beta,\gamma} \nabla^a\!\left( \Gamma^\beta\Theta\, \partial_\tau\nabla\cdot \widetilde\Gamma^\gamma u^{\mathrm{cf}} \right)+C^2_{\beta,\gamma}\nabla^{1+a}\!\left( \nabla\Gamma^\beta\Theta\, \Gamma^\gamma\Theta \right)\right]dxd\tau\\
			&+\sum_{\substack{|\alpha|\leq\kappa-3,\ |a|\leq1\\
					|\beta|\leq\kappa-4}}\left\|\nabla\nabla^a\Gamma^\alpha\Theta\right\| _{L_\tau^2 L_x^2}\left\| \nabla\nabla^a\Gamma^\beta\Theta\right\| _{L_\tau^2 L_x^2}.
		\end{aligned}
	\end{equation}
	
	We next estimate the boundary terms and the remainders arising from integration by parts separately,
	\begin{equation}\label{eq:J-cn-initial-cubic}
		\begin{aligned}
			&\sup_{0\le \tau \le t}\sum_{\substack{|\alpha|\leq\kappa-3\\|a|\leq1}}
			\int_{\mathbb R^3} \left|\partial_l\nabla^a\widetilde\Gamma^\alpha u_i
			\partial_m\nabla^a\widetilde\Gamma^\alpha u_j
			\biggl(g^{ijk}_{lmn}\partial_nu_k+\frac12g^{ij}_{lm}\Theta\biggr)\right|dx\lesssim \left(\mathcal E^{low}_{\kappa-2}(t)\right)^{\frac32}.
		\end{aligned}
	\end{equation}
	We denote the remainders arising from integration by parts by $B_{1}(u,u,u)$ and estimate them case by case.
	\begin{equation}	
		\begin{aligned}
			B_{1}(u,u,u)\triangleq\sum_{\substack{|\alpha|\leq\kappa-3\\|a|\leq1}} \int_0^t\!\int_{\mathbb R^3} g^{ijk}_{lmn} \partial_l\nabla^a\widetilde\Gamma^\alpha u_i \partial_m\nabla^a\widetilde\Gamma^\alpha u_j \partial_\tau\partial_nu_kdxd\tau.
		\end{aligned}
	\end{equation}
	We first consider $B_1(u^{df}, u^{df},u^{df})$. Using \eqref{equ_nu_df_df}, the null structure, and Lemma \ref{lem_sob_1}, we obtain
	\begin{equation}\label{equ_4.90}
		\begin{aligned}
			\sum_{|\alpha|\le \kappa-3}&N(\nabla^{\le 1}\widetilde{\Gamma}^\alpha u^{df},\nabla^{\le 1}\widetilde{\Gamma}^\alpha u^{df},\partial_\tau u^{df})\\
			\lesssim{}& \langle\tau \rangle^{-\frac32}\sum_{|\alpha|\le \kappa-3}\left\| \nabla\nabla^{\le 1}\widetilde{\Gamma}^\alpha u^{df}\right\|_{L_x^2}
			\bigg(\left\| \langle c_2 \tau -r\rangle ^\frac12 \nabla\nabla^{\le 1}\widetilde{\Gamma}^\alpha u^{df}\right\|_{L_x^3}\left\|\langle c_2\tau -r\rangle \partial \nabla u^{df}\right\|_{L_x^6}\\
			&+\left\| \langle r \rangle ^\frac12 \nabla^{\le 1}\widetilde{\Gamma}^{\alpha+1} u^{df}\right\|_{L_x^{\infty,2}}\left\|\partial \nabla^2u^{df}\right\|_{L^{2,\infty}_x}+\left\|\nabla \nabla^{\le 1}\widetilde{\Gamma}^{\alpha+1} u^{df}\right\|_{L_x^2}\left\| \langle r\rangle ^\frac12 \partial \nabla u^{df}\right\|_{L_x^\infty}\bigg)\\
			\lesssim{}&\langle \tau \rangle^{-\frac32+\frac\delta2} \widetilde{ \mathcal{J}}_{\kappa-2}(t)\widetilde{\mathcal I}_\kappa^\frac12(t),
		\end{aligned}
	\end{equation}
	i.e.
	\begin{equation}\label{equ_B11_2}
		\begin{aligned}
			B_1&(u^{df},u^{df},u^{df})\lesssim{}\int_0^t \langle\tau\rangle^{-\frac32+\frac{\delta}{2}}\widetilde{\mathcal J}_{\kappa-2}(\tau )\widetilde{\mathcal I}_\kappa^{\frac12}(\tau)\,d\tau.
		\end{aligned}
	\end{equation}
	For $B_1(u,u^{cf},u^{df})$ and $B_1(u,u^{df},u^{cf})$, we have
	\begin{equation}\nonumber
		\begin{aligned}
			B_{1}&(u,u^{cf},u^{df}) + B_1(u,u^{df},u^{cf}) \\
			&\lesssim{}\sum_{\substack{|\alpha|\leq\kappa-3}} \int_0^t \langle\tau\rangle^{-\frac32+\frac{\delta}{2}} \|\partial\nabla^{\le 1}\widetilde\Gamma^\alpha u\|_{L_x^2} 
			\bigg(\left\| \langle\tau\rangle^{\frac12-\frac{\delta}{2}}\nabla\nabla^{\le1}\widetilde{\Gamma}^\alpha u^{cf} \right\| _{L_x^2} \left\|\langle \tau\rangle \partial\nabla u^{df}\right\|_{L_x^\infty}\\
			&\qquad + \left\|\langle\tau\rangle^{-\frac{\delta}{2}} \langle c_2\tau-r\rangle \partial\nabla^{\le1}\widetilde\Gamma^\alpha u^{df}\right\|_{L_x^6} \left\|\langle\tau\rangle^{\frac12}\partial\nabla u^{cf}\right\|_{L_x^3}
			+ \left\|\langle\tau\rangle^{-\frac{\delta}{2}} \partial\nabla^{\le1}\widetilde\Gamma^\alpha u^{df}\right\|_{L_x^2} \left\|\langle\tau\rangle^{\frac12} \langle r\rangle\partial\nabla u^{cf}\right\|_{L_x^\infty} \bigg)d\tau\\
			& \lesssim\int_0^t \langle\tau\rangle^{-\frac32+\frac{\delta}{2}}\widetilde{\mathcal J}_{\kappa-2}(\tau )\widetilde{\mathcal I}_\kappa^{\frac12}(\tau)\,d\tau.
		\end{aligned}
	\end{equation}
	Noting that we do not have an energy bound for $u^{cf}$ with a first-order time weight, we apply Lemma \ref{lem_Gamma_1} to obtain
	\begin{equation}\nonumber
		\begin{aligned}
			B_1(u^{cf}&,u^{cf},u^{cf})\lesssim{}\sum_{\substack{|\alpha|\leq\kappa-3\\ |\beta|\leq\kappa-4}} \sup_{0\le \tau\le t} \|\nabla\nabla^{\le1}\widetilde\Gamma^\alpha u\|_{L_x^2} \|\langle \tau \rangle^{\frac12}\nabla \nabla^{\le1}\widetilde\Gamma^\beta u^{cf}\|_{L_x^2} \|\langle \tau \rangle^{\frac12} \partial\nabla u^{cf}\|_{L_x^\infty}\\
			&+\sum_{\substack{|\alpha|\leq\kappa-2,|\gamma|\le 1\\ |\beta|\leq\kappa-4}} \int_0^t \langle\tau\rangle^{-\frac32+\frac{\delta}{2}} \|\langle\tau\rangle^{\frac12-\frac{\delta}{2}} \nabla\nabla^{\le1}\widetilde\Gamma^\alpha u^{cf}\|_{L_x^2} \|\langle\tau\rangle^{\frac12} \nabla\nabla^{\le 1}\widetilde\Gamma^\beta u^{cf}\|_{L_x^2} \|\langle\tau\rangle^{\frac12} \partial\nabla \widetilde{\Gamma}^\gamma u^{cf}\|_{L_x^\infty} \,d\tau\\
			\lesssim{}& \left(\mathcal E^{low}_{\kappa-2}(t)\right)^{\frac32} +\int_0^t \langle\tau\rangle^{-\frac32+\frac{\delta}{2}}\widetilde{\mathcal J}_{\kappa-2}(\tau )\widetilde{\mathcal I}_\kappa^{\frac12}(\tau)\,d\tau .
		\end{aligned}
	\end{equation}
	The remaining terms can be reduced to the cases above by symmetry. Thus,
	\begin{equation}\nonumber
		B_1(u,u,u) \lesssim \left(\mathcal E^{low}_{\kappa-2}(t)\right)^{\frac32} +\int_0^t \langle\tau\rangle^{-\frac32+\frac{\delta}{2}}\widetilde{\mathcal J}_{\kappa-2}(\tau )\widetilde{\mathcal I}_\kappa^{\frac12}(\tau)\,d\tau.
	\end{equation}
	
	We next consider the terms whose nonlinear parts contain
	no highest-order factors, i.e.
	\begin{equation}\label{eq:J-cn-fourth-integral}
		\begin{aligned}
			&\sum_{\substack{|\beta|,|\gamma|<|\alpha|\leq\kappa-3,\\ |\beta|+|\gamma|\leq|\alpha|, |a|\leq1}} 
			\left|  \int_0^t\!\int_{\mathbb R^3}\partial_\tau\nabla^a\widetilde\Gamma^\alpha u_i g^{ijk}_{lmn} \nabla^a\partial_l\!\left( \partial_m\widetilde\Gamma^\beta u_j\, \partial_n\widetilde\Gamma^\gamma u_k \right)dxd\tau\right|\\
			&\qquad\lesssim{} B_{2}(u,u,u)+B_{3}(u,u,u)+B_{4}(u,u,u)+B_5(u^{df},u^{df},u^{df}).
		\end{aligned}
	\end{equation}
	where $\nu=(\nu_1,\nu_2,\nu_3)$ specifies the Helmholtz component of each factor in the first three sums, and
	\begin{equation}\nonumber
		\begin{aligned}
			&B_{2}(u,u,u)={}\sum_{\substack{\nu\in\{cf,df\}^3,\ \nu\neq(df,df,df)\\|\beta|,|\gamma|<|\alpha|\leq\kappa-3,\\ |\beta|+|\gamma|\leq|\alpha|}} 
			\int_0^t\!\int_{\mathbb R^3}\left|g^{ijk}_{lmn} \partial\nabla\widetilde\Gamma^\alpha u^{\nu_1}\nabla^2\widetilde\Gamma^\beta u^{\nu_2}\nabla^2\widetilde\Gamma^\gamma u^{\nu_3}\right|dx\,d\tau\\
			&B_{3}(u,u,u)={}\sum_{\substack{\nu\in\{cf,df\}^3,\ \nu\neq(df,df,df)\\|\alpha|\leq\kappa-3,|\gamma|\le \kappa-4\\|\beta|+1\leq \kappa-6,|\beta|+|\gamma|\le \kappa-3}} 
			\int_0^t\!\int_{\mathbb R^3}\left|\partial\nabla\widetilde\Gamma^\alpha u^{\nu_1}   \nabla \nabla^{\le 2}\widetilde{\Gamma}^\beta u^{\nu_2}\nabla\widetilde{\Gamma}^\gamma u^{\nu_3}\right|dx\,d\tau\\
			&B_{4}(u,u,u)={}\sum_{\substack{\nu\in\{cf,df\}^3,\ \nu\neq(df,df,df)\\|\alpha|\leq\kappa-3,|\beta|\le \kappa-4\\ |\gamma|\leq \kappa-6}} 
			\int_0^t\!\int_{\mathbb R^3}\left|\partial\nabla\widetilde\Gamma^\alpha u^{\nu_1}\nabla^3\widetilde\Gamma^\beta u^{\nu_2}_j \nabla\widetilde\Gamma^\gamma u^{\nu_3}_k\right|dx\,d\tau\\
			&B_5(u^{df},u^{df},u^{df})={}\sum_{\substack{|\beta|,|\gamma|<|\alpha|\leq\kappa-3,\\ |\beta|+|\gamma|\leq|\alpha|, |a|\leq1}} 
			\int_0^t\!\int_{\mathbb R^3}\left| \partial_\tau\nabla^a\widetilde\Gamma^\alpha u^{df}_i g^{ijk}_{lmn} \nabla^a\partial_l\!\left( \partial_m\widetilde\Gamma^\beta u^{df}_j\, \partial_n\widetilde\Gamma^\gamma u^{df}_k \right)\right|dxd\tau.
		\end{aligned}
	\end{equation}
	
	We estimate the $df$ and $cf$ combinations separately.
	We first consider $B_{5}(u^{df},u^{df},u^{df})$.
	Using the null-form estimate \eqref{eq:2.21} in Lemma \ref{lem_null},
	the definition of $N_2(u,v)$ in \eqref{equ_N_al},
	and the weighted Sobolev estimates, we obtain
	
	\begin{equation}\nonumber
		\begin{aligned}
			&\sum_{\substack{|\beta|+|\gamma|\le \kappa-3\\ |\beta |,|\gamma|\le \kappa-4}}\left\| N_2(\nabla^{\le 1}\widetilde{\Gamma}^\beta u^{df},\widetilde{\Gamma}^{\gamma}u^{df})\right\|_{L_x^2}\\
			&\quad\lesssim{}\langle\tau \rangle^{-\frac32}\Bigg[ \sum_{\substack{|\beta|\le\kappa-6\\ |\gamma|\le \kappa-4}}\left\|\langle r\rangle^\frac12\nabla^{\le 2} \widetilde{\Gamma}^{\beta+1}u^{df}\right\|_{L_x^\infty}\left\| \nabla\nabla^{\le 1}\widetilde{\Gamma}^\gamma u^{df}\right\|_{L_x^2}\\
			&\qquad+\sum_{\substack{|\beta|\le\kappa-4\\ |\gamma|\le \kappa-6}}\bigg(\left\|\langle r\rangle^{\frac12}\widetilde{\Gamma}^{\beta+1}u^{df}\right\|_{L^{\infty,2}_x}\left\|\nabla^2\widetilde{\Gamma}^{\gamma}u^{df}\right\|_{L_x^{2,\infty}}
			+\left\|\nabla\nabla^{\le 1}\widetilde{\Gamma}^{\beta+1}u^{df}\right\|_{L_x^2}\left\| \langle r\rangle^\frac12 \nabla \widetilde{\Gamma}^\gamma u^{df}\right\|_{L_x^\infty} \bigg)\\
			&\qquad+\sum_{\substack{|\beta|\le \kappa-3\\ |\gamma|\le \kappa-6}} \Big(\left\| \langle c_2\tau -r \rangle\nabla^2\widetilde{\Gamma}^\beta u^{df} \right\|_{L_x^2}\left\|\langle c_2 \tau -r\rangle\nabla\widetilde{\Gamma}^\gamma u^{df}\right\|_{L_x^\infty}\\
			&\qquad\quad+\left\| \langle c_2\tau -r \rangle\nabla^2\widetilde{\Gamma}^\gamma u^{df} \right\|_{L_x^3}\left\|\langle c_2 \tau -r\rangle\nabla\widetilde{\Gamma}^\beta u^{df}\right\|_{L_x^6} \Big)\Bigg]\\
			&\quad\lesssim{}\langle\tau \rangle^{-\frac32+\frac\delta2}\widetilde{\mathcal{I}}_{\kappa}^\frac12(t)\widetilde{\mathcal{J}}_{\kappa-2}^\frac12(t),
		\end{aligned}
	\end{equation}
	Thus,
	\begin{equation}\label{equ_N2_dd}
		\begin{aligned}
			B_5(u^{df},u^{df},u^{df})\lesssim{}&\sum_{\substack{|\alpha|\le \kappa-3,\\|\beta|+|\gamma|\leq\kappa-3}} \lesssim{} \int_0^t \langle\tau\rangle^{-\frac32+\frac{\delta}{2}} \widetilde{\mathcal J}_{\kappa-2}(\tau) \widetilde{\mathcal I}_{\kappa}^{\frac12}(\tau) \,d\tau .
		\end{aligned}
	\end{equation}
	Next, by the symmetry of $B_{2}$ in its last two arguments,
	$B_{2}(u,u^{cf},u^{df})$ and $B_{2}(u,u^{df},u^{cf})$
	can be estimated in the same way, so we only consider the former.
	We obtain
	\begin{equation}\label{equ-B21-2}
		\begin{aligned}
			B_{2}(u,u^{cf},u^{df}) \lesssim{}& \sum_{\substack{ |\alpha|\leq\kappa-3,\;|\beta|\leq\kappa-4\\|\gamma|\leq\kappa-6 }} 
			\bigl\| \partial\nabla\widetilde\Gamma^\alpha u \bigr\|_{L_\tau^\infty L_x^2} \bigl\| \langle\tau\rangle^{\frac12} \nabla^2\widetilde\Gamma^\beta u^{cf} \bigr\|_{L_\tau^\infty L_x^2} \bigl\| \langle \tau\rangle \nabla^2\widetilde\Gamma^\gamma u^{df} \bigr\|_{L_\tau^\infty L_x^\infty} \\
			&+ \sum_{\substack{ |\alpha|\leq\kappa-3,\;|\gamma|\leq\kappa-4\\|\beta|\leq\kappa-6 }} 
			\bigl\| \partial\nabla\widetilde\Gamma^\alpha u \bigr\|_{L_\tau^\infty L_x^2} \bigl\| \langle\tau\rangle^{\frac12}\langle r\rangle \nabla^2\widetilde\Gamma^\beta u^{cf} \bigr\|_{L_\tau^\infty L_x^\infty} \bigl\| \langle c_2\tau-r\rangle \nabla^2\widetilde\Gamma^\gamma u^{df} \bigr\|_{L_\tau^\infty L_x^2} \\
			\lesssim{}& \left(\mathcal E^{low}_{\kappa-2}(t)\right)^{\frac32}.
		\end{aligned}
	\end{equation}
	
	For $B_{2}(u^{cf},u^{df},u^{df})$,
	applying the weighted estimates for the $cf$ and $df$ components
	to the respective factors, we obtain
	\begin{equation}\nonumber
		\begin{aligned}
			B_{2}(u^{cf},u^{df},u^{df}) \lesssim{}& \sum_{\substack{ |\alpha|\leq\kappa-3,|\beta|\leq\kappa-4\\
					|\gamma|\leq\kappa-6 }} \int_0^t \langle\tau\rangle^{-\frac32+\frac{\delta}{2}} \bigl\| \langle\tau\rangle^{\frac12-\frac\delta2} \partial\nabla\widetilde\Gamma^\alpha u^{cf} \bigr\|_{L_x^2} \bigl\| \nabla^2\widetilde\Gamma^\beta u^{df} \bigr\|_{L_x^2}
			\bigl\| \langle \tau\rangle \nabla^2\widetilde\Gamma^\gamma u^{df} \bigr\|_{L_x^\infty}\,d\tau \\
			\lesssim{}& \int_0^t \langle\tau\rangle^{-\frac32+\frac{\delta}{2}} \widetilde{\mathcal J}_{\kappa-2}(\tau) \widetilde{\mathcal I}_{\kappa}^{\frac12}(\tau) \,d\tau .
		\end{aligned}
	\end{equation}
	Finally, we consider $B_{2}(u,u^{cf},u^{cf})$.
	The Sobolev estimates yield
	\begin{equation}\label{equ-B21-1}
		\begin{aligned}
			B_{2}(u,u^{cf},u^{cf}) \lesssim{}& \sum_{\substack{ |\alpha|\leq\kappa-3\\|\beta|\leq\kappa-4\\
					|\gamma|\leq\kappa-6 }} \bigl\| \partial\nabla\widetilde\Gamma^\alpha u\bigr\|_{L_\tau^\infty L_x^2} \bigl\| \langle\tau\rangle^{\frac12} \nabla^2\widetilde\Gamma^\beta u^{cf} \bigr\|_{L_\tau^2 L_x^2} \bigl\| \langle\tau\rangle^{\frac12} \nabla^2\widetilde\Gamma^\gamma u^{cf} \bigr\|_{L_\tau^2 L_x^\infty} 
			\lesssim{}\left(\mathcal E^{low}_{\kappa-2}(t)\right)^{\frac32}.
		\end{aligned}
	\end{equation}
	
	We next estimate $B_{3}(u,u,u)$.
	When all three arguments are $df$ components, the different null condition allows us to obtain the same estimate by the argument used in \eqref{equ-B21-2}.
	Similarly, the method used in \eqref{equ-B21-1} also applies to $B_{3}(u^{df},u^{cf},u^{df})$.
	Since we do not have an energy bound for $u^{cf}$ with a first-order time weight, we must use weighted estimates for all three factors to obtain the temporal decay $\langle\tau\rangle^{-\frac32+\frac\delta2}$.
	We therefore estimate $B_{3}(u^{cf},u^{cf},u^{cf})$ separately and obtain
	\begin{equation}\nonumber
		\begin{aligned}
			B_{3}&(u^{cf},u^{cf},u^{cf})\\
			&\lesssim{}\sum_{\substack{ |\alpha|\leq\kappa-3,\;|\gamma|\leq\kappa-4\\ |\beta|\leq\kappa-6 }} 
			\int_0^t \langle\tau\rangle^{-\frac32+\frac{\delta}{2}} \bigl\| \langle\tau\rangle^{\frac12-\frac\delta2} \partial\nabla\widetilde\Gamma^\alpha u^{cf} \bigr\|_{L_x^2} 
			\bigl\| \langle\tau\rangle^{\frac12} \nabla\nabla^{\le 2}\widetilde\Gamma^\beta u^{cf} \bigr\|_{L_x^3}\bigl\| \langle\tau\rangle^{\frac12} \nabla\widetilde\Gamma^\gamma u^{cf} \bigr\|_{L_x^6} d\tau \\
			&\lesssim{} \int_0^t \langle\tau\rangle^{-\frac32+\frac{\delta}{2}} \widetilde{\mathcal J}_{\kappa-2}(\tau) \widetilde{\mathcal I}_{\kappa}^{\frac12}(\tau) \,d\tau .
		\end{aligned}
	\end{equation}
	Next, we estimate $B_{3}(u^{df},u,u^{cf})$ and $B_{3}(u^{cf},u,u^{df})$ together to obtain
	\begin{equation}\nonumber
		\begin{aligned}
			B_{3}(u,&u^{df},u^{cf}) +B_{3}(u,u^{cf},u^{df}) \\
			\lesssim{}&\sum_{\substack{ |\alpha|\leq\kappa-3,|\gamma|\leq\kappa-4\\ |\beta|\leq\kappa-6}} 
			\int_0^t \langle\tau\rangle^{-\frac32+\frac\delta2} \bigl\|\partial\nabla\widetilde\Gamma^\alpha u \bigr\|_{L_x^2} \\
			&\times \Biggl[ \bigl\| \langle\tau\rangle^{1-\frac\delta2}\nabla\nabla^{\le 2}\widetilde\Gamma^\beta u^{df} \bigr\|_{L_x^\infty} \bigl\| \langle \tau\rangle^{\frac12}\nabla\widetilde\Gamma^\gamma u^{cf} \bigr\|_{L_x^2}+ \bigl\| \langle\tau\rangle^{1-\frac\delta2} \nabla\widetilde\Gamma^\gamma u^{df} \bigr\|_{L_x^\infty} \bigl\| \langle \tau\rangle^\frac12 \nabla\nabla^{\le 2}\widetilde\Gamma^\beta u^{cf} \bigr\|_{L_x^2} \Biggr]\,d\tau \\
			\lesssim{}& \int_0^t \langle\tau\rangle^{-\frac32+\frac{\delta}{2}} \widetilde{\mathcal J}_{\kappa-2}(\tau) \widetilde{\mathcal I}_{\kappa}^{\frac12}(\tau) \,d\tau .
		\end{aligned}
	\end{equation}
	For the terms $B_{3}(u^{cf},u^{df},u^{df})$ and $B_{3}(u^{df},u^{cf},u^{cf})$, we obtain
	\begin{equation}\nonumber
		\begin{aligned}
			B_{3}&(u^{cf},u^{df},u^{df})+B_{3}(u^{df},u^{cf},u^{cf})\\
			&\lesssim{}\sum_{\substack{ |\alpha|\leq\kappa-3, |\gamma|\leq\kappa-4\\|\beta|\leq\kappa-7 }} 
			\int_0^t \langle\tau\rangle^{-\frac32+\frac\delta2} \bigl\| \nabla\widetilde\Gamma^\gamma  u\bigr\|_{L_x^2}  \bigl\| \langle\tau\rangle \nabla\nabla^{\le 2}\widetilde\Gamma^\beta u^{df} \bigr\|_{L_x^\infty} \bigl\| \langle\tau\rangle^{\frac12-\frac\delta2} \partial\nabla\widetilde\Gamma^\alpha u^{cf} \bigr\|_{L_x^2}\\
			&\qquad+\sum_{\substack{ |\alpha|\leq\kappa-3,|\gamma|\leq\kappa-4\\|\beta|\leq\kappa-6 }}
			\int_0^t \langle\tau\rangle^{-\frac32+\frac\delta2} \bigl\| \nabla\widetilde\Gamma^\gamma  u\bigr\|_{L_x^2} \left\| \langle \tau \rangle^{-\frac\delta2 }\langle c_2 \tau -r\rangle \partial \nabla \widetilde{\Gamma}^\alpha u^{df}\right\|_{L_x^6}\left\| \langle r\rangle \langle \tau \rangle^\frac12 \nabla\nabla^{\le 2}\widetilde{\Gamma}^\beta u^{cf}\right\|_{L_x^2} \\
			&\qquad+\sum_{\substack{ |\alpha|\leq\kappa-3, |\gamma|\le 3 \\|\beta|\leq\kappa-6 }}
			\int_0^t \langle\tau\rangle^{-\frac32+\frac\delta2} \bigl\| \langle\tau\rangle^{\frac12-\frac\delta2} \partial\nabla\widetilde\Gamma^\alpha u^{cf} \bigr\|_{L_x^2}\\
			&\qquad\quad\times\bigg(\bigl\|\langle r\rangle  \nabla\widetilde\Gamma^\gamma  u\bigr\|_{L_x^\infty}  \bigl\|\nabla\nabla^{\le 2}\widetilde\Gamma^\beta u^{df} \bigr\|_{L_x^2} +\bigl\|\nabla\widetilde\Gamma^\gamma  u\bigr\|_{L_x^3}  \bigl\|\langle c_2\tau -r\rangle\nabla\nabla^{\le 2}\widetilde\Gamma^\beta u^{df} \bigr\|_{L_x^6} \bigg)\\
			&\lesssim{}\int_0^t \langle\tau\rangle^{-\frac32+\frac{\delta}{2}} \widetilde{\mathcal J}_{\kappa-2}(\tau) \widetilde{\mathcal I}_{\kappa}^{\frac12}(\tau) d\tau .
		\end{aligned}
	\end{equation}
	
	The term $B_{4}(u,u,u)$ is estimated in the same way as $B_{1}(u,u,u)$, yielding the bound in \eqref{equ_B1}.
	Combining the estimates above, we obtain
	\begin{equation}
		\begin{aligned}
			\sum_{\substack{|\beta|,|\gamma|<|\alpha|\leq\kappa-3,\\\beta|+|\gamma|\leq|\alpha|, |a|\leq1}} &\int_0^t\!\int_{\mathbb R^3} \partial_\tau\nabla^a\widetilde\Gamma^\alpha u_i g^{ijk}_{lmn} \nabla^a\partial_l\!\left( \partial_m\widetilde\Gamma^\beta u_j\, \partial_n\widetilde\Gamma^\gamma u_k \right)dxd\tau\\
			\lesssim{}& \left(\mathcal E^{low}_{\kappa-2}(t)\right)^{\frac32} + \int_0^t \langle\tau\rangle^{-\frac32+\frac{\delta}{2}} \widetilde{\mathcal J}_{\kappa-2}(\tau) \widetilde{\mathcal I}_{\kappa}^{\frac12}(\tau) \,d\tau .
		\end{aligned}
	\end{equation}
	
	It remains to estimate the last term involving higher-order derivatives.
	\begin{equation}\label{equ_4.97}
		\begin{aligned}
			\sum_{\substack{ |\alpha|\leq\kappa-3\\|a|\leq1 }} &
			\int_0^t\int_{\mathbb R^3} g^{ij}_{lm}\, \partial_\tau\nabla^a \widetilde\Gamma^\alpha u_i\, \nabla^a\partial_l \bigl(\partial_m u_j\,\Gamma^\alpha\Theta\bigr) \,dx\,d\tau 
			\lesssim{} \sum_{|\alpha|\le \kappa-3} \left\| \partial\nabla^{\le 1}\widetilde{\Gamma}^\alpha u\right\|_{L_\tau ^\infty L_x^2}\\
			&\quad\times\Big( \left\|\nabla\nabla^{\le1}\Gamma^\alpha\Theta\right\|_{L_\tau^2 L_x^2}
			\left\| \langle\tau\rangle^{\frac12} \nabla\nabla^{\le 1}u^{cf} \right\|_{L_\tau^2 L_x^\infty} 
			+\left\|\Gamma^\alpha\Theta\right\|_{L_\tau^2 L_x^6} \left\| \langle\tau\rangle^{\frac12} \nabla^2\nabla^{\le 1}u^{cf} \right\|_{L_\tau^2 L_x^3}\Big)\\
			&+\sum_{|\alpha|\le \kappa-3}\int_0^t\langle \tau \rangle^{-\frac32+\frac\delta2}\left\| \partial\nabla \widetilde{\Gamma}^\alpha u\right\|_{L_x^2}\left\| \langle \tau \rangle ^{\frac12 -\frac\delta 2}\nabla^{\le 2}\Gamma^\alpha \Theta \right\|_{L_x^2}\left\|\langle \tau \rangle \nabla \nabla^{\le 2}u ^{df}\right\|_{L_x^\infty}   \\
			\lesssim{}& \int_{0}^t\langle \tau\rangle ^{-\frac32+\frac\delta2}\widetilde{\mathcal J}_{\kappa-2}(\tau)\widetilde{\mathcal I}_{\kappa}^{\frac12}(\tau)d\tau + \left(\mathcal E^{low}_{\kappa-2}(t)\right)^{\frac32}.
		\end{aligned}
	\end{equation}
	
	Since $\Theta$ and $u^{cf}$ have similar structures, the same argument, together with \eqref{equ_4.97}, yields
	\begin{equation}\label{equ_B1}
		\begin{aligned}
			\sum_{\substack{|\alpha|\le \kappa-3,|a|\le1\\ \beta+\gamma=\alpha}}&\int_0^t\int_{\mathbb R^3} \partial_\tau\nabla^a\widetilde\Gamma^\alpha u_i\biggl[ g^{ijk}_{lmn}\nabla^a\partial_l \bigl( \partial_m\widetilde\Gamma^\beta u_j\, \partial_n\widetilde\Gamma^\gamma u_k \bigr) +g^{ij}_{lm}\nabla^a\partial_l \bigl( \partial_m\widetilde\Gamma^\beta u_j\, \Gamma^\gamma\Theta \bigr) \biggr]\,dx\,d\tau \\    
			\lesssim{}& \int_{0}^t\langle \tau\rangle ^{-\frac32+\frac\delta2}\widetilde{\mathcal J}_{\kappa-2}(\tau)\widetilde{\mathcal I}_{\kappa}^{\frac12}(\tau)d\tau + \left(\mathcal E^{low}_{\kappa-2}(t)\right)^{\frac32}.
		\end{aligned}
	\end{equation}
	
	For the quadratic nonlinear terms in the elastic equation involving only the temperature, H\"older's inequality yields
	\begin{equation}\label{equ_low_8.51}
		\begin{aligned}
			&\sum_{\substack{|\alpha|\leq\kappa-3,\ |a|\leq1\\|\beta|+|\gamma|\leq|\alpha|}}
			\int_0^t\int_{\mathbb R^3}\left| \partial_\tau\nabla^a\widetilde\Gamma^\alpha u \nabla^{a} \bigl(\nabla\Gamma^\beta\Theta\Gamma^\gamma\Theta\bigr)\right| dxd\tau\\
			&\quad\lesssim \sum_{\substack{|\alpha|\leq\kappa-1\\
					|\gamma|+|\beta|\leq\kappa-1}}\bigl\|\partial_\tau\nabla^{\le 1} \widetilde\Gamma^\alpha u^{cf}\bigr\|_{L_\tau^\infty L_x^2} \bigl\|\nabla^{\le 1}(\nabla\Gamma^\beta\Theta\Gamma^\gamma\Theta)\bigr\|_{L_\tau^1 L_x^2}\\
			&\quad\lesssim \mathcal E^{low}_{\kappa-2}(t)
			\Big( \sum_{\substack{|\beta|\le \kappa-3\\ |\gamma|\le \kappa-6}}\left\|\nabla\nabla^{\le 1} \Gamma^\beta \Theta \right\| _{L_\tau^2 L_x^2}\left\| \nabla^{\le 1}\Gamma^\gamma \Theta \right\| _{L_\tau^2 L_x^\infty}
			+\sum_{\substack{|\beta|+1\le \kappa-6\\ |\gamma|\le \kappa-3}}\left\| \nabla\nabla^{\le 1}\Gamma^\beta \Theta\right\| _{L_\tau^2 L_x^3}\left\| \nabla^{\le 1}\Gamma^\gamma \Theta \right\| _{L_\tau^2 L_x^6}\Big)\\
			&\quad\lesssim\left(\mathcal E^{low}_{\kappa-2}(t)\right)^{\frac32}.
		\end{aligned}
	\end{equation}
	
	Finally, we consider the $M_1$ term in the heat equation.
	Since $\widetilde\Gamma^\gamma u^{cf}$ already carries two derivatives, directly expanding by the product rule with $|a_1|+|a_2|=|a|$ and then applying Sobolev embedding may require more derivatives than the lower-order energy controls.
	We therefore treat $|a|=0$ and $|a|=1$ separately, obtaining
	\begin{equation}
		\begin{aligned}
			&\sum_{\substack{ |\alpha|\leq\kappa-3\\|\beta|+|\gamma|\leq\kappa-3\\|a|\leq1 }}
			\int_0^t\int_{\mathbb R^3} \nabla^a\Gamma^\alpha\Theta\, \nabla^a\bigl( \Gamma^\beta\Theta\, \partial_\tau\nabla\cdot \widetilde{\Gamma}^{\gamma}u^{cf} \bigr) \,dx\,d\tau \\
			&\quad\lesssim{} \sum_{\substack{ |\alpha|,|\beta|\leq\kappa-3\\|\gamma|+2\leq\kappa-5 }} 
			\left\|\Gamma^\beta\Theta\right\|_{L_\tau^\infty L_x^2}\Bigl( \left\|\nabla^{2}\Gamma^\alpha\Theta\right\|_{L_\tau^2 L_x^2}\left\|\langle \tau\rangle^\frac12 \partial\nabla\widetilde{\Gamma}^{\gamma}u^{cf} \right\|_{L_\tau^2 L_x^\infty} + \left\|\Gamma^\alpha\Theta\right\|_{L_\tau^2 L_x^6}\left\|\langle \tau\rangle^\frac12 \partial\nabla\widetilde{\Gamma}^{\gamma}u^{cf} \right\|_{L_\tau^2 L_x^3} \Bigr) \\
			&\qquad+ \sum_{\substack{ |\alpha|,|\gamma|\leq\kappa-3\\|\beta|\leq\kappa-5 }} 
			\left\| \partial_\tau\nabla^{\le1}\widetilde{\Gamma}^{\gamma}u^{cf} \right\|_{L_\tau^\infty L_x^2}\Bigg( \left\|\nabla\nabla^{\le1}\Gamma^\alpha\Theta\right\|_{L_\tau^2 L_x^2} \left\|\Gamma^\beta\Theta\right\|_{L_\tau^2 L_x^\infty}+ \left\|\Gamma^\alpha\Theta\right\|_{L_\tau^2 L_x^6} \left\|\nabla\Gamma^\beta\Theta\right\|_{L_\tau^2 L_x^3}\Bigg) \\
			&\quad\lesssim{} \left(\mathcal E^{low}_{\kappa-2}(t)\right)^{\frac32}.
		\end{aligned}
	\end{equation}
	
	For the last nonlinear term $M_2$ in the heat equation, we use \eqref{eq:F1-alpha-theta-estimate} to obtain
	\begin{equation}\label{eq:J-low-term-10}
		\begin{aligned}
			\sum_{\substack{ |\alpha|\leq\kappa-3,|a|\le 1\\
					|\beta|+|\gamma|\leq\kappa-3 }}& \int_0^t\!\int_{\mathbb R^3} \nabla^{1+a}\Gamma^\alpha\Theta\cdot \nabla^a\bigl( \nabla\Gamma^\beta\Theta\, \Gamma^\gamma\Theta \bigr)\,dx\,d\tau 
			\lesssim{}\left(\mathcal E^{low}_{\kappa-2}(t)\right)^{\frac32}.
		\end{aligned}
	\end{equation}
	
	Substituting the estimates in \eqref{eq:J-cn-initial-cubic}, \eqref{equ_B1}-\eqref{eq:J-low-term-10} into \eqref{eq:J-low-expanded} and applying the quadratic inequality, we obtain
	
	\begin{equation}
		\begin{aligned}
			\mathcal J_{\kappa-2}(t) \lesssim{}& \int_0^t \langle\tau\rangle^{-\frac32+\frac{\delta}{2}}\widetilde{\mathcal J}_{\kappa-2}(\tau)\, \widetilde{\mathcal I}_{\kappa}^{\frac12}(\tau) \,d\tau  + \left(\mathcal E^{low}_{\kappa-2}(t)\right)^{\frac32}+ \mathcal J_{\kappa-2}(0) +\mathcal J^{\frac12}_{\kappa-2}(t)\mathcal J^{\frac12}_{\kappa-3}(t).
		\end{aligned}
	\end{equation}
	
	\subsection{Conclusion of the proof}
	
	Combining Lemmas \ref{lem_J_1}, \ref{lem_J_2},
	\ref{lem_i_1}, and \ref{lem_i_2} with the higher- and
	lower-order energy estimates above, we obtain
	\begin{equation}\nonumber
		\left\{
		\begin{aligned}
			&\mathcal E_{\kappa}(t) \lesssim\langle t\rangle^{-\delta}\int_0^t \langle\tau\rangle^{-1+\delta} \mathcal E_{\kappa}(\tau) \left(\mathcal E^{low}_{\kappa-2}(\tau)\right)^{\frac12} d\tau+\mathcal E_\kappa(t) \left(\mathcal E^{low}_{\kappa-2}(t)\right)^{\frac12} +\mathcal I_{\kappa-1}(t) +\widetilde{\mathcal I}_{\kappa}(0), \\
			&\mathcal E^{low}_{\kappa-2}(t) \lesssim\int_0^t \langle\tau\rangle^{-\frac32+\frac{\delta}{2}} \mathcal E^{low}_{\kappa-2}(\tau) \mathcal E_{\kappa}^{\frac12}(\tau) d\tau+\left(\mathcal E^{low}_{\kappa-2}(t)\right)^{\frac32} +\mathcal J_{\kappa-3}(t) +\mathcal J_{\kappa-2}(0),
		\end{aligned}
		\right.
	\end{equation}
	
	We iterate the two inequalities separately, reducing the energy orders
	from $\kappa$ and $\kappa-2$ to $1$.
	At the lowest order, no vector fields are applied to the system,
	so no linear remainder terms arise from commutators.
	We therefore obtain
	\begin{equation}
		\left\{
		\begin{aligned}
			&\mathcal E_{\kappa}(t) \leq{} C\int_0^t \langle\tau\rangle^{-1} \mathcal E_{\kappa}(\tau) \left(\mathcal E^{low}_{\kappa-2}(\tau)\right)^{\frac12} d\tau+C\mathcal E_\kappa(t) \left(\mathcal E^{low}_{\kappa-2}(t)\right)^{\frac12} +C\widetilde{\mathcal I}_{\kappa}(0), \\
			&\mathcal E^{low}_{\kappa-2}(t) \leq{} C\int_0^t \langle\tau\rangle^{-\frac32+\frac{\delta}{2}} \mathcal E^{low}_{\kappa-2}(\tau) \mathcal E_{\kappa}^{\frac12}(\tau) d\tau+C\left(\mathcal E^{low}_{\kappa-2}(t)\right)^{\frac32} +C\widetilde{\mathcal J}_{\kappa-2}(0),\\
			&\widetilde{\mathcal I}_{\kappa}(0)\le C\widetilde E_{\kappa}(0),\\
			&\widetilde{\mathcal J}_{\kappa-2}(0)\le C\widetilde E_{\kappa-2}(0).
		\end{aligned}
		\right.
	\end{equation}
	Here, $C>2$ is a uniform constant.
	
	Assume that the initial data satisfy
	\begin{equation}
		\label{equ:initial-condition} C^{-1}\widetilde{\mathcal J}_{\kappa-2}(0)\exp\left(C_0 \sqrt{C^{-1}}\widetilde{\mathcal I}_{\kappa}^{\frac12}(0)\right) \le \widetilde E_{\kappa-2}(0)\exp\left(C_0 \widetilde E_{\kappa}(0)^{\frac12} \right) \le \varepsilon^2.
	\end{equation}
	Let $T>0$ be the maximal time for which the following
	bootstrap assumption holds:
	\begin{equation}\label{equ:bootstrap-assumption}
		\mathcal E^{low}_{\kappa-2}(t) \leq 4C^2\varepsilon^2, \qquad 0\leq t\leq T.
	\end{equation}
	Choose $\varepsilon_0>0$ sufficiently small so that,
	for every $0<\varepsilon\leq\varepsilon_0$,
	\begin{equation}
		4C^3\varepsilon \le \frac12.
	\end{equation}
	Using the bootstrap assumption to absorb the nonlinear terms
	involving the energies at time $t$ into the left-hand sides, we obtain
	\begin{equation}	\label{equ:coupled-energy}
		\left\{
		\begin{aligned}
			&\mathcal E_{\kappa}(t) \leq{}2C\int_0^t \langle\tau\rangle^{-1} \mathcal E_{\kappa}(\tau) \left(\mathcal E^{low}_{\kappa-2}(\tau)\right)^{\frac12} d\tau +2C\widetilde{\mathcal I}_{\kappa}(0), \\
			&\mathcal E^{low}_{\kappa-2}(t) \leq{} 2C\int_0^t \langle\tau\rangle^{-\frac32+\frac{\delta}{2}} \mathcal E^{low}_{\kappa-2}(\tau) \mathcal E_{\kappa}^{\frac12}(\tau) d\tau +2C\widetilde{\mathcal J}_{\kappa-2}(0).
		\end{aligned}
		\right.
	\end{equation}
	Applying the bootstrap assumption and Gronwall's inequality
	to the first inequality in \eqref{equ:coupled-energy}, we obtain
	\begin{equation}
		\label{equ:high-order-growth} \mathcal E_{\kappa}(t) \leq 2C\widetilde{\mathcal I}_{\kappa}(0) \langle t\rangle^{ 4C^2\varepsilon }.
	\end{equation}
	Substituting \eqref{equ:high-order-growth} into the second inequality
	in \eqref{equ:coupled-energy}, we obtain
	\begin{equation}\nonumber
		\mathcal E^{low}_{\kappa-2}(t) \leq{}2\sqrt{2} C^{\frac32} \widetilde{\mathcal I}_{\kappa}^{\frac12}(0) \int_0^t \langle\tau\rangle^{ -\frac32+\frac{\delta}{2} +2C^2\varepsilon } \mathcal E^{low}_{\kappa-2}(\tau) d\tau+2C\widetilde{\mathcal J}_{\kappa-2}(0).
	\end{equation}
	For $0<\delta<\frac14$, we further decrease $\varepsilon_0$ so that $2C^2\varepsilon_0\leq\frac14$.
	The time weight on the right-hand side is then integrable over $[0,\infty)$.
	Applying Gronwall's inequality once more, and using the initial condition \eqref{equ:initial-condition}, we obtain
	\begin{equation}\nonumber
		\mathcal E^{low}_{\kappa-2}(t) \leq 2 C\widetilde{\mathcal J}_{\kappa-2}(0) \exp\left( C_0 \sqrt{C^{-1}}\widetilde{\mathcal I}_{\kappa}^{\frac12}(0) \right)\leq 2C^2\varepsilon^2 < 4C^2\varepsilon^2,
	\end{equation}
	where the time integral is bounded uniformly in $t$. We may take
	\begin{equation}\nonumber
		2\sqrt{2}C^2 \int_0^t \langle\tau\rangle^{ -\frac32+\frac{\delta}{2} +2C^2\varepsilon } d\tau\leq 32C^2\triangleq C_0.
	\end{equation}
	This strictly improves the bootstrap assumption
	\eqref{equ:bootstrap-assumption}.
	By continuity, the above estimates hold for all $t\ge0$.
	
	Finally, we obtain
	\begin{equation}\nonumber
		\widetilde E_{\kappa-2}(t)\leq C_0^2\varepsilon^2, \qquad \widetilde E_{\kappa}(t)\leq C_0^2\widetilde E_{\kappa}(0)\langle t\rangle^{C_0\varepsilon+\delta}.
	\end{equation}
	This completes the bootstrap argument.

\bigskip
\begingroup
\small
\noindent\textsc{School of Mathematical Sciences, Fudan University, Shanghai 200433, China.}\par
\noindent\textit{Email address:} \texttt{25110180056@m.fudan.edu.cn}\par
\medskip
\noindent\textsc{School of Mathematical Sciences, Fudan University, Shanghai 200433, China.}\par
\noindent\textit{Email address:} \texttt{yizhou@fudan.edu.cn}\par
\endgroup
	
\end{document}